\documentclass[a4paper,fleqn]{cas-sc}

\usepackage[authoryear,longnamesfirst]{natbib}

\def\tsc#1{\csdef{#1}{\textsc{\lowercase{#1}}\xspace}}
\tsc{WGM}
\tsc{QE}

\usepackage{graphicx}%
\usepackage{multirow}%
\usepackage{amsmath,amssymb,amsfonts}%
\usepackage{soul, color, xcolor}
\usepackage{textcomp}%
\usepackage{booktabs}%
\usepackage{amsthm}%
\usepackage{mathrsfs}%
\usepackage[title]{appendix}%
\usepackage{manyfoot}%
\usepackage{algorithm}%
\usepackage{algorithmicx}%
\usepackage{algpseudocode}%

\usepackage{listings}%
\usepackage[justification=centering]{caption}
\usepackage{caption}

\usepackage[utf8]{inputenc} 
\usepackage[T1]{fontenc}    
\usepackage{hyperref}       
\usepackage{url}            
\usepackage{nicefrac}       
\usepackage{microtype}      
\usepackage{lipsum}
\usepackage{fancyhdr}       
\usepackage{adjustbox}
\usepackage{fancybox}

\usepackage{enumerate}
\usepackage[nolist, nohyperlinks, printonlyused]{acronym} \usepackage{epsfig} 
\usepackage{epstopdf}
\usepackage{subfigure}
\usepackage{verbatim}
\usepackage{array}
\usepackage{caption}
\usepackage{amssymb}
\usepackage{fancyhdr}
\usepackage{enumitem}
\usepackage{tabularx}
\usepackage{bm}
\graphicspath{{media/}}
\usepackage{indentfirst} 
\newcommand{\figlegend}[1]{\footnotesize #1}

\theoremstyle{thmstyleone}%
\newtheorem{theorem}{Theorem}
\theoremstyle{thmstyletwo}%

\theoremstyle{thmstylethree}%

\begin{document}
\let\WriteBookmarks\relax
\def\floatpagepagefraction{1}
\def\textpagefraction{.001}

\shorttitle{}    

\shortauthors{}  

\title [mode = title]{QC-OT: Optimal Transport with Quasiconformal Mapping}  



%

\author[1]{Yuping Lv}



\ead{xiaobu.lv@stu.xjtu.edu.cn}



\affiliation[1]{organization={School of Mathematics and Statistics Xi’an Jiaotong University},
	addressline={}, 
	city={Xi’an},
	postcode={710049}, 
	state={Shaanxi},
	country={China}}

\author[1]{Qi Zhao}\ead{shmilyqi@stu.xjtu.edu.cn}
\author[1]{Xuebin Chang}\ead{xuebin\_chang@stu.xjtu.edu.cn}
\author[1]{Wei Zeng}\ead{wz@xjtu.edu.cn}
\cormark[1] 

\cortext[1]{Corresponding author}




\begin{abstract}
The optimal transport (OT) map offers the most economical way to transfer one probability measure distribution to another. 
Classical OT theory does not involve a discussion of preserving topological connections and orientations in transmission results and processes. 
Existing numerical and geometric methods for computing OT seldom pays specific attention on this aspect. Especially, when dealing with the triangular mesh data, the known semi-discrete geometric OT (sd-OT) method employs critical operation of Delaunay triangulation (DT) to adapt topology to ensure the convexity of the energy function and the existence of the solution. This change in topology hampers the applicability of OT in modeling non-flip physical deformations in real-world tasks such as shape registration and editing problems in computer vision and medical imaging fields.
This work introduces the topology structure-preserving optimal transport (QC-OT) map for the triangular mesh input. 
The computational strategy focuses on the two components: relaxing DT and convexity check in sd-OT and integrating quasiconformal (QC) correction. Here, quasiconformal mapping is employed to correct the regions unexpected distortions,  
and guarantee the topological preserving property of the transport. Furthermore, the spatial-temporal topology-preserving OT map is presented based t-OT to study the dynamics of the transportation. 
Multiple experiments have validated the efficiency and effectiveness of the proposed method and demonstrated its potential in the applications of mesh parameterization and image editing.
\end{abstract}


\maketitle

\section{Introduction}\label{sec:introduction}
With the advancement of deep learning, optimal transport (OT) theory has gained increasing attention and has been widely employed in many applications, such as natural language processing \citep{sun2024spectr,huang2016supervised,kusner2015word}, computer vision \citep{izquierdo2024optimal,bonneel2011displacement,solomon2015convolutional}, generative adversarial network \citep{arjovsky2017wasserstein}, clustering \citep{ho2017multilevel}, domain adaptation \citep{flamary2016optimal}, anomaly detection \citep{tong2022fixing}, and so on. 
Introduced by French mathematician \cite{monge1781memoire}, OT aims to find the most economical way to transfer one probability measure distribution to another. \cite{kantorovich2006problem} proposed the relaxation version of OT method which allows mass splitting from a source towards several targets.
Based on the classical methods mentioned above, various methodologies have emerged and can be categorized as follows.

\vspace{1mm}
\noindent\textbf{Numerical methods.} Since OT maps are highly nonlinear, improving the computational efficiency of optimal transport has become a concern. \cite{peyre2019computational} introduced an alternative numerical method, the well-known Sinkhorn method, for the L2 Monge-Kantorovich problem based on the reinterpretation of the mass transfer problem in a continuum mechanics framework. 
The Sinkhorn method accelerates computation by adding an entropic regularization term to the original transport problem and solving it using a simple alternate minimization scheme.  
However, this speedup comes at the cost of sacrificing accuracy; in addition, it generates an OT plan and has to lose information when an OT map is required. 
This method and its variations deal with discrete points directly without specially considering the topology of the data. 

\begin{figure}[t]
\centering
\includegraphics[width=0.7\textwidth]{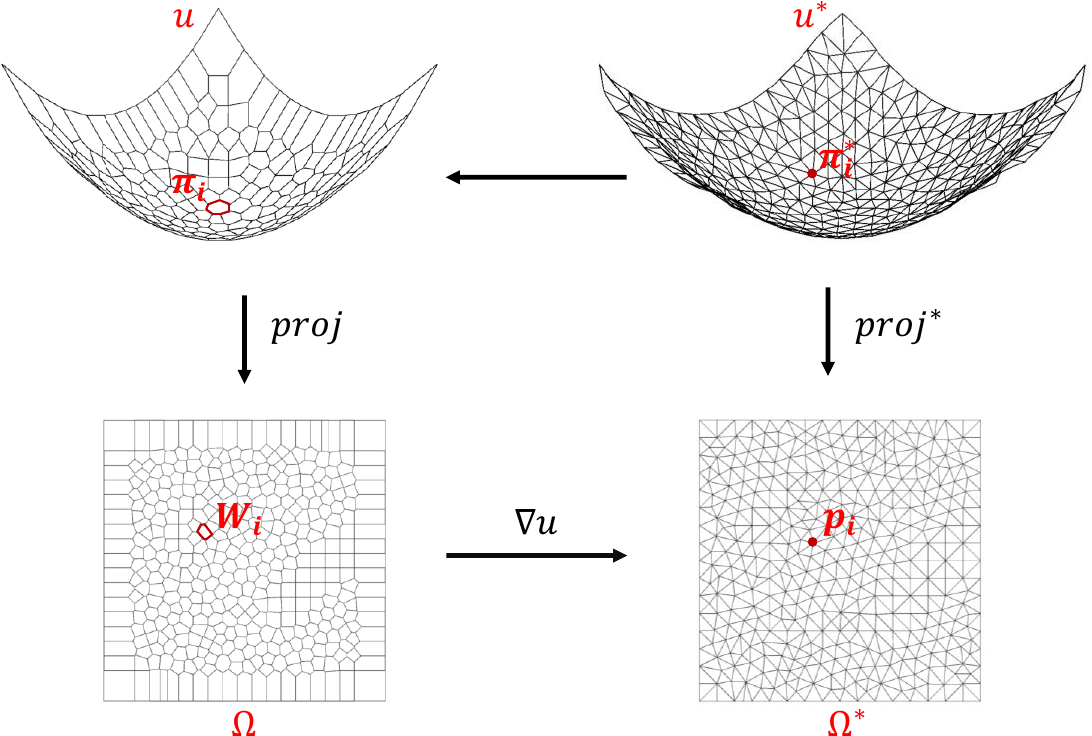} \\
\captionsetup{justification=justified, singlelinecheck=false}
\caption{Geometric semi-discrete OT method \citep{gu2013variational}.
The graph of $u$ and $u^{*}$ are the upper envelope and the lower convex hull of the hyperplane $\pi_{i}$ and its dual point $\pi_{i}^{*}$. 
$u$ denotes the Brenier potential function and the gradient of it $\bigtriangledown u$ is the OT map between domains $\Omega$ and $\Omega^{*}$, which transports the power cell center of $W_{i}$ to $p_{i}$.}
\label{SDOT2}
\end{figure}

\vspace{1mm}
\noindent\textbf{Geometric methods.} In recent years, geometric variational algorithms \citep{gu2013variational,su2015optimal,lei2021optimal} for computing semi-discrete OT (sd-OT) have been presented. These methods rely on the gradient of the Brenier potential function to determine the OT map between the projections of the upper envelope and the lower convex hull (see Fig. \ref{SDOT2}). However, such geometric variational algorithms have stringent geometric restrictions, such as satisfying Delaunay triangulation (DT) and ensuring non-empty power cells to guarantee the convexity of the energy and thus the convergence of algorithm and the existence and uniqueness of the solution. This approach uses the topological information of the input data but inevitably adjust it during the OT calculation. 

\vspace{1mm}
\noindent\textbf{Fluid dynamic methods.} \cite{benamou2000computational} proposed a method to solve the Monge-Ampere equation by minimizing the kinetic energy of a flow field. 
In fluid mechanics, OT maps are obtained through fluid dynamics by processing the vector field, such as eliminating the curl component \citep{haker2004optimal} and imposing specific constraints on the vector field \citep{chen2021spatiotemporal}. However, these approaches may fail in high dimensions.
Moreover, although methodologies in fluid mechanics can illustrate dynamic changes during transportation process, they mainly focus on particle transport and fail to preserve the topological structure \citep{gu2013variational}. 

\vspace{1mm}
It is worth noting that existing methods treat the OT problem as a transportation plan or map between \textit{discrete points}, without taking into account preserving the original topological orientation and connectivity, especially in geometric category of methods. We call this as the \textit{topology-preserving} property in this work. 
However, the dynamics of many natural systems are fundamentally constrained by their underlying structure \citep{stankovic2020data,pan2019deep}, e.g., energy, topology or geometry, thus preserving structure in dynamics is crucial and highly desired in related research and application fields \citep{arroyo2015evaluation,clough2020topological,yu2020c2fnas}. With this inspiration, we aim to investigate the topology-preserving OT map under the geometric variational framework of sd-OT for the triangular mesh input data.

\begin{figure*}
\centering
\footnotesize
\begin{tabular}
{@{\hspace{1em}}m{0.02\linewidth}<{\centering}@{\hspace{0em}}m{0.23\linewidth}<{\centering}@{\hspace{-0.5em}}m{0.23\linewidth}
<{\centering}@{\hspace{1em}}m{0.23\linewidth}
<{\centering}@{\hspace{-0.5em}}m{0.23\linewidth}<{\centering}@{}}
A &
\includegraphics[height=0.23\textwidth]{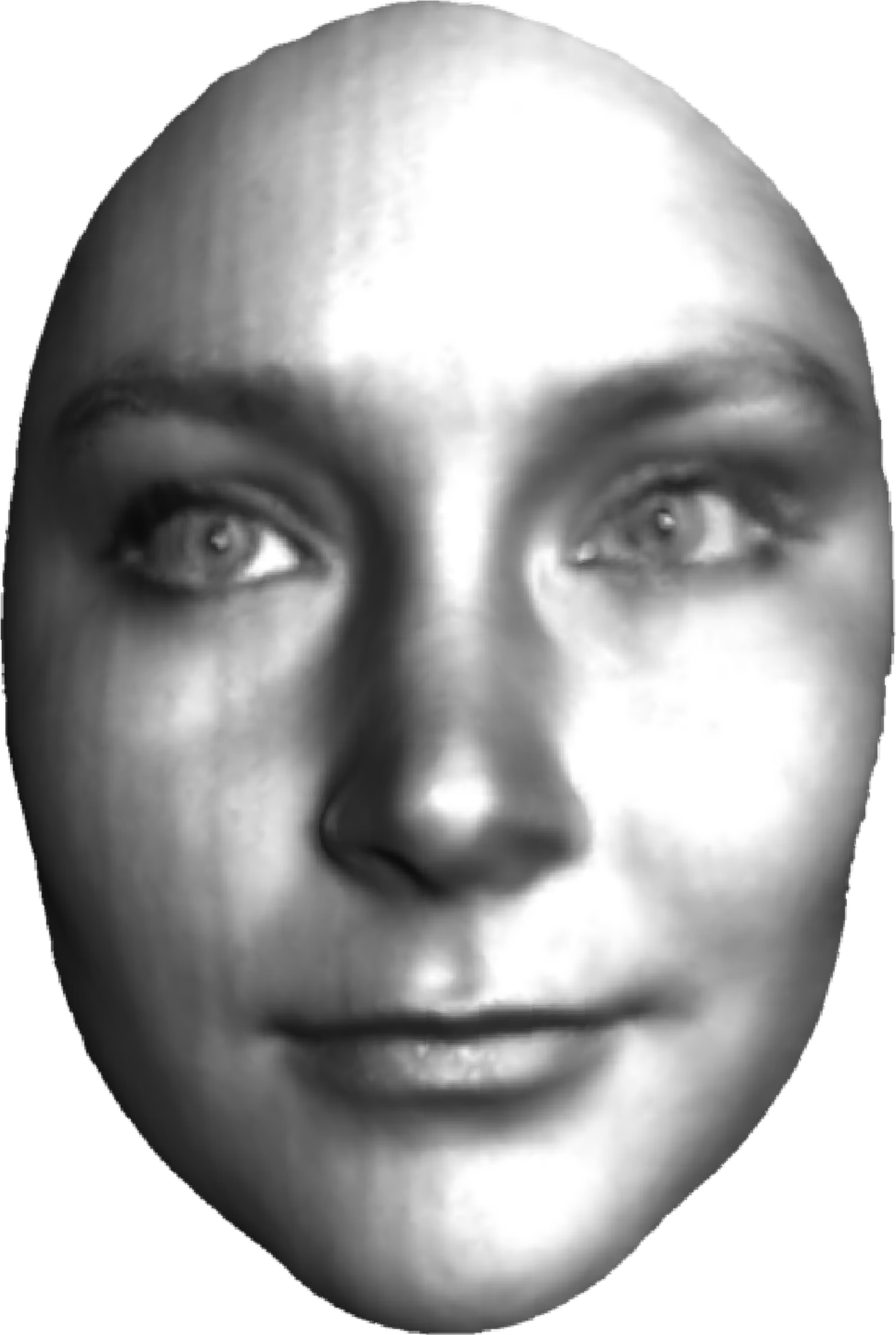} & 
\includegraphics[width=0.23\textwidth]{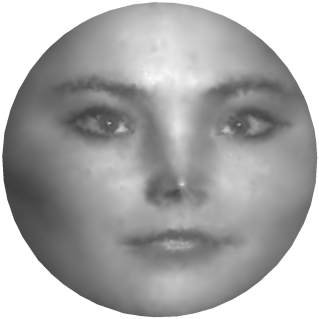} & 
\includegraphics[width=0.23\textwidth]{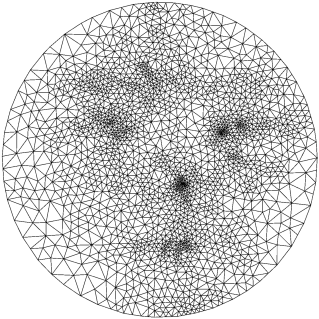} &
\includegraphics[height=0.23\textwidth]{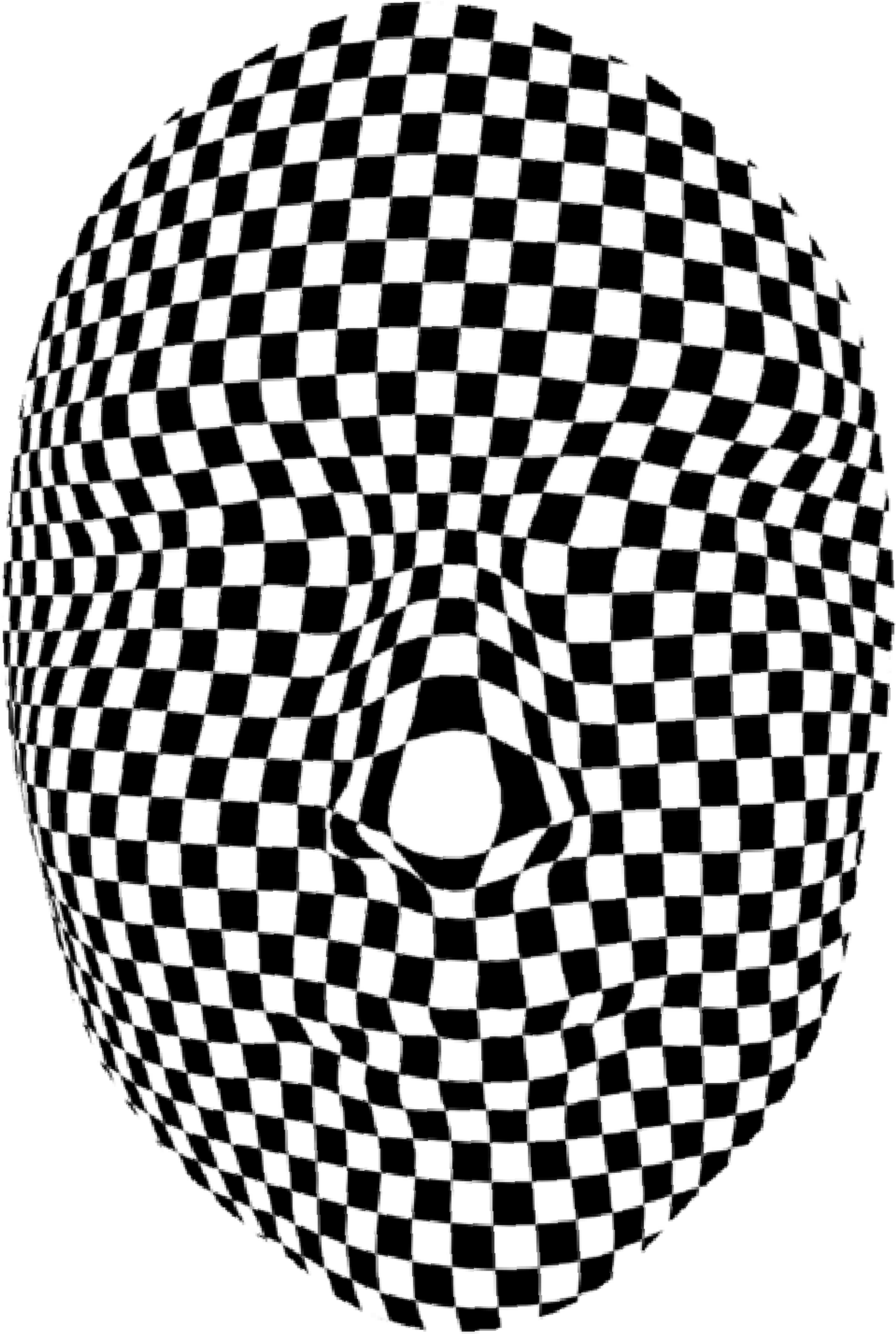} 
\\
B&
\includegraphics[height=0.23\textwidth]{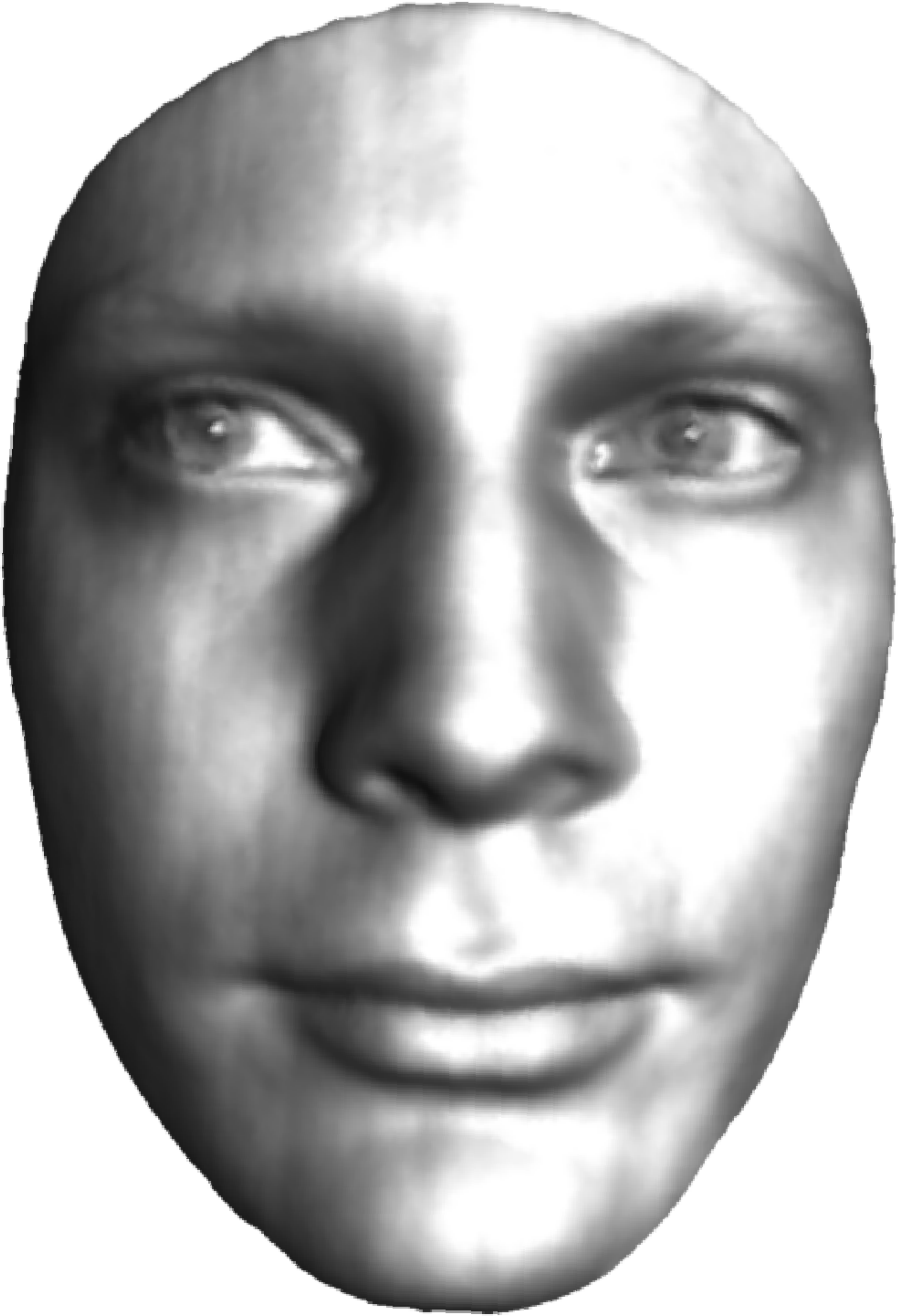}&
\includegraphics[width=0.23\textwidth]{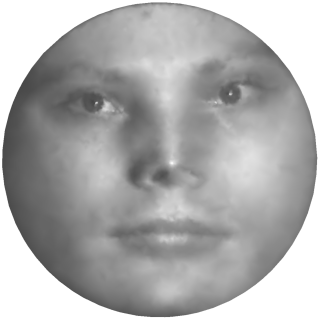}&
\includegraphics[width=0.23\textwidth]{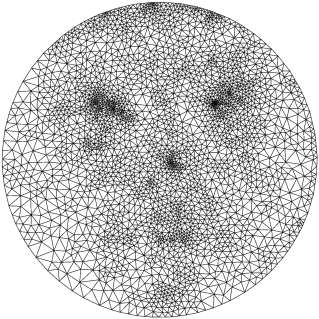}&
\includegraphics[height=0.23\textwidth]{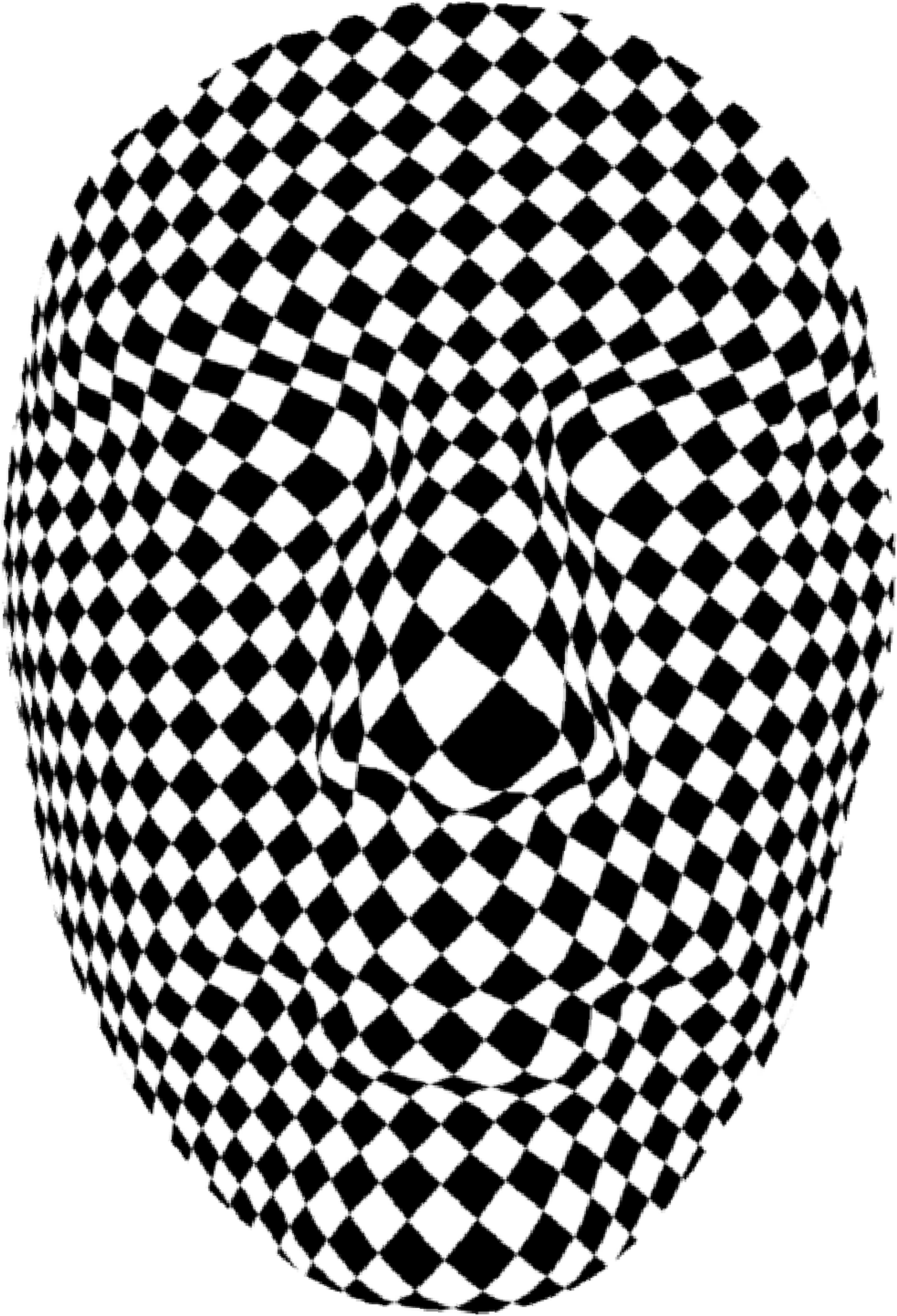}
\\
C&
\includegraphics[height=0.23\textwidth]{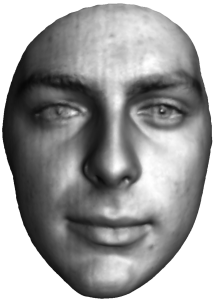}&
\includegraphics[width=0.23\textwidth]{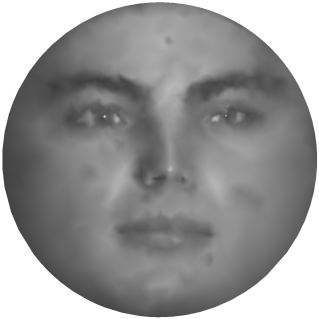}&
\includegraphics[width=0.23\textwidth]{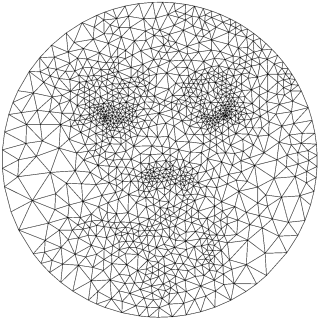}&
\includegraphics[height=0.23\textwidth]{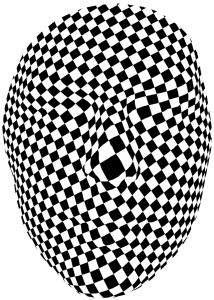}
\\
\figlegend{&(a) 3D Surface & (b) Parameter Domain & (c) Parameter Mesh & (d) Texture Mapping} \\

\end{tabular}
\caption{Conformal parameterization of 3D surfaces A, B and C. 
 }
\label{fig:faces-conformal}
\end{figure*}

\begin{figure}[h]
\centering
\footnotesize
\begin{tabular}{cccc}
\includegraphics[height=0.23\textwidth]{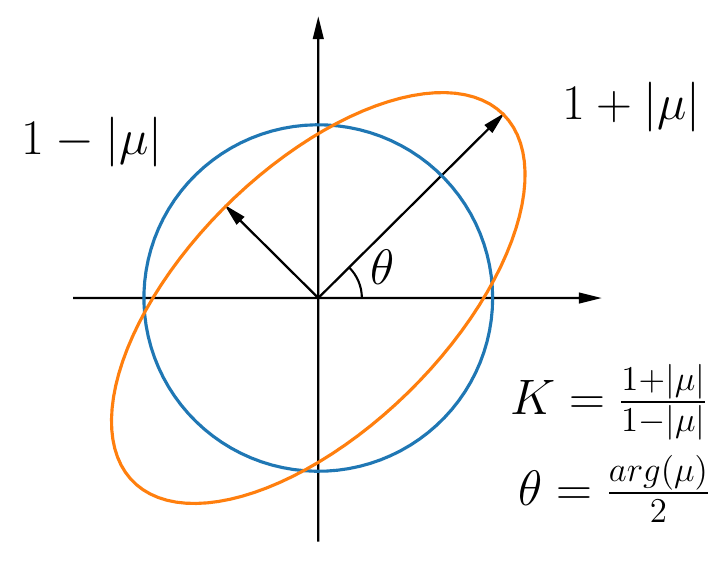}&\includegraphics[height=0.23\textwidth]{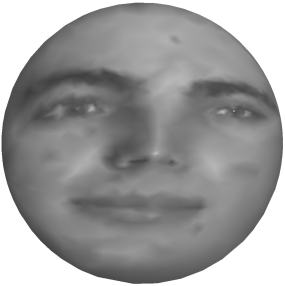}&
 \includegraphics[height=0.23\textwidth]{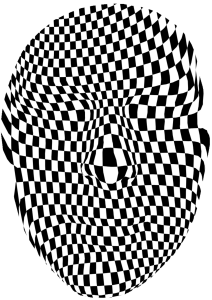}&
  \includegraphics[height=0.23\textwidth]{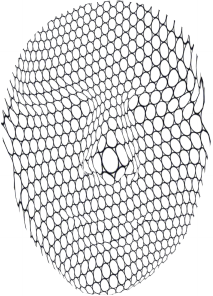}\\
 (a) Beltrami Coefficient& (b) $\mu=0.8 + 0.5i$ & (c) Checker-Board & (d) Circle-Packing 
\end{tabular}
\captionsetup{justification=justified, singlelinecheck=false}
\caption{Quasiconformal mapping associated with a Beltrami coefficient function $\mu$. (a) The QC mapping maps infinitesimal circles to infinitesimal ellipses with the maximal dilation $K$. (b) The QC mapping of surface C in Fig. \ref{fig:faces-conformal} shows obvious angle distortions from its conformal mapping. (c-d) The angle distortions are visualized in both checker-board and circle-packing texture mapping results. }
\label{fig:QC-mu}
\end{figure}

We integrate differential geometric mapping into sd-OT framework by leveraging the fact that geometric mapping generates diffeomorphism and preserves mesh topology , e.g., the well-known conformal mapping for 3D surfaces (see Fig. \ref{fig:faces-conformal}). 
Surface parameterization is expected to be a one-to-one and onto geometric mapping from the surface to a suitable parameter domain, and is fundamental in geometric modeling, computer-aided design, computer graphics and vision \citep{floater2005surface,levy2006parameterization}. Typical geometric mappings are tailored to have specific properties, usually including angle-preserving (conformal) \citep{mullen2008spectral}, area-preserving \citep{choi2022adaptive}, distance-preserving  \citep{estellers2012efficient}, and so on. As a unified form of general mappings, quasiconformal (QC) mapping is a generalization of conformal mapping that describes angle distortions of a mapping by a complex-valued function, Beltrami coefficient (BC) $\mu$ $(0\leq\| \mu \| _{\infty}<1)$ (see Fig. \ref{fig:QC-mu}). 
When $\mu$ equals zero everywhere, the mapping becomes conformal. When $\mu$ satisfies $\| \mu \| _{\infty} > 1$ that is equivalent to the Jacobian determinant $\| J\| < 0$,  there must be self-flip(s), i.e., topology is not preserved. 
We can use the metric $\| \mu \| _{\infty} > \varepsilon$ with a threshold $0\leq \varepsilon \leq 1$ to control the extent of deformations where the angle distortions exceeds the specified range. This would help us efficiently discover the areas with flips or large twists.

Computing a QC mapping associated with a given Beltrami coefficient function $\mu$ is equivalent to solving the Beltrami equation. For a simply connected domain, each mapping has a unique Beltrami coefficient up to a M\"obius transformation, and vice versa. The existence and uniqueness of the solution and the diffeomorphism property are guaranteed in theory \citep{ahlfors2006lectures}. 
The existing computational methods for QC mapping include Beltrami holomorphic flow \citep{lui2012optimization}, linear Beltrami solver \citep{lui2013texture}, holomorphic 1-form \citep{zeng2009surface}, and Ricci flow \citep{zeng2012computing}. Among them, an efficient computing strategy is to apply the auxiliary metric to convert a QC mapping to be a conformal mapping, therefore, any conformal mapping methods can be employed to compute QC mapping with the new metric \citep{zeng2009surface}. QC mapping is capable of representing a general mapping associated with various BC. 
It has been successfully used in medical diagnosis \citep{wong2023quasiconformal}, surface parameterization \citep{lyu2023bijective}, retinotopic map \citep{ta2022quantitative,tu2020diffeomorphic}, decorated surface parameterization \citep{zeng2014surface}, and surface registration \citep{zeng2011registration,zeng2014surface2,lam2015landmark,ma2017robust}. 

In this work, the concept of topology-preserving optimal transport (t-OT) is introduced. 
We propose the relaxed version of semi-discrete OT by removing DT operations and convexity checks to keep the topology during OT and then append QC mapping to it to achieve the t-OT map, which leverages the diffeomorphism property of QC mapping and the geometric sd-OT algorithmic framework. Because of the relaxation of convexity guarantee in sd-OT, unexpected distortions (flips, degenerated triangles, user-customized conditions) may occur 
which are then corrected by QC mapping \citep{ahlfors2006lectures} in our framework. The obtained t-OT map for a triangular mesh input preserves topological relationship and has no self-flips. Such a mapping is also referred to as an orientation-preserving mapping \citep{liu2024one}. 

Besides the property of optimal  transport cost, the t-OT method has flexibility of control on mapping distortions, through setting the angle distortion level by the threshold $\epsilon$ and prescribing the area distortion by specifying a more intuitive measure (density) function. It provides a comprehensive tool for generating geometric mappings. Experiments were conducted to evaluate the efficiency of the proposed algorithmic framework and demonstrate the potential as a practical approach to dealing with geometry processing and medical imaging applications. 

In summary, the contributions of this work include:

\begin{enumerate}

 \item \textit{Topology-Preserving OT (t-OT)}: We present a relaxed version of semi-discrete OT that preserves the original topology by integrating quasiconformal mapping. It is the first work that maintains topological connectivity and orientation for the triangular mesh inputs, overcoming the limitations of traditional OT methods.

 \item \textit{Spatial-Temporal Topology-Preserving OT (tt-OT)}: We provide a temporal model of the t-OT to dynamically and continuously display the topology-preserving transmission process. This facilitates the simulation of physical evolution.
 
 \item \textit{Measure-Driven OT Variants}: 
We perform a series of design of probability measure functions, such as scalar, intensity and area, to drive the calculation of OT maps. This provides users with flexible mapping customization capabilities.
  
 \item \textit{Applicability}: Our proposed technique has demonstrated its efficiency and potential for general geometric processing and medical imaging applications such as mesh parameterization and image-driven editing. 
\end{enumerate}

The rest of the paper is organized as follows: 
Section \ref{sec:background} reviews the theoretical background of closely related geometric mapping and OT theory. Section \ref{sec:algorithm} addresses the proposed algorithms of t-OT and Section \ref{sec:tt-OT} further presents a temporal model based on the t-OT framework. Section \ref{sec: Algorithm Analysis} performs the algorithm analysis through example illustrations. Section \ref{sec:experiment} demonstrates various potential applications of t-OT on 3D surface parameterization and image-driven editing tasks. Finally, the work was concluded in Section \ref{sec:conclusion} along with future research directions outlined.

\section{Theoretical Background}
\label{sec:background}
This work proposes a method for computing topology-preserving OT map for simply connected surfaces. It employs the sd-OT map by \cite{gu2013variational} as a baseline and utilizes the diffeomorphism property of quasiconformal geometric mapping \citep{zeng2012computing,Yau2016ComputationalCG}. Here, we briefly review the relevant prerequisites.

\subsection{Optimal Transport Theory}
Let $\Omega$ and $\Omega^{*}$ be two compact convex domains in Euclidean space $\mathbb{R}^{d}$, with probability measures $\tau$ and $\nu$, respectively, satisfying the equal mass condition $\tau(\Omega) = \nu(\Omega^{*})$. 
The corresponding density functions $f, g$ are given by $d\tau = f(x)dx$ and $d\nu = g(y)dy$. 
The OT map $T: \Omega \to \Omega^{*}$ is a measure-preserving mapping with the minimization of the transportation cost,  
\begin{equation*}
\min \int_{\Omega} c(x,T(x))d\tau(x), 
s.t., T_{\#}\tau=\nu , 
\end{equation*}
where $ c : \Omega \times \Omega^{*} \rightarrow \mathbb{R}^{+}$ is the cost function, and $T_{\#}\tau = \nu$ (symbol $_{\#}$ denotes the push forward operator) means that for any Borel set $B \subset \Omega^{*}$, $\int_{T^{-1}(B)}d\tau(x)=\int_{B}d\nu(y)$, which gives the measure preserving constraint.

\cite{brenier1991polar} solved the OT map using the Monge-Ampere equation 
\begin{equation*}
	\det D^{2}u(x) = \frac{f(x)}{g\circ u(x)}, \text{ s.t. } \nabla u(\Omega) = \Omega^{*},
\end{equation*}
and stated that with the cost function $c(x,y) = \frac{1}{2} \left|  x-y \right|^{2} $, the OT map is the gradient of the Brenier potential function $u:\Omega \rightarrow \mathbb{R}$, given by $T = \nabla u$. Here, $u \in C^{2}$ (twice continuously differentiable). The operator $\det D^{2}u(x)$ represents the determinant of the Hessian matrix of $u$. In theory, the OT map exists and is unique. 

\vspace{1mm}
\noindent\textbf{Semi-Discrete Optimal Transport.}
\cite{gu2013variational} proposed variational principles for discrete Monge-Ampere equation and gave a more general geometric variational approach to the semi-discrete OT problem, i.e., the original triangular mesh is organized by a set of discrete points with prescribed probability measure defined on each point and the density function on the reference domain $\Omega$ is smooth. Here the probability measure in the reference domain is set to be uniform everywhere.

\begin{theorem}
(Alexandrov \citep{gu2013variational}) Suppose $\Omega$ is a compact convex domain in $\mathbb{R}^{d}$, and $p_{1}, p_{2}, \dots , p_{n}$ are $n$ distinct points with $A_{1}, A_{2}, \dots , A_{n} > 0$ such that $\sum_{i=1}^{n} A_{i} = \text{vol}(\Omega)$. Then there exists a vector $\textbf{h} = \left( h_{1}, h_{2}, \dots , h_{n}\right) \in \mathbb{R}^{n}$, unique up to adding the constant $(c, c, \dots , c)$, such that the piecewise linear convex function
\begin{equation*}
	u_{\textbf{h}}(x)= \max_{x \in \Omega} \left\lbrace \left\langle x, p_{i}\right\rangle + h_{i} \right\rbrace,
\end{equation*}
satisfies $vol(\left\lbrace x \in \Omega | \nabla u(x) = p_{i} \right\rbrace ) = A_{i}$.
\end{theorem}

The functions $u$ and $\nabla u(x)$ in this theorem are called the Alexandrov potential and Alexandrov map, respectively. Here, $u$ is also the upper envelope of hyperplanes $\pi_{i}(x)= \left\lbrace \left\langle x, p_{i}\right\rangle + h_{i} \right\rbrace$, $i=1,\dots, n$ , denoted as $Env(\left\lbrace \pi_{i} \right\rbrace_{i=0}^{n})$.
For each hyperplane $\pi_{i}(x)$, its upper envelope is the graph of the function $u_{\textbf{h}}(x)$. The Legendre dual of $u_{\textbf{h}}(x)$ is defined as 
\begin{equation*}
	u^{*}(y) := \max_{x\in \mathbb{R}^{d}}\left\lbrace \left\langle x,y \right\rangle -u(x) \right\rbrace ,
 \label{dual}
\end{equation*}
which introduces the discrete Brenier potential. 

As illustrated in Fig. \ref{SDOT2}, each hyperplane $\pi_{i}(x)$ has a dual point in $\mathbb{R}^{d+1}$, denoted as $\pi_{i}^{*} := (p_{i},-h_{i})$, and the graph of $u^{*}$ is the lower convex hull of $\left\lbrace \pi_{i}^{*} \right\rbrace _{i=1}^{n}$. The projection of the upper envelope produces a nearest power diagram $\mathcal{T}$ of $\Omega$, while the projection of the lower convex hull induces a nearest weighted Delaunay triangulation $\mathcal{T}^{*}$ of $\Omega^{*}$. Based on the Legendre transformation theory \citep{gu2013variational}, the lower convex hull and the upper envelope are dual to each other, namely vertex $p_{i}$ connects to $p_{j}$ if and only if the power cell $W_{i}(\textbf{h})$ is adjacent to $W_{j}(\textbf{h})$. The power diagram introduced by projecting the upper envelope 
$Env(\left\lbrace \pi_{i} \right\rbrace_{i=1}^{n})$ into $\Omega$ is defined as
\begin{equation*}
	\Omega = \bigcup _{i=1} ^{n} W_{i}(\textbf{h}), W_{i}(\textbf{h}) := \left\lbrace \textbf{x} \in \mathbb{R}^{d}: \nabla u_{h} = p_{i} \right\rbrace.
\end{equation*}

\begin{theorem}
    \citep{gu2013variational} Let $\Omega$ be a compact convex domain in $\mathbb{R}^{d}$, $\left\lbrace p_{1},p_{2}, \dots ,p_{n} \right\rbrace $ be a set of $n$ distinct points in $\mathbb{R}^{d}$, and $f:\Omega \longrightarrow \mathbb{R}$ be a continuous function. For any $\nu_{1},\nu_{2}, \dots ,\nu_{n}>0$ with $\sum_{i=1}^{n}\nu_{i}=\int_{\Omega}f(x)dx$, there exists $\textbf{h} = \left( h_{1},h_{2}, \dots ,h_{n}   \right) \in \mathbb{R}^{n}$, unique up to adding  a constant $(c, c, \dots , c)$, such that $\forall \ i \in \left\lbrace  1,\dots ,n\right\rbrace$,
\begin{equation}\label{measure}
	\tau( W_{i}(\textbf{h})\cap\Omega) = \omega _{i} (\textbf{h}) := \int_{W_{i}(\textbf{h})\cap\Omega}f(x)dx=\nu_{i},
\end{equation}
where the $\tau$-volume of each cell denotes its probability measure. 
\end{theorem} 
The height vector $\textbf{h}$ is exactly the optimal solution of the convex energy function
\begin{equation*}
	E(\textbf{h}) = \int_{0}^{\textbf{h}} \sum_{i=1}^{n}\tau( W_{i}(\textbf{h})\cap\Omega)dh_{i}-\sum_{i=1}^{n}h_{i}\nu_{i},
\end{equation*}
on the open convex set (the admissible solution space) 
$
	H = H_{1}\cap H_{2},
$
where $ H_{1}$ and $ H_{2}$ are defined as
\begin{equation}\label{spaceH1}
	H_{1} =\left\lbrace \textbf{h}\in \mathbb{R}^{n} | 	\tau( W_{i}(\textbf{h})\cap\Omega) > 0 , \forall \ i \in \left\lbrace  1,\dots ,n\right\rbrace \right\rbrace ,
\end{equation}
\begin{equation}\label{spaceH2}
	H_{2} = \left\lbrace \sum_{i=1}^{n}h_{i}=0 \right\rbrace .
\end{equation}
Then the semi-discrete OT map is given by $T(x)=\nabla u_{\textbf{h}}(x)$. 
The mass center of each power cell $W_{i}$ is calculated as $m_{i}=\int_{\Omega}xd\tau(x)/\nu_{i}$, then we have a bijective map $\hat{T}$ induced by the OT map from $\tau$ to $\nu$ : $\hat{T}(m_{i}) = p_{i}, \forall \ i =\left\lbrace  1,\dots ,n\right\rbrace$.

\subsection{Geometric Mapping}
Here, we introduce the fundamental geometric mappings involved in this work such as harmonic mapping, conformal mapping and quasiconformal mapping, especially focusing on their discrete cases on triangular meshes. 

\vspace{1mm}
\noindent\textbf{Harmonic Mapping.}
A triangular mesh is denoted as \( M=(V,E,F) \), where \( V, E, F \) represent the vertex, edge and triangular face sets of the mesh, respectively. The harmonic function is defined as \( f:V\longrightarrow \mathbb{\mathbb{R}}^2 \) 
and is achieved by minimizing the harmonic energy,
\begin{equation}\label{eqn:HarmonicEnergy}
	E(f) = \frac{1}{2} \sum_{ [p_{i},p_{j}] \in E} w_{ij}(f(p_{i})-f(p_{j}))^{2},
\end{equation}
where $w_{ij}$ is the cotangent edge weight given by
\begin{equation*}
	w_{ij}=\left\{
	\begin{aligned}
		&cot\theta_{ij}^{k} + cot\theta_{ij}^{l} , &e_{ij} \notin \partial M, \\
		&cot\theta_{ij}^{k}, &e_{ij} \in \partial M,
	\end{aligned}
	\right.
\end{equation*}
 \( \partial M \) indicates the boundary of the triangular mesh, and \( \theta_{ij}^{k} \) represents the corner angle in face \( [ p_{i},p_{j},p_{k} ] \) at \( p_{k} \). 
Then, the discrete Laplace operator $\bigtriangleup f$ is defined as follows:
\begin{equation*}
	\bigtriangleup f(p_{i}) = \sum_{ [p_{i},p_{j}] \in E} w_{ij}(f(p_{i})-f(p_{j})),
 \label{eqn:harmonic_fundtion}
\end{equation*}
and $f$ is obtained by solving the sparse linear system. When \( \Delta f = 0 \), \( E(f) \) defined in Eqn. (\ref{eqn:HarmonicEnergy}) reaches a minimum value.
For a topological disk triangular mesh, when the target domain is convex, its harmonic mapping exists and is unique, and is guaranteed to be a diffeomorphism \citep{gu2008computational}.

\vspace{1mm}
\noindent\textbf{Conformal and Quasiconformal Mappings.}
 Given a Riemann surface $M$ and a target parameter domain $\mathbb{D}$, conformal mapping (parameterization) $\phi:M\longrightarrow \mathbb{D}$ transforms surface metric $g$ to metric $\tilde{g}$ on parameter domain, 
 $\tilde{g}=ge^{\lambda}$, where $\lambda$ is the conformal factor function. 
Quasiconformal mapping $\phi: M\longrightarrow \mathbb{D}$ is defined by the Beltrami equation,
\begin{equation*}
	\dfrac{\partial \phi}{\partial \bar{z}} = \mu(z) \dfrac{\partial \phi}{\partial z},
\end{equation*}
where Beltrami coefficient $\mu$ is a complex-valued function, satisfying $\left\| \mu \right\|_{\infty} <1$ (see Fig. \ref{fig:QC-mu}).  
 Conformal mapping maps infinitesimal circles to infinitesimal circles, while quasiconformal mapping maps infinitesimal circles to infinitesimal ellipses. Intuitively, conformal mapping preserves angle structures, and quasiconformal mapping permits bounded angle distortions. Conformal mapping is a special case of quasiconformal mapping when $\mu=0$. 
  
In discrete case, a triangular mesh $M = \left( V, E, F \right) $ and the Beltrami coefficient function $\mu$ defined on vertices are given. Consider a triangular face $\left[p_{i},p_{j},p_{k}\right]$ with the corresponding intrinsic coordinates $z_{i}$,$z_{j}$, and $z_{k}$ and Beltrami coefficients $\mu(p_{i})$, $\mu(p_{j})$, and $\mu(p_{k})$. The Beltrami coefficient on edge $e_{ij} = \left[ p_{i},p_{j} \right] $ is given by 
$
\mu_{ij} = \frac{ \mu(p_{i}) + \mu(p_{j}) }{2}.
$
The \textit{auxiliary metric} \citep{zeng2009surface} associated with $\mu_{ij}$ is defined as
\begin{equation}\label{auxiliary metric}
	\hat l_{ij} = \left| dz_{ij} + \mu_{ij} \cdot d\bar{z}_{ij} \right|, 
\end{equation}
where $dz_{ij} = z_{i}-z_{j}$ and $d\bar{z}_{ij}$ is the conjugate of $dz_{ij}$. 

Therefore, a quasiconformal mapping associated with a Beltrami coefficient function can be converted to a conformal mapping with the corresponding auxiliary metric \citep{zeng2009surface,zeng2012computing}. Therefore, under the auxiliary metric, all existing conformal mapping methods can be employed to compute quasiconformal mapping, and quasiconformal mapping is guaranteed to be a diffeomorphism in theory.

\section{Topology-Preserving  Optimal Transport (t-OT)}
\label{sec:algorithm}

\begin{figure}[h]
\centering
\includegraphics[width=1\textwidth]{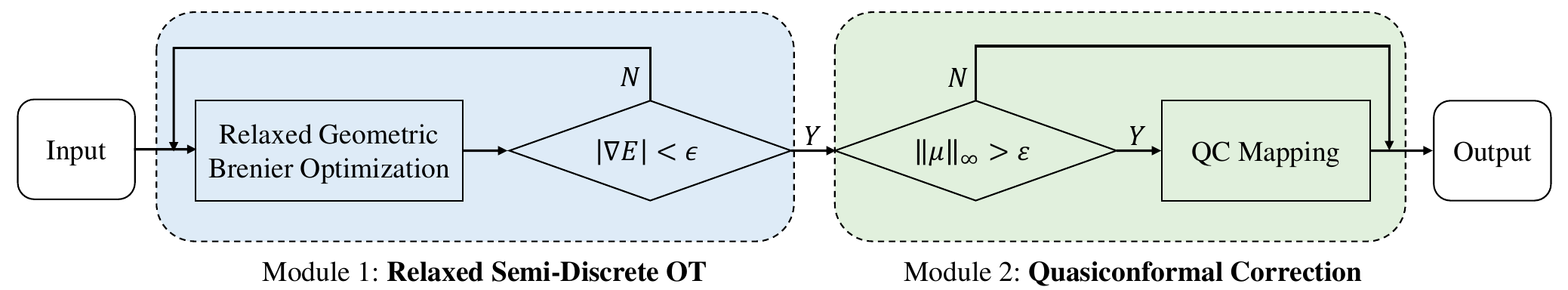}\\
\caption{The workflow of the t-OT method.
}
\label{Flow2}
\end{figure}

Suppose two simply-connected domains, 
$(\Omega,M,\tau)$, $(\Omega^{*},M^{*},\nu)$, are given with density functions $f$ and $g$, respectively. In the semi-discrete OT framework, the mesh structure of $M$ is not used in the computation since  the probability measure is set to be uniform everywhere. Therefore, $M$ can be chosen flexibly. Intuitively, the mesh $M$ here is regarded as an auxiliary reference. For a given surface, the computation of OT starts in the conformal parameter domain which intrinsically preserves the geometric structure. Therefore, we first compute the conformal mapping $\phi: (\Omega^{*},M^{*}) \to (\mathbb{D}, M_0)$ of the input triangular mesh to canonical domain $\mathbb{D}$, e.g., disk, square and rectangle. Specially, if the object to be operated is an image, we then directly converting it to a triangular mesh with square or rectangle boundary as the input of t-OT.  
Next, we try to compute the optimal transport map $\hat{T}_1: (\Omega,M,\tau) \rightarrow (\mathbb{D},M_0,\nu)$ that preserves the overall probability measure and topological structure, while minimizing the transportation cost. Based on the baseline variational sd-OT algorithm \citep{gu2013variational}, we relax the admissible space of the solution and then append QC mapping to refine unexpected distortions and keep the original topology. As shown in Fig. \ref{Flow2}, the workflow includes two key modules: Relaxed semi-discrete OT and Quasiconformal correction. 

Accordingly, the computational procedure is presented in Algorithm \ref{alg:1}, in which two modules are further expanded in Algorithm \ref{alg:2} (containing Algorithm \ref{alg:2.1}) and Algorithm \ref{alg:3}, respectively. The whole pipeline is briefly configured as follows: 
\begin{enumerate}
	\item Compute  the relaxed semi-discrete OT map $\hat{T}_1$ 
 without changing the original topology, $T_1:=\hat{T}^{-1}_1$ 
 (see Algorithm \ref{alg:2});
	\item Compute the QC mapping $T_2$ 
 with the corrected Beltrami coefficient $\mu_{T_1}$, 
 $T:=T_2\circ \phi$
 (see Algorithm \ref{alg:3}). 
\end{enumerate}

\setcounter{algorithm}{0}
\renewcommand{\thealgorithm}{1}
\begin{algorithm}
\caption{Topology-Preserving Optimal Transport (t-OT)} 
\label{alg:1}
\begin{algorithmic}[1]
\Require The domains with triangular meshes $(\Omega,M,\tau)$, $(\Omega ^{*},M^*,\nu)$, satisfying $\tau(\Omega)=\nu(\Omega ^{*})$; conformal mapping $\phi:(\Omega ^{*},M^*)\to(\mathbb{D}, M_0)$, 2D domain $\mathbb{D} \subset \Omega$; the hyper-parameters $\epsilon$, $\varepsilon$, and $\lambda$;
\Ensure The map $T$ and the generated mesh $\hat M$.
\State Construct the measure for mesh $M_0$ from $\nu$, $\nu(p_i)=\nu_i$ for mesh vertex $p_i$;
\State Compute the Relaxed Semi-Discrete OT map $\hat{T}_1: (\Omega,M) \rightarrow (\mathbb{D}, M_0)$ and update $\hat{M}$ by Algorithm \ref{alg:2} ;
\State Compute Beltrami coefficient $\mu_{T_1}$ of $T_1:=\hat{T}^{-1}_1:M_0 \rightarrow \hat{M}$;
\If{$\left\| \mu_{T_1} \right\|_{\infty} > \varepsilon$}
\State Compute the QC mapping $T_2: (\mathbb{D}, M_0) \rightarrow (\Omega ^{*},\hat{M})$ associated with the corrected $\mu_{T_1}$ by Algorithm \ref{alg:3};
\State Update $\hat{M} \leftarrow T_2(M_0)$;
\EndIf
\State \textbf{Return} The map $T:= T_2\circ \phi$ and the generated mesh $\hat{M}$.
\end{algorithmic}
\end{algorithm}

\setcounter{algorithm}{0}
\renewcommand{\thealgorithm}{2}
\begin{algorithm}
\caption{Relaxed Semi-Discrete OT} 
\label{alg:2}
\begin{algorithmic}[1]
\Require The domains with triangular meshes $(\Omega,M,\tau)$, $(\mathbb{D}, M_0,\nu)$; 
the prescribed measure $\nu(p_i)$; the hyper-parameters $\epsilon$ and $\lambda$;
\Ensure The map $\hat{T}_1: (\Omega,M) \rightarrow (\mathbb{D}, M_0)$ and the generated mesh $\hat{M}$.
\State Initialize $\hat{M} \leftarrow M$, $\textbf{h} = \frac{1}{2} ( \left| p_{i} \right| ^{2} -1 )$;
\While $\left| \nabla E \left( \textbf{h} \right)  \right| > \epsilon$
\State $\textbf{h}_0 \leftarrow \textbf{h}$; 
\State Conduct the relaxed geometric Brenier optimization process by Algorithm \ref{alg:2.1} for computing $\hat{T}_1: \hat{M} \rightarrow M_0$ and updating the height vector $\textbf{h}$, $\left| \nabla E \left( \textbf{h} \right) \right|$, $\textbf{d}_0$ and $\hat{M}$;
\EndWhile
\State \textbf{Return} The map $\hat{T}_1: (\Omega,M) \rightarrow (\mathbb{D}, M_0)$ and the generated mesh $\hat{M}$;
\end{algorithmic}
\end{algorithm}

In detail, for the regions with unexpected distortions (including flips) in $\hat{M}$, described by the Beltrami coefficient $\left\| \mu_{T_1} \right\|_{\infty} > \varepsilon$ ($\varepsilon \in [a ,1]$, $a>0$), we use harmonic function to refill Beltrami coefficients in the bad regions, and then compute the QC mapping associated with the updated $\mu_{T_1}$ to smooth the distortions.

\subsection{Relaxed Semi-Discrete Optimal Transport}\label{ROT}
\setcounter{algorithm}{0}
\renewcommand{\thealgorithm}{2.1}
\begin{algorithm}
\caption{Relaxed Geometric Brenier Optimization} 
\label{alg:2.1}
\begin{algorithmic}[1]
\Require The domains with triangular meshes $(\Omega,M,\tau)$, $(\mathbb{D}, M_0,\nu)$;  
the previous height vector $\textbf{h}_0$, $\left| \nabla E \left( \textbf{h}_0 \right) \right|$; the previous vector $\textbf{d}_0$; the prescribed measure $\nu _i = \nu(p_i)$, and hyper-parameter $\lambda$;
\Ensure The map $\hat{T}_1: \hat{M} \rightarrow M_0$, the height vector $\textbf{h}$, $\left| \nabla E \left(\textbf h \right) \right|$, $\textbf{d}_0$ and the generated mesh $\hat{M}$.
\While{$\widetilde{\Delta} \left( \left| \nabla E \right| \right) > 0$}
\State $\textbf{h} \leftarrow \textbf{h}_0 + \lambda \textbf{d}_0$;
\State Compute the lower convex hull of $\left\lbrace \left( p_{i},-h_{i} \right)  \right\rbrace _{i=1}^{n}$;
\State  Compute the upper envelope $\left\lbrace x\cdot p_{i} + h_{i} \right\rbrace _{i=1}^{n}$;
\State Project the upper envelope to the plane to get the power diagram $\cup _{i=1}^{n} W_{i}(\textbf{h})$;
\State  Compute the probability measure $w_{i}=\tau(W_{i}(\textbf{h})\cap\Omega)$;
\State Compute the gradient of $E(\textbf{h})$, $\nabla E(\textbf{h}) = (\nu_{1}-w_{1},\dots,\nu_{n}-w_{n})$; 
\State Compute the gradient difference $\widetilde{\Delta} \left( \left| \nabla E \right| \right) = \left( \left| \nabla E(\textbf{h}) \right| \right) - \left( \left| \nabla E(\textbf{h}_0) \right| \right)$;
\State $\lambda \leftarrow \frac{\lambda}{2}$;
\EndWhile
\State Solve the linear equation system $Hess(E(\textbf{h}))\textbf{d}=\nabla E(\textbf{h})$;
\State Compute the mass centers $\left\lbrace m_{i} \right\rbrace _{i=1}^{n}$ of $\left\lbrace W_i(\textbf{h})\right\rbrace$ to update $\hat{M}$ with $\left\lbrace m_{i} \right\rbrace _{i=1}^{n}$;
\State Update the vector $\textbf{d}_0\leftarrow \textbf{d}$ and the map $\hat{T}_1(m_{i})=p_i$;
\State \textbf{Return} The map $\hat{T}_1: \hat{M} \rightarrow M_0$, the height vector $\textbf{h}$, $\left| \nabla E(\textbf h) \right|$, $\textbf{d}_0$ and the generated mesh $\hat{M}$.
\end{algorithmic}
\end{algorithm}

The computational algorithm for relaxed sd-OT map is depicted in Algorithm \ref{alg:2}. The key operation of the relaxed geometric Brenier optimization is given in Algorithm \ref{alg:2.1}. Suppose a simply connected domain $\Omega^{*}$ with relative density function $g$, there are $n$ distinct vertices $\left\lbrace p_{i} \right\rbrace _{i=1}^{n}$ in $\Omega ^{*}$, then the probability measure of $\Omega^{*}$ is defined as $\nu = \sum_{i=1}^{n}\nu_{i}\delta (y-p_{i})$, where $\delta$ is the dirac delta function.
For any vertex $p_{i} \in \left\lbrace p_{i}  \right\rbrace _{i=1}^{n} $, we have $\nu_{i} = g(p_{i})$ and then normalize it as $\nu_{i} = \frac{\nu_{i}}{\sum_{i=1}^{n}\nu_{i}}$ for each vertex $p_{i}$. $\Omega$ is larger than $\Omega^{*}$ and the density function $f$ is specified to be the uniform.

The relaxed sd-OT method mainly aims at minimizing the energy function $E(\textbf{h})$ by using Newton's method. In the very beginning, for each vertex $p_{i}$, we initialize the height vector as $\textbf{h}_{0} = \frac{1}{2} ( \left| p_{i} \right| ^{2} -1 )$. Then, the height vector $h$ in each interaction can be used to form the lower convex hull $\left\lbrace \left( p_{i},-h_{i} \right)  \right\rbrace _{i=1}^{n}$. By projecting the lower convex hull, a triangulation of the $\left\lbrace p_{i} \right\rbrace _{i=1}^{n}$ is induced while maintaining the original topological structure.

The gradient of the energy  $E(\textbf{h})$ is given by 
\begin{equation*}
	\nabla E(\textbf{h}) = \left( \nu_{1} - \omega _{1} (\textbf{h}), \nu_{2}-\omega _{2} (\textbf{h}), \dots , \nu_{n}-\omega _{n} (\textbf{h})\right) ,
 \label{eqn:gradient_energy}
\end{equation*}
where $\omega_i$ is defined in Eqn. (\ref{measure}).
Based on the previous height vector $\textbf{h}_0$ and the current height vector $\textbf{h}$, the difference of $\left| \nabla E \right|$ is defined by $\widetilde{\Delta} \left( \left| \nabla E \right| \right) = \left( \left| \nabla E(\textbf{h}) \right| \right) - \left( \left| \nabla E(\textbf{h}_0) \right| \right)$.
The step length parameter $\lambda$ is updated until $\widetilde{\Delta} \left( \left| \nabla E \right| \right) < 0$.

The Hessian matrix $Hess(E(\textbf{h}))$ of the energy $E(\textbf{h})$ for the off diagonal elements can be constructed as the following equation:
\begin{equation*}
	Hess(E(\textbf{h}))_{ij} =-\frac{\tau \left( W_{i}(\textbf{h}) \cap  W_{j}(\textbf{h})\right) }{\left| p_{i}-p_{j} \right| }= \frac{-\tau (\bar{e}_{ij})}{\left| {e}_{ij} \right| },
\end{equation*}
where ${e}_{ij}$ is the edge in $\Omega^{*}$ that connects the two vertices $p_{i}$ and $p_{j}$, and $\bar{e}_{ij}$ is the dual edge in the power diagram $\mathcal{T}(\textbf{h})$. The diagonal elements of $Hess(E(\textbf{h}))$ are defined as
\begin{equation*}
	Hess(E(\textbf{h}))_{ii} =\sum_{i\sim j} \frac{\tau \left( W_{i}(\textbf{h}) \cap  W_{j}(\textbf{h})\right) }{\left| p_{i}-p_{j} \right| } = \sum_{i\sim j}\frac{\tau (\bar{e}_{ij})}{\left| {e}_{ij} \right| },
\end{equation*}
where the symbol $\sim$ represents the neighborhood relationship between the power diagrams  $W_{i}$ and $W_{j}$. Then we have a linear function system 
\begin{equation}
	Hess(E(\textbf{h}))\textbf{d}=\nabla E(\textbf{h}).
 \label{eqn:linear_function_system}
\end{equation}
By solving Eqn. (\ref{eqn:linear_function_system}), the update direction vector $\textbf{d}$ is given and the height vector is updated by $\textbf{h}\leftarrow \textbf{h}+\lambda \textbf{d}$. Compared with Eqn. (\ref{spaceH1}), our admissible space $H = \hat{H_{1}} \cap H_{2} $, in which $\hat{H_{1}}$ is defined as 
\begin{equation*}
	\hat{H_{1}} =\left\lbrace \textbf{h}\in R^{n} | 	\tau( W_{i}(\textbf{h})\cap\Omega) \geqslant0 , \forall i \in \left\lbrace  1, \dots ,n\right\rbrace \right\rbrace,
\end{equation*}
and $H_{2}$ is given in Eqn. (\ref{spaceH2}).
The preceding steps will be repeated until $\left| \nabla E(\textbf{h}) \right|<\epsilon$, where the hyper-parameter $\epsilon$ has been specified beforehand. Then the mass center $m_{i}$ of each convex polygon $W_{i}$ and the one-to-one correspondence relationship $\hat{T}_1: \hat{T}_1 (m_{i}) = p_{i}$ are generated.
Subsequently a new adaptive mesh is created.

\subsection{Quasiconformal Correction}\label{QC_Correction}

\setcounter{algorithm}{0}
\renewcommand{\thealgorithm}{3}
\begin{algorithm}[t]
\caption{Quasiconformal Correction} 
\label{alg:3}
\begin{algorithmic}[1]
\Require The triangular meshes $M_0$ and $\hat{M}$ with the same topology, i.e., 
$T_1 : M_0 \to \hat{M}$, with the associated Beltrami coefficient 
$\mu_{T_1}$; the hyper-parameter $\varepsilon$;
\Ensure The QC mapping $T_2: M_0 \to \hat{M}$ and generated mesh $\hat{M}$ with correction.
\For{\textbf{all} $p_i$} 
\State Compute the edge weight $\omega$ of $M_0$;
\If{$|\mu_{T_1}({p_i}) | > \varepsilon $}
\State Find the $\gamma-$ ring neighboring patch $P_i$ of $p_i$ with corresponding convex patch $D_i$ on $\hat M$, $\partial P_i = \{p_i^j\} $;  
\State Fill the Beltrami coefficients $\hat \mu_i$ for $P_i$ by harmonic diffusion  
  using edge weight $\omega$ and $\hat \mu_i(\partial P_i)= \mu_{T_1}(\partial P_i)$;
\State Compute the auxiliary metric $\hat l_i$ under $\hat \mu_i$ and the corresponding edge weight $\hat \omega_i$ for the patch $P_i$;
\State Compute the harmonic map $\bar{T}_i: P_i\to D_i$ using edge weight $\hat \omega_i$
 with $\partial D_i = \bar{T}_i(\partial P_i)$;
\State Update $\hat M$ with $\bar{T}_i(P_i)$;
\EndIf
\EndFor
\State \textbf{Return} The QC mapping $T_2: M_0\to \hat M$ and the corrected mesh $\hat{M}$.
\end{algorithmic}
\end{algorithm}

The QC mapping is used for mesh correction (see Algorithm \ref{alg:3}). 
Suppose a mapping $T_1: M_0 \to \hat{M}$ with a given Beltrami coefficient $\mu$. 
Then for the vertex $p_{i}$ in $M_0$ with $|\mu({p_i}) | > \varepsilon $, $P_{i}$ on $M_0$ is the $\gamma-$ring neighbor patch of $p_{i}$ and $D_{i}$ is the corresponding convex patch on $\hat{M}$, where $\gamma \in N^{+}$ is a positive integer ($\gamma = 2$ in the experiments). The boundary $\partial P_i = \{p_i^j\} $ serves as boundary landmark constraints. 
The Beltrami coefficient of the interior patch is filled by the Beltrami coefficient of the boundary landmarks. 
Specifically, for \(P_i\), the boundary satisfies \(\hat \mu_i(\partial P_i) = \mu(\partial P_i)\), while the interior relative Beltrami coefficients \(\hat \mu\) are determined by harmonic diffusion using edge weight \(\omega\).
The auxiliary metric $\hat l_i$ is computed under $\hat \mu_i$ using Eqn. (\ref{auxiliary metric}). Therefore, the edge weight $\hat \omega$ is computed with the auxiliary metric and the harmonic map $\bar{T}_i: P_i\to D_i, \partial D_i = \bar{T}_i(\partial P_i)$ is computed by using $\hat \omega_i$. The result leads to a QC mapping of $P_i$ with $\hat \mu_i$, where the boundary is mapped to the specified positions. The mapping is guaranteed to be a diffeomorphism, and the unexpected distortions in $\hat{M}$ are then smoothed. 
\section{Temporal Topology-Preserving OT (tt-OT)}\label{sec:tt-OT}

We further study the transportation process based on the proposed t-OT method. We aim to keep the topology-preserving property during the transportation process, which would facilitate analyze physical diffusion or evolution and impact broader applications. 
As we know, dynamic optimal transport is a major theory in the fluid mechanics that relevant to the velocity field and geodesic. Many methodologies on dynamic OT have emerged based on the approach for solving the Monge-Ampere equation proposed by \cite{benamou2000computational}. The methods illustrate the dynamic changes of particles throughout the transportation process, and has continuity, allowing measures to vary continuously over time to accommodate the movement and evolution of the fluid. However, they do not preserve the topological structure.

\begin{figure}[h]
\centering
\includegraphics[width=1\textwidth]{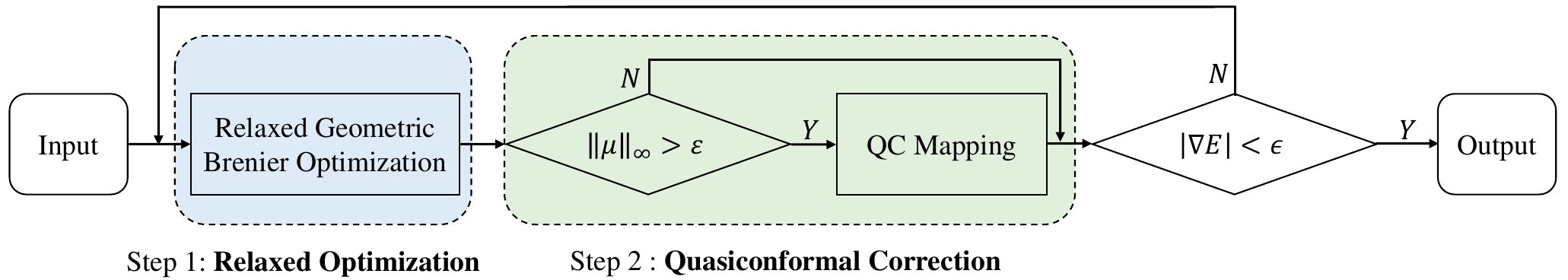}\\
\caption{The workflow of the tt-OT method.} 
\label{fig:TemporalWorkFlow}
\end{figure}

In order to display the transportation process dynamically and continuously, we develop the temporal topology-preserving optimal transport (tt-OT) methodological framework based on t-OT. The computational flow of tt-OT is shown in Fig. \ref{fig:TemporalWorkFlow}. Two main steps are performed in each iteration, namely (1) Relaxed Geometric Brenier Optimization, which generates topology-preserving spatial-temporal results, and (2) QC correction, which detects local unexpected angle distortions and eliminates these distortions through QC mapping while preserving the topology. Therefore, during the optimization process, we obtain a sequence of t-OT maps (with different measures). 
The computational details are described in Algorithm \ref{alg:temporal t-OT}, where two main steps refer to Algorithm \ref{alg:2.1} and Algorithm \ref{alg:3}, respectively.

\setcounter{algorithm}{0}
\renewcommand{\thealgorithm}{4}
\begin{algorithm}[t]
\caption{Temporal Topology-Preserving Optimal Transport (tt-OT)} 
\label{alg:temporal t-OT}
\begin{algorithmic}[1]
\Require The domains with triangular meshes $(\Omega,M,\tau)$, $(\Omega ^{*},M^*,\nu)$, satisfying $\tau(\Omega)=\nu(\Omega ^{*})$; conformal mapping $\phi:(\Omega^{*},M^*)\to(\mathbb{D}, M_0)$, 2D domain $\mathbb{D} \subset \Omega$; the hyper-parameters $\epsilon$, $\varepsilon$, and $\lambda$;
\Ensure The temporally generated mesh list $\{\hat M_{t}\}$.

\State Construct the measure for mesh $M$ from $\nu$, $\nu(p_i)=\nu_i$ for mesh vertex $p_i$;
\State Initialize $\hat{M} \leftarrow M_0$, $h_{i}=\frac{1}{2} ( \left| p_{i} \right| ^{2} -1 )$, $t=1$; 
\While{$\left| \nabla E(\textbf h_0) \right| > \epsilon$}
\State $\textbf h_0 \leftarrow \textbf h$; 
\State Conduct the relaxed geometric Brenier optimization process by Algorithm \ref{alg:2.1} and compute $\hat{T}_1: \hat{M} \rightarrow M_0$ with $\textbf h_0$, and then update $\textbf h$ and the generated mesh $\hat{M}$; 
\State Compute Beltrami coefficient $\mu_{T_1}$ of $T_1:=\hat{T}_1^{-1} : M_0 \to \hat{M}$;
\If{$\left\| \mu_{T_1} \right\|_{\infty} > \varepsilon$} 
\State Compute the QC mapping $T_2: (\mathbb{D}, M_0) \rightarrow (\Omega ^{*},\hat{M})$ associated with the corrected $\mu_{T_1}$  by Algorithm \ref{alg:3}; 
\State Update $\hat{M} \leftarrow T_2(M_0)$;
\EndIf
\State Update $\hat{M}_{t} \leftarrow \hat{M}$, $t \leftarrow t + 1$;
\EndWhile
\State \textbf{Return} The temporal mesh list $\{ \hat{M}_{t} \}$.
\end{algorithmic}
\end{algorithm}

\section{Algorithm Analysis}\label{sec: Algorithm Analysis}

\subsection{Example}

\begin{figure}[h]
\setlength{\fboxsep}{0pt}
\centering
\footnotesize
\begin{tabular}{ccc}		
\includegraphics[height=0.2\linewidth]{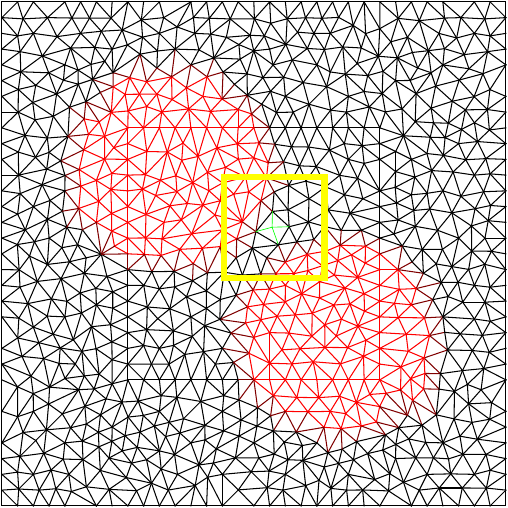}&
\includegraphics[height=0.2\linewidth]{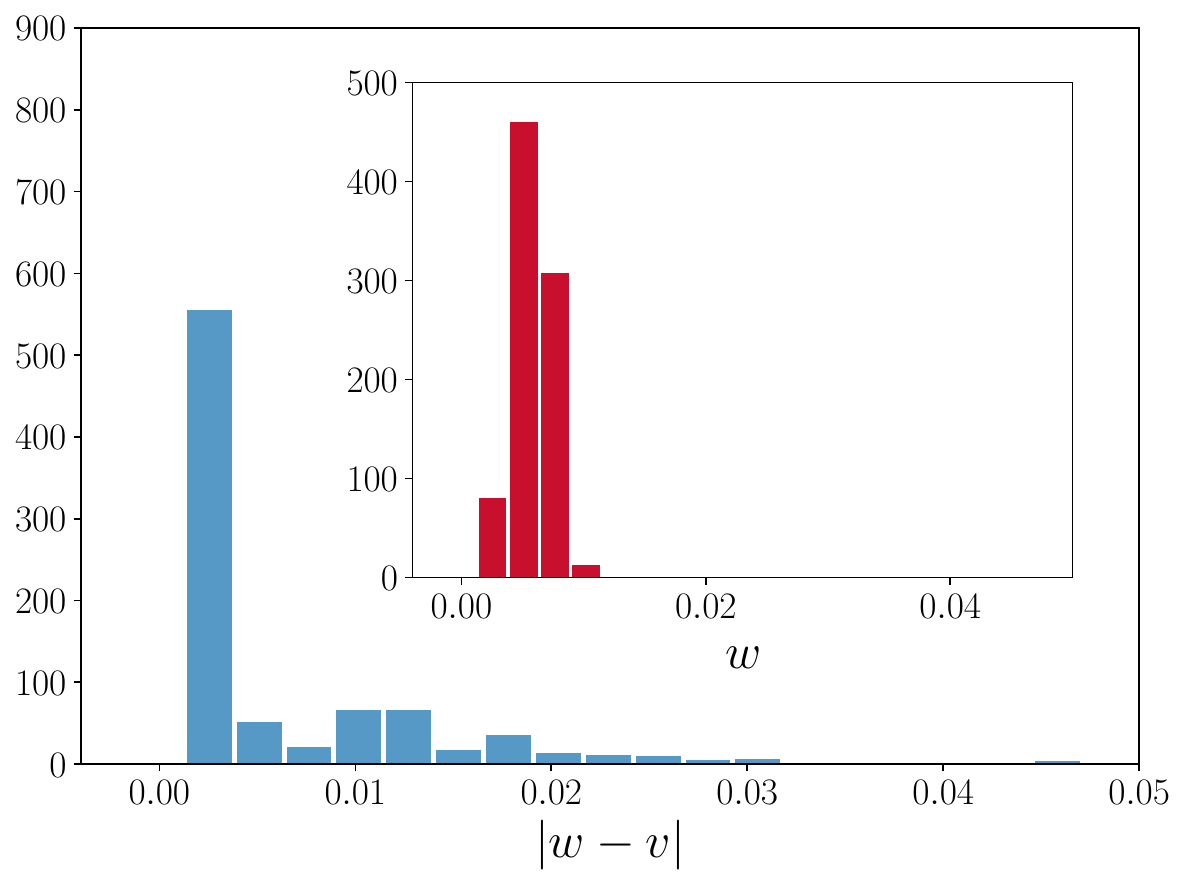}&
\fbox{\includegraphics[trim=0 0 0 0 , clip, height=0.2\linewidth]{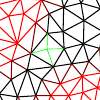}}\\
(a) Original Mesh &(b) Measure Distribution &(c) Zoomed-in Focus  \\
\end{tabular}
\captionsetup{justification=justified, singlelinecheck=false}
\caption{The original mesh and measure distribution of $M$. The focal vertex of $M$ and adjacent edges are  highlighted in green for further analyzing the transport map performance (note that its $rgb=0$).}
\label{t-OT_input_2Area}
\end{figure}

The t-OT map is carried out by integrating relaxed sd-OT and quasiconformal correction. Here, we use an example to illustrate the computation of t-OT. As shown in Fig. \ref{t-OT_input_2Area}, the input mesh $M$ consists of $n=861$ vertices and $1603$ triangular faces, 
defined on the source domain $\Omega ^{*}=\left[ -1, 1 \right] ^{2}$, with the prescribed measure distribution. 
The measure $\nu$ is defined by the rgb values of mesh vertices, 
\begin{equation}\nu_{i} =k * (\delta + rgb) * a_{i}, 
\end{equation}
where we set the customizable density scalar $k=4.0$, the parameter $\delta=0.02$ to perturb the value $rgb$ considering the case of $rgb=0$, and the local area 
	$a_{i} = \frac{1}{3} \sum_{f \in \sigma _{i}} Area(f)$,
in which $\sigma _{i}$ denotes all the triangular faces adjacent to vertex $v_{i}$.
Therefore, for a relatively uniform mesh, the probability measure for red vertex is larger. 
The measure is then normalized as $\nu_i := \frac{\nu_i}{\sum_{i=1}^{n} \nu_{i}}$, and the tolerance in Algorithm \ref{alg:1} is set to be $\epsilon = 1e-5$. The reference domain $\Omega$ is set slightly larger than $\Omega ^{*}$ to ensure that the power diagram of it is fully contained within the reference domain, set as 
$\Omega =\left[ -1.2, 1.2 \right]^{2}$. 
The density function defined on $\Omega$ is \textit{uniform}. 

\begin{figure*}[t]
\centering
\footnotesize
\setlength{\fboxsep}{0pt}
\setlength{\tabcolsep}{2pt}
\begin{tabular}{ccccc}
Output Mesh&Angle Distortion $\left| \mu\right|$&& Measure Distribution&Local Structure\\

\includegraphics[height=0.19\textwidth]{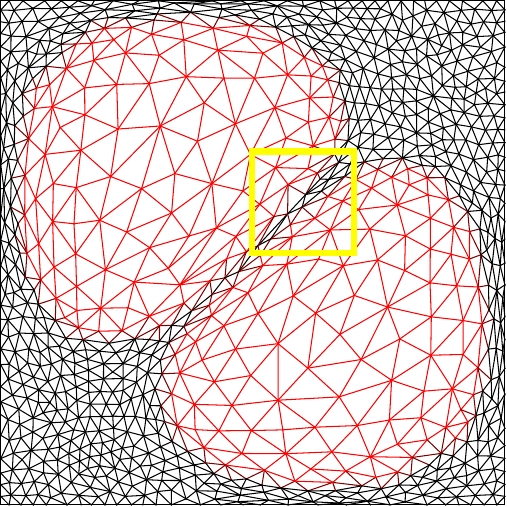}&\includegraphics[height=0.19\textwidth]{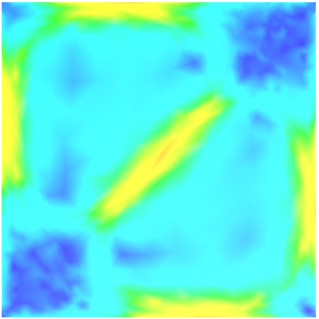}&
\hspace{-1em}
\includegraphics[height=0.19\textwidth]{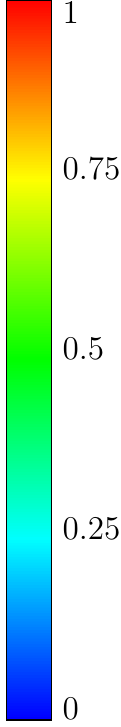}&
\includegraphics[height=0.19\textwidth]{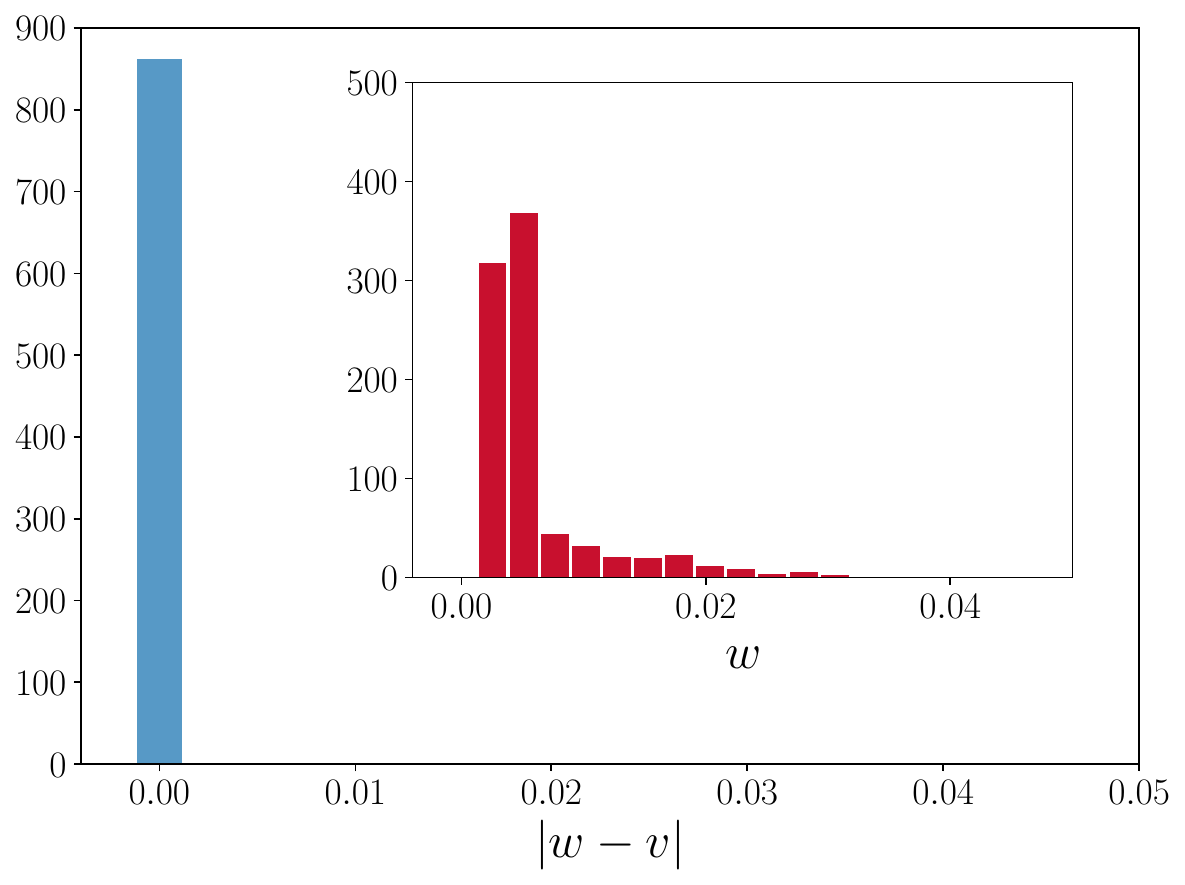}&
\fbox{\includegraphics[height=0.19\textwidth]{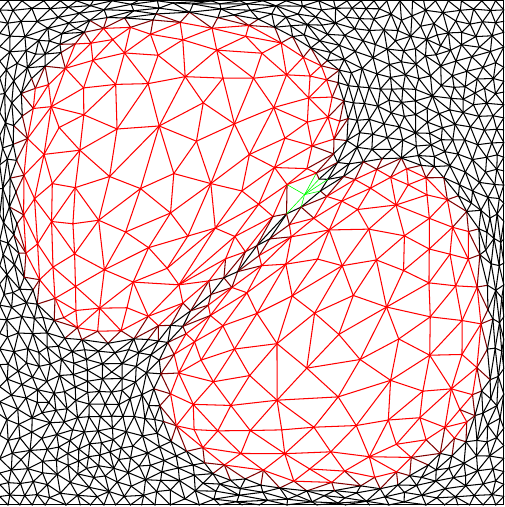}}\\
\multicolumn{5}{c}{(a) sd-OT: Output mesh $\bar{M}$}\\
 
\includegraphics[height=0.19\textwidth]{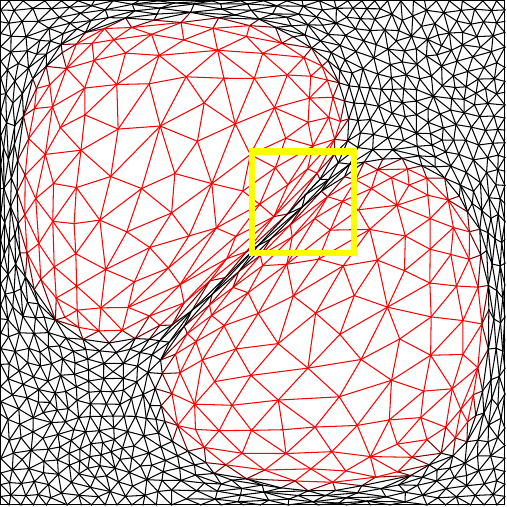}&		\includegraphics[height=0.19\textwidth]{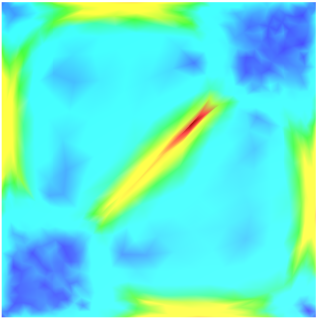}&\hspace{-1em}
\includegraphics[height=0.19\textwidth]{img/boxer5Tex5/colorbar1.pdf} &		\includegraphics[height=0.19\textwidth]{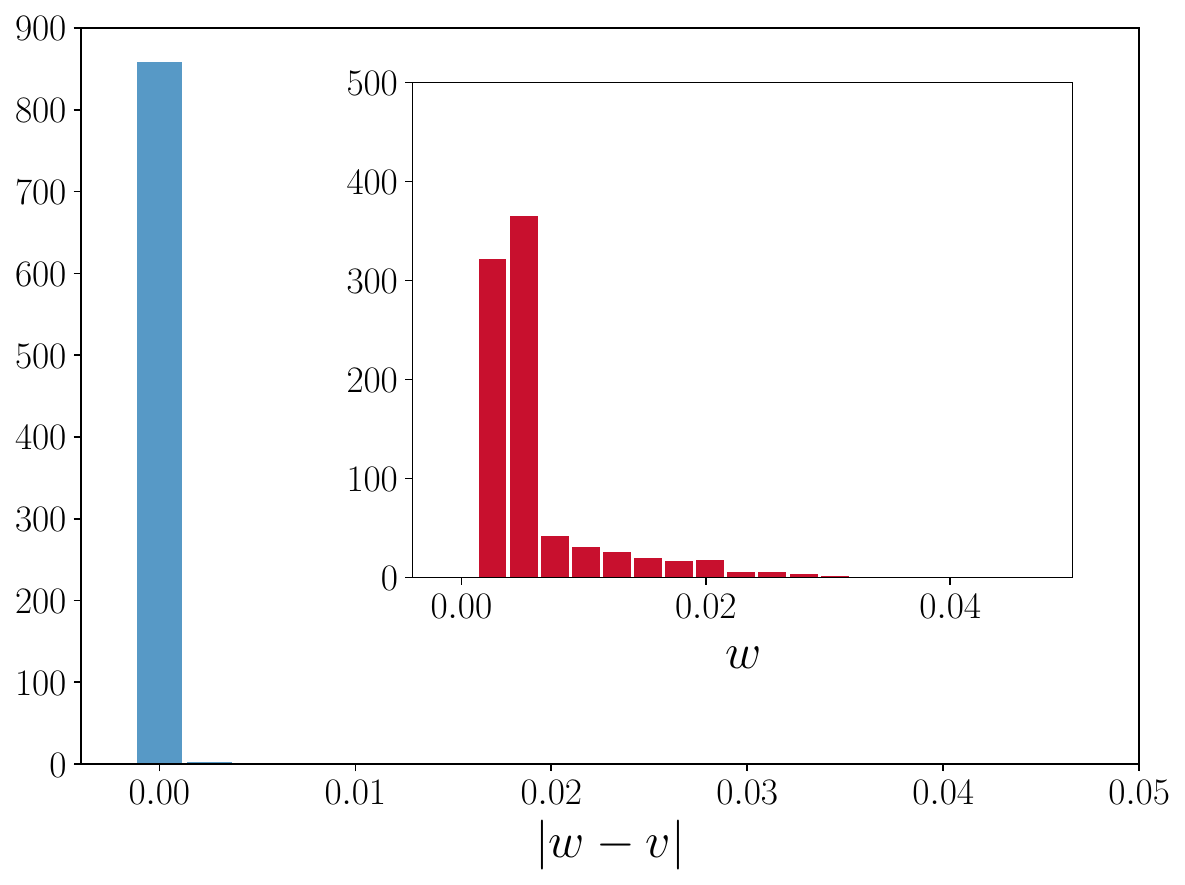}&\fbox{\includegraphics[width=0.19\textwidth]{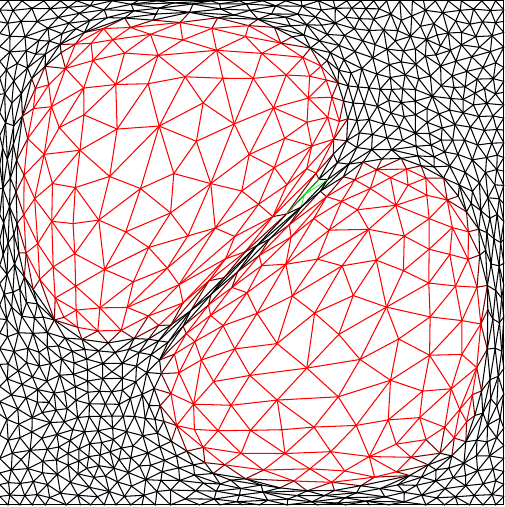}}\\
\multicolumn{5}{c}{(b) sd-OT with original topology: $\check{M}$}\\

\includegraphics[height=0.19\textwidth]{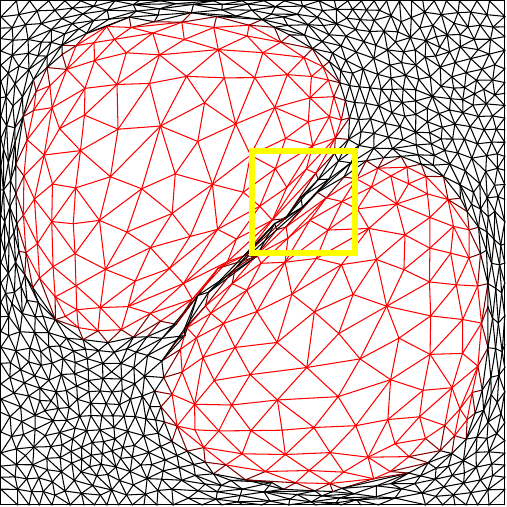}&
\includegraphics[height=0.19\textwidth]{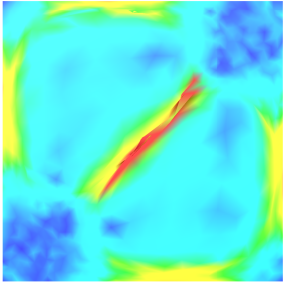}&
\hspace{-1em}
\includegraphics[height=0.19\textwidth]{img/boxer5Tex5/colorbar1.pdf} &	\includegraphics[height=0.19\textwidth]{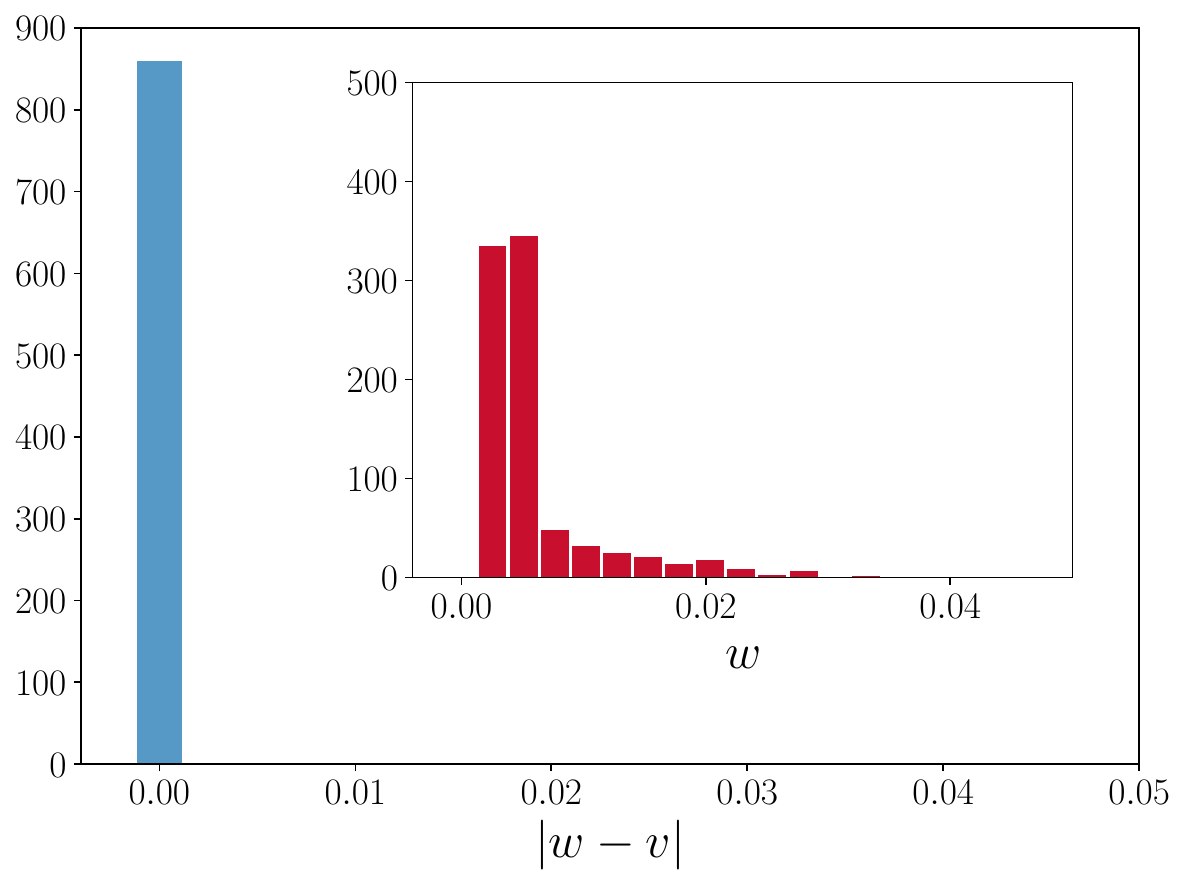}&
\fbox{\includegraphics[width=0.19\textwidth]{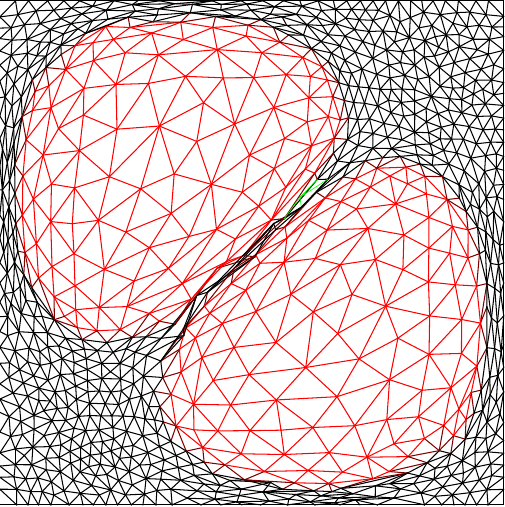}}\\
\multicolumn{5}{c}{(c) Relaxed sd-OT: output mesh $\hat{M}$}\\

\includegraphics[height=0.19\textwidth]{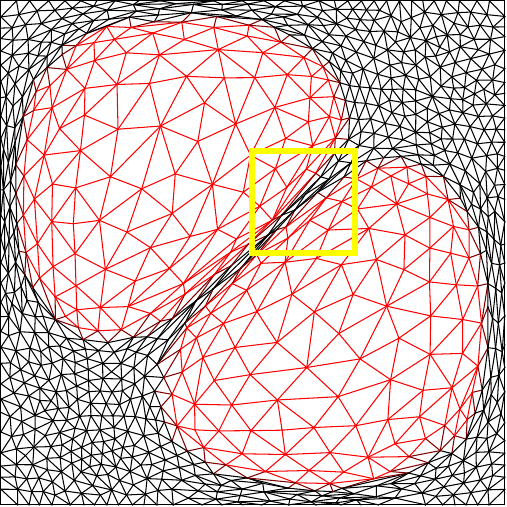}&
\includegraphics[height=0.19\textwidth]{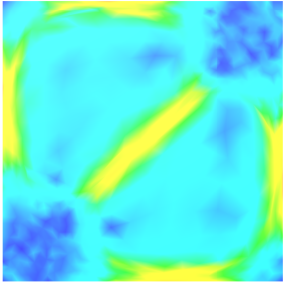}&
\hspace{-1em}
\includegraphics[height=0.19\textwidth]{img/boxer5Tex5/colorbar1.pdf} &	\includegraphics[height=0.19\textwidth]{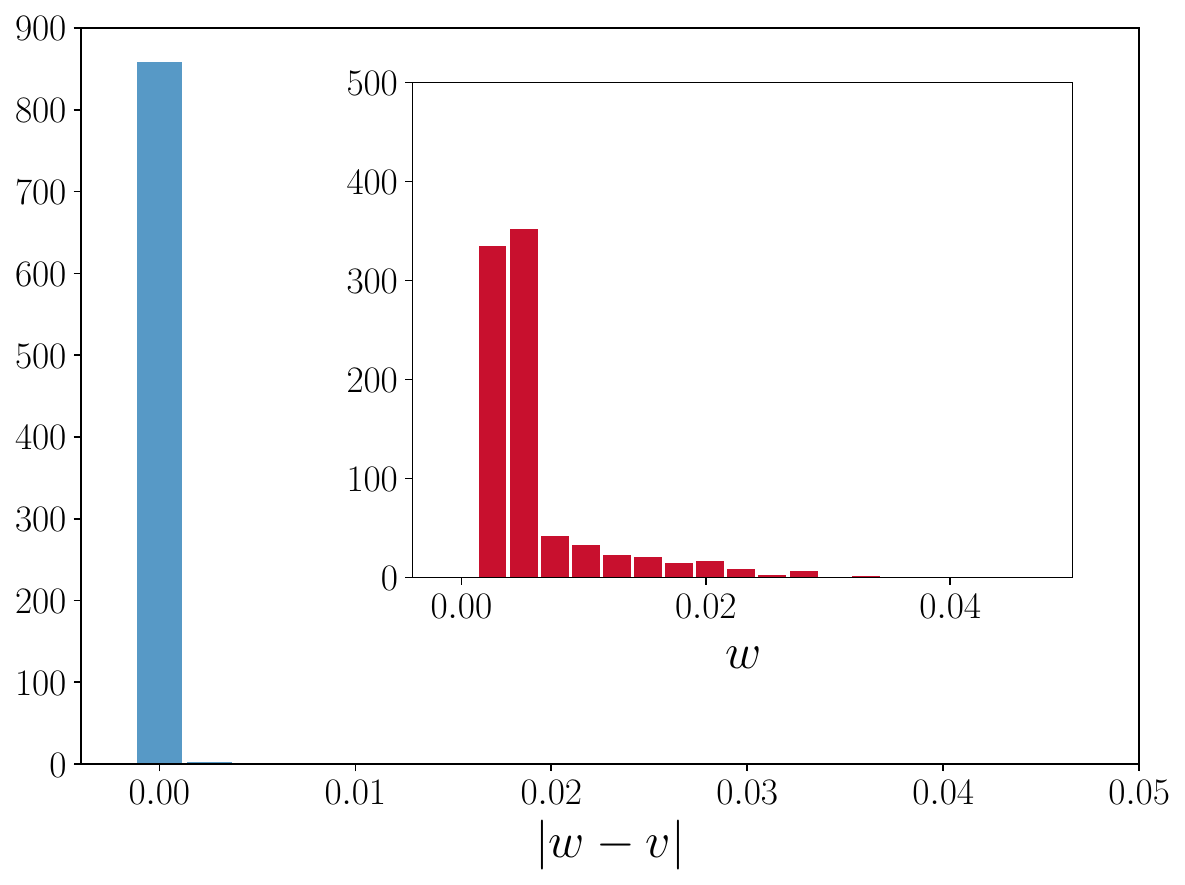}&
\fbox{\includegraphics[width=0.19\textwidth]{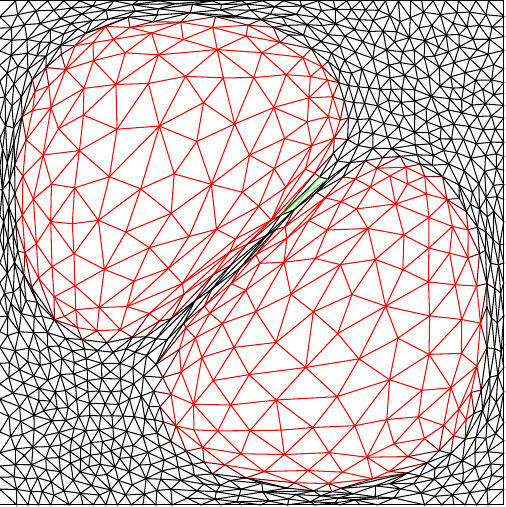}}\\
\multicolumn{5}{c}{(d) t-OT: output mesh $\hat{M}_{\varepsilon = 1.0}$}\\
\includegraphics[height=0.19\textwidth]{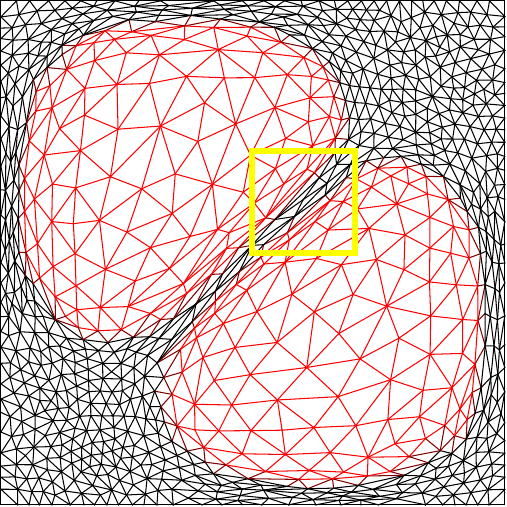}&
\includegraphics[height=0.19\textwidth]{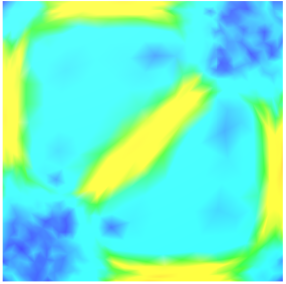}&
\hspace{-1em}
\includegraphics[height=0.19\textwidth]{img/boxer5Tex5/colorbar1.pdf} &	\includegraphics[height=0.19\textwidth]{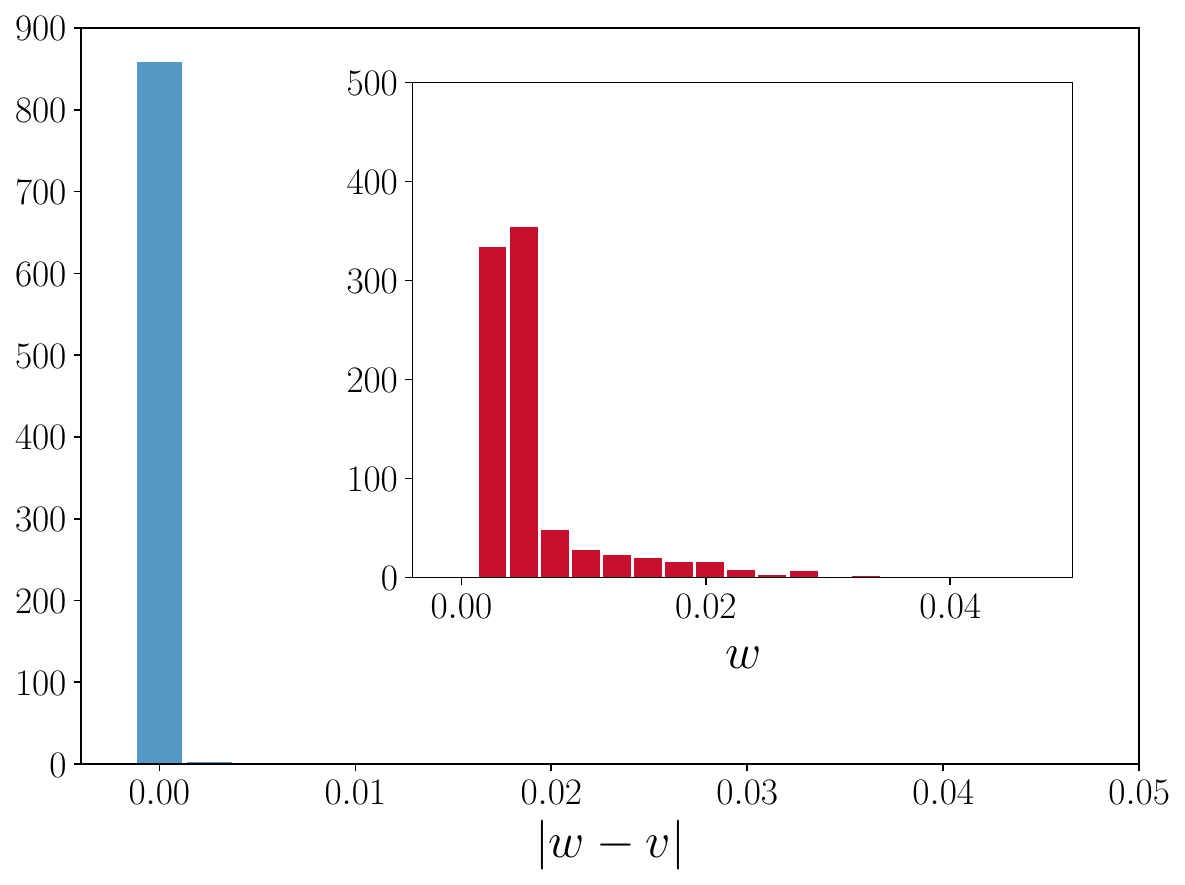}&
\fbox{\includegraphics[width=0.19\textwidth]{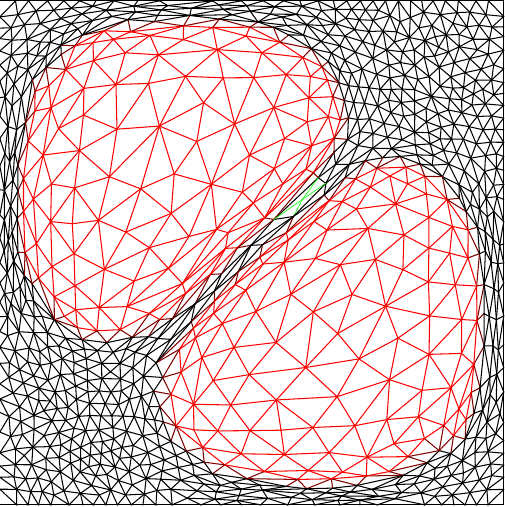}}\\
\multicolumn{5}{c}{(e) t-OT: output mesh $\hat{M}_{\varepsilon = 0.7}$}\\
\end{tabular}\\
\caption{Comparison: sd-OT, sd-OT with original topology, relaxed sd-OT, and t-OT.} 
\label{fig:t-OT_result_2Area}
\end{figure*}

Figure \ref{fig:t-OT_result_2Area} illustrates the comparable results for related methods, including: (a) the baseline sd-OT, (b) the option of transferring original topology to sd-OT result, (c) the relaxed sd-OT without QC correction, and (d-e) the t-OT combing relaxed sd-OT and QC correction with different distortion thresholds $\varepsilon$.  
We use both the output mesh structure and the angle distortion (denoted by Beltrami coefficient distribution) to analyze the topology preservation property, and use the measure distribution to evaluate the algorithm convergence. 

First of all, we observe that the sd-OT, relaxed sd-OT and t-OT algorithms converge to the specified measure distribution with very light differences. Then, looking at the meshes, the four methods generates similar mapping deformations; for further details, they exhibit different angle distortions and mesh structure, especially at the long merging region between red areas. 
The region encoded to red denotes that its Beltrami coefficient is equal to or greater than 1, which corresponds to the undesired distortions, such as skinny, degenerated or flipped triangles. 
We can see that sd-OT generates almost smooth distortions with light red in the central congested area; sd-OT with original topology generates obvious large distortions, while relaxed sd-OT generates significant distortions. Compared with them, t-OT achieves the most smooth distortions (without red). 
The enlarged mesh structure around the focal vertex indicates whether the topology is preserved or not, as well as the quality of the mapped mesh. 
Compared to the original focal vertex structure (degree 4), sd-OT changes the topology with focal vertex of degree 7; sd-OT with the original topology maintains the original triangular connectivity on the mapped vertex positions, but generates skinny and degenerated triangles; and relaxed-OT generates flipped triangles. In contrast to them, our t-OT maintains the topology while obtaining acceptable triangles.

In the following, we will explain the main procedures of t-OT. 

\vspace{1mm}
\noindent\textbf{Relaxed sd-OT.} Relaxed sd-OT does not include DT adaption, therefore it preserves the topological structure. However, due to the relaxation of the admissible space, the power cell might be empty and significant local distortions can occur. It causes ``bad'' triangles, e.g., the flips in Fig. \ref{fig:t-OT_result_2Area} (c) with the Beltrami coefficients $| \mu(v_i)| > 1$.

\vspace{1mm}
\noindent\textbf{Quasiconformal correction.} We detect the bad vertices with undesired distortions by BC with threshold $\varepsilon$, and follow the strategy of selecting the one-ring of bad vertex as unknown variables and two-ring vertices serve as boundary constraints for computing QC mapping. If the boundary is not convex, we then need to further expand until reaching convex condition. We apply this setting to harmonically diffuse the boundary Beltrami coefficients to fill the interior based on the original geometry of $M$, and then compute the QC mapping on $M$ by the auxiliary metric method to get the positions for the interior vertices. 
It is obvious that the local angle distortions in Fig. \ref{fig:t-OT_result_2Area} (c) are erased in (d-e). 
Meanwhile, the measure distributions of $\omega$ and $|\omega- \nu|$ 
remain consistent due to the faces detected under thresholds $1.0$ and $0.7$ have no much difference. 

\vspace{1mm}
\noindent\textbf{t-OT optimization.} 
By applying the iterative optimization until $\left| \nabla E(\textbf{h}) \right|<\epsilon$ with $\epsilon = 1e-5$ and distortion control $\varepsilon=1.0$ or $0.7$ in Algorithm \ref{alg:1}, t-OT can achieve similar performance on measure convergence compared to the baseline sd-OT and have the final mesh $\hat{M}$ with topology preserved. 
Thus, for a given 3D mesh, starting with its conformal domain, t-OT can work as a mesh parameterization which is able to flexibly control two key components: \textit{area distortion} (customized by measure functions) and \textit{angle distortion} (set by Beltrami coefficient threshold).  

\begin{figure*}[h]
\centering
\footnotesize
 \setlength{\fboxsep}{0pt}
 \setlength{\tabcolsep}{2pt}
 \begin{tabular}{cc} \includegraphics[height=0.23\linewidth]{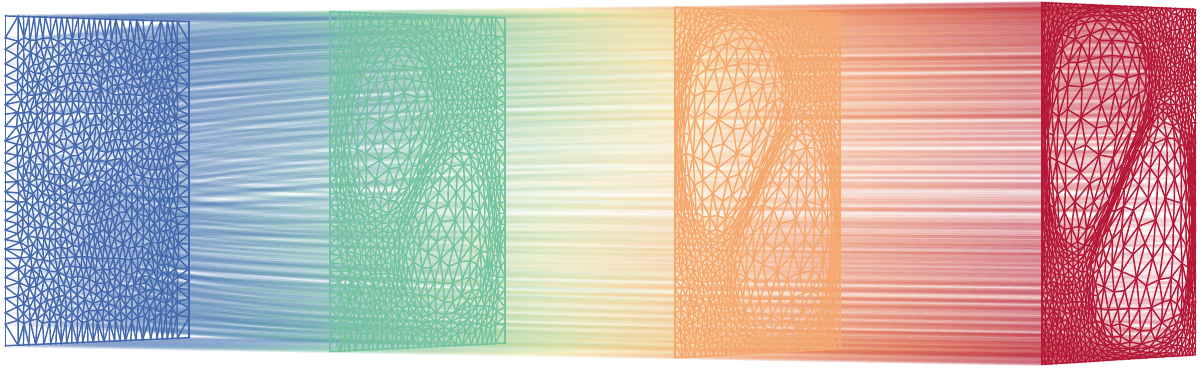}&
\includegraphics[height=0.23\linewidth]{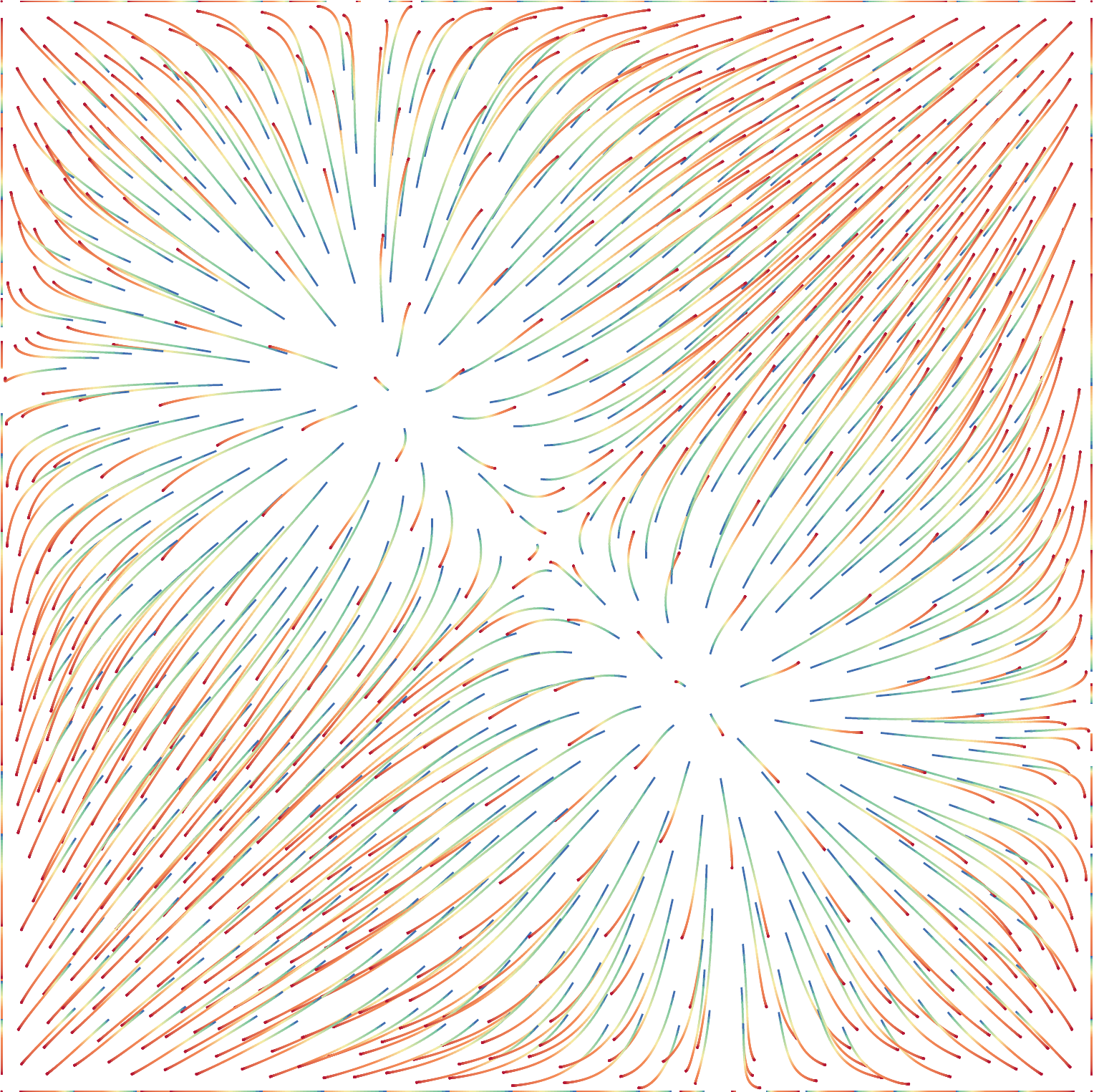}\\
(a) Side View&(b) Front View\\
\end{tabular}
\caption{Trajectory of the temporal t-OT from different perspectives.} 
\label{fig:TQCOT_results_1Color_Mesh}
\end{figure*}

\vspace{1mm}
\noindent\textbf{Temporal t-OT.}
We further apply temporal t-OT to the same input example. 
As illustrated in Fig. \ref{fig:TQCOT_results_1Color_Mesh}, each temporal triangular mesh during the transportation process is assigned a different color. In the side view, each color gradient curve represents the dynamic trajectory of a vertex of the whole transportation. There is no intersection of the trajectories, intuitively representing the topology preservation property. 
The front view vividly visualizes the expansion movement of two red regions of the original input. Although the dynamic trajectory curves visually intersect, the colors at the intersection points are different denoting that they reach the intersecting location at different times. Both of them demonstrate the tt-OT has no intersecting trajectories for all vertices (particles in dynamics) in the whole process. The tt-OT result at any intermediate moment has topology preserved.

\subsection{Analysis}
\label{sec:experiment}
We analyze the proposed algorithm in terms of topology preservation, convergence, and computational efficiency.

\vspace{1mm}
\noindent\textbf{Topology guarantee.} 
As a discretization method, DT helps reduce the discretization interpolation error of the Brenier potential function $u$ , thereby ensuring the accuracy of the geometric duality relationships in the algorithm. 
The baseline sd-OT method requires DT throughout the optimization process
\citep{gu2013variational,helfer2013secondary,lei2019secondary}. 
In our t-OT method, the relaxed sd-OT procedure avoids DT operations and maintains vertex connectivity, and the QC correction ensures a smooth, flip-free diffeomorphism, preserving the mesh topological structure during the process.

\begin{figure*}[t]
\centering
\footnotesize
\begin{tabular}{cc}
\includegraphics[width=0.45\textwidth]{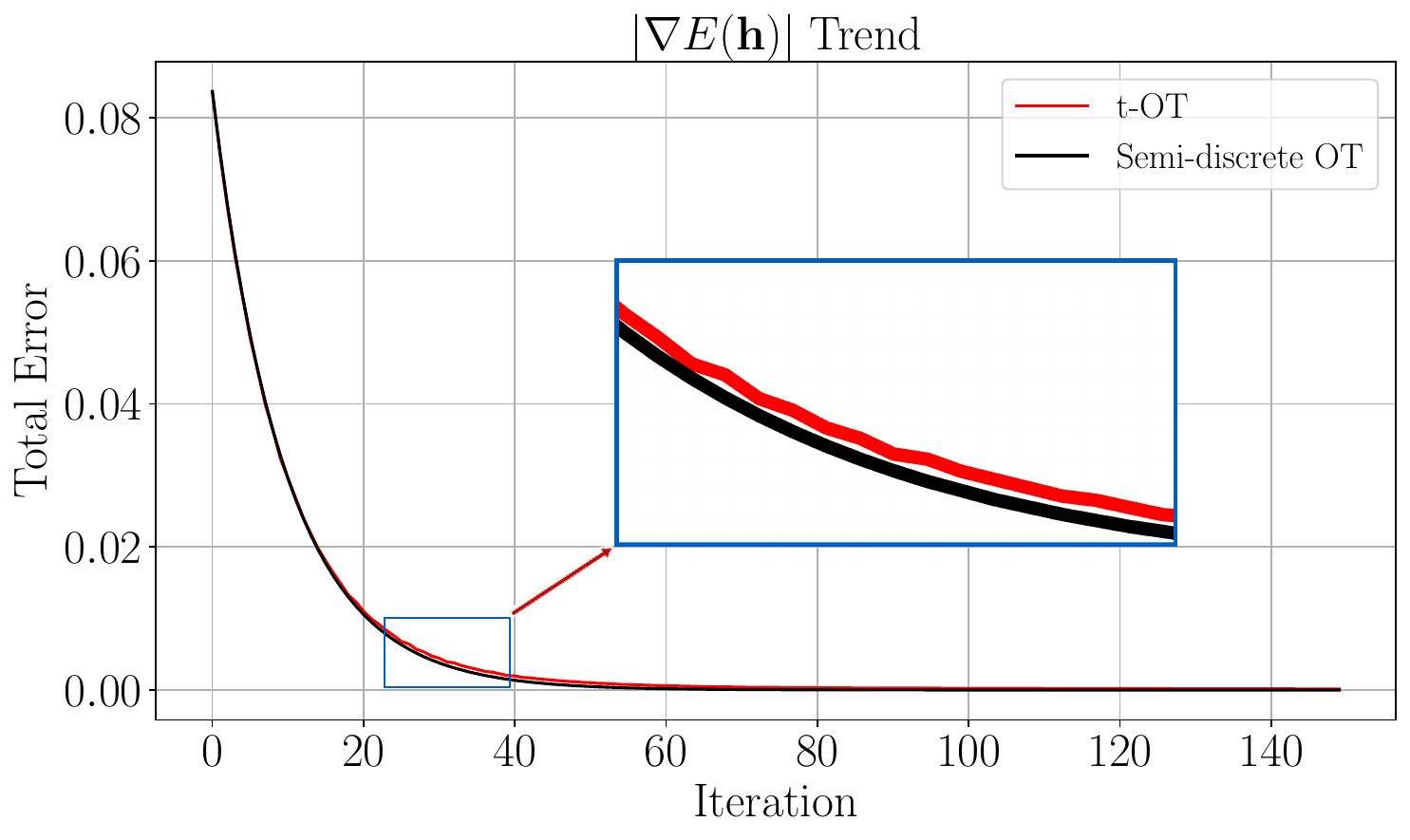}&
\includegraphics[width=0.45\textwidth]{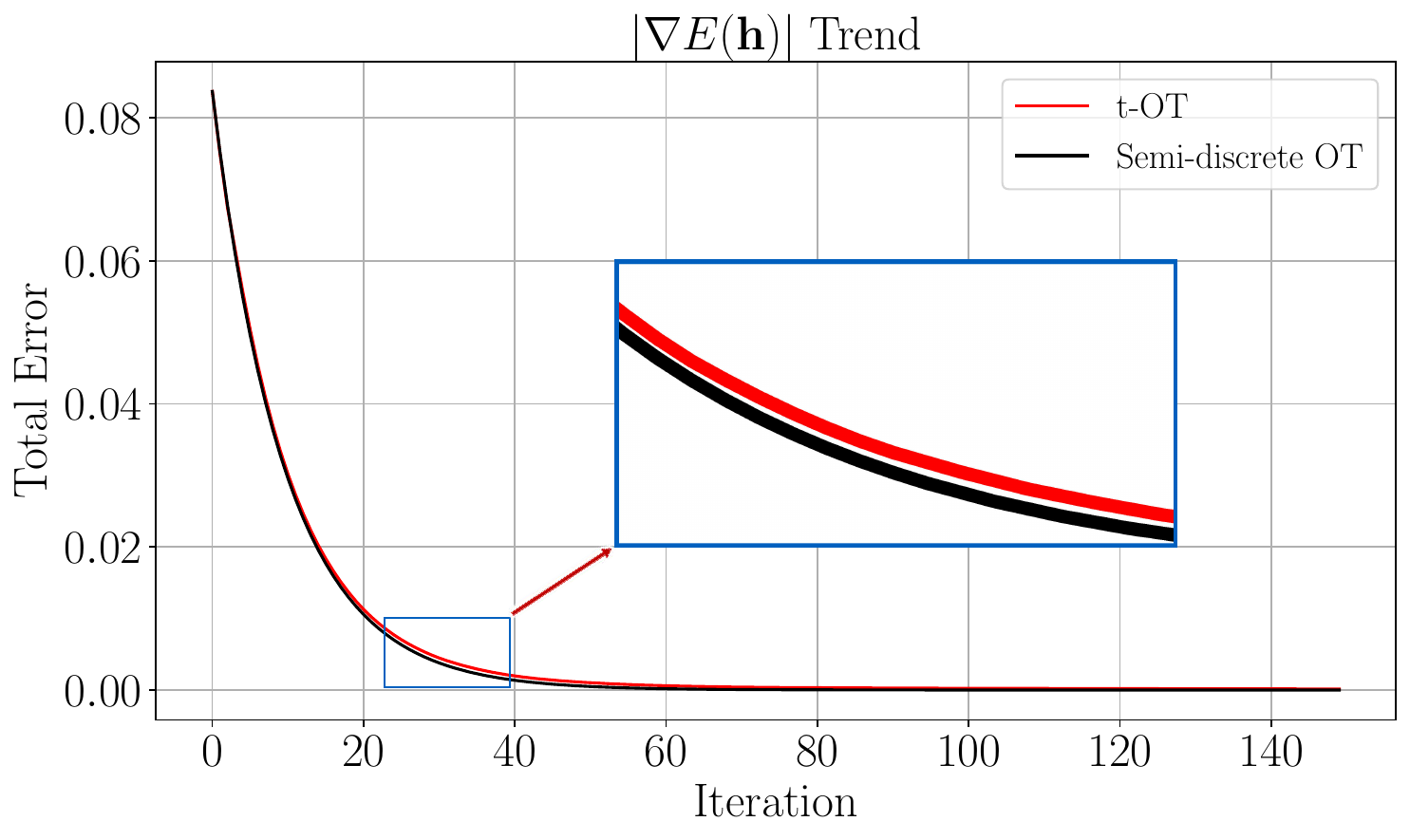}\\
\includegraphics[width=0.45\textwidth]{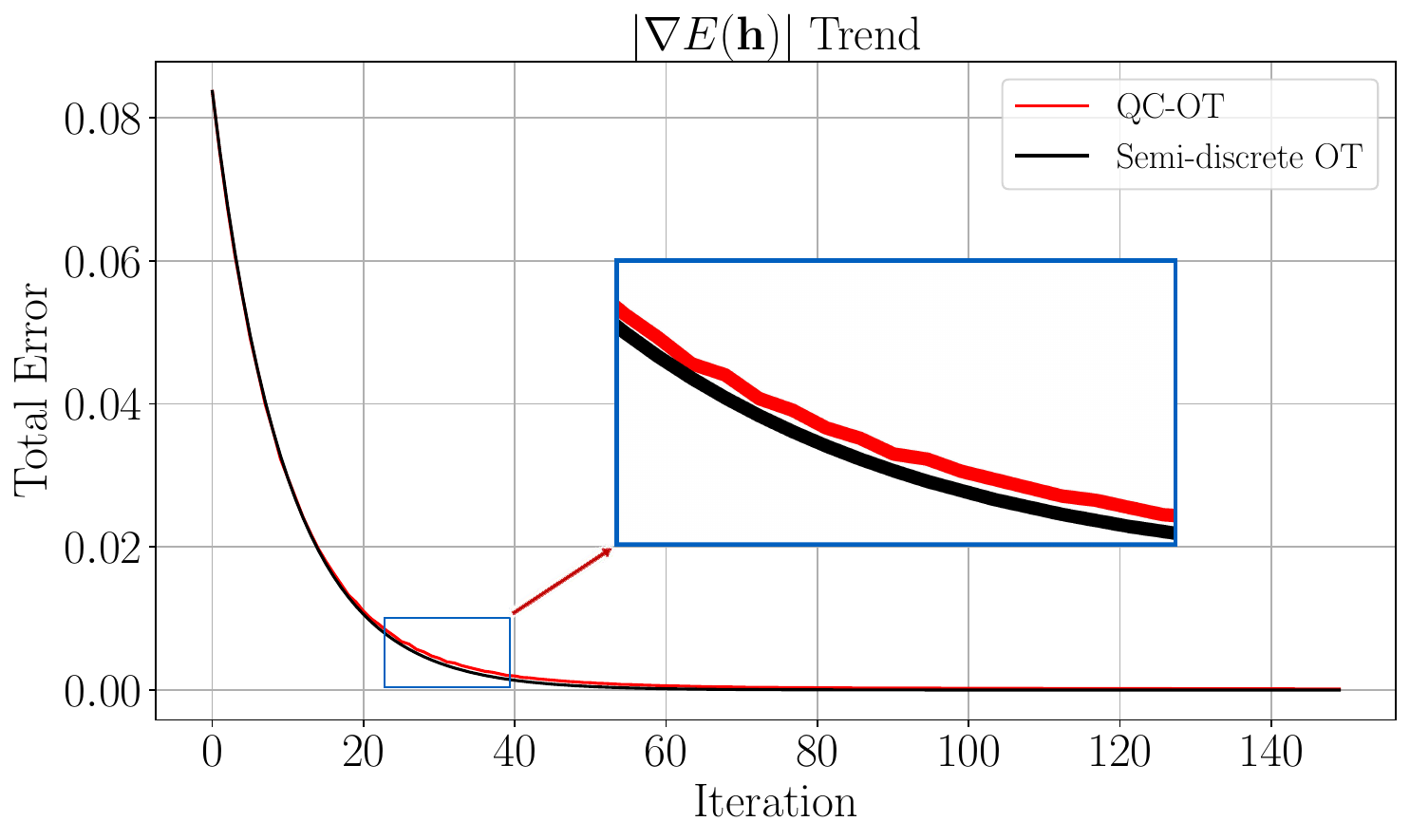}&
\includegraphics[width=0.45\textwidth]{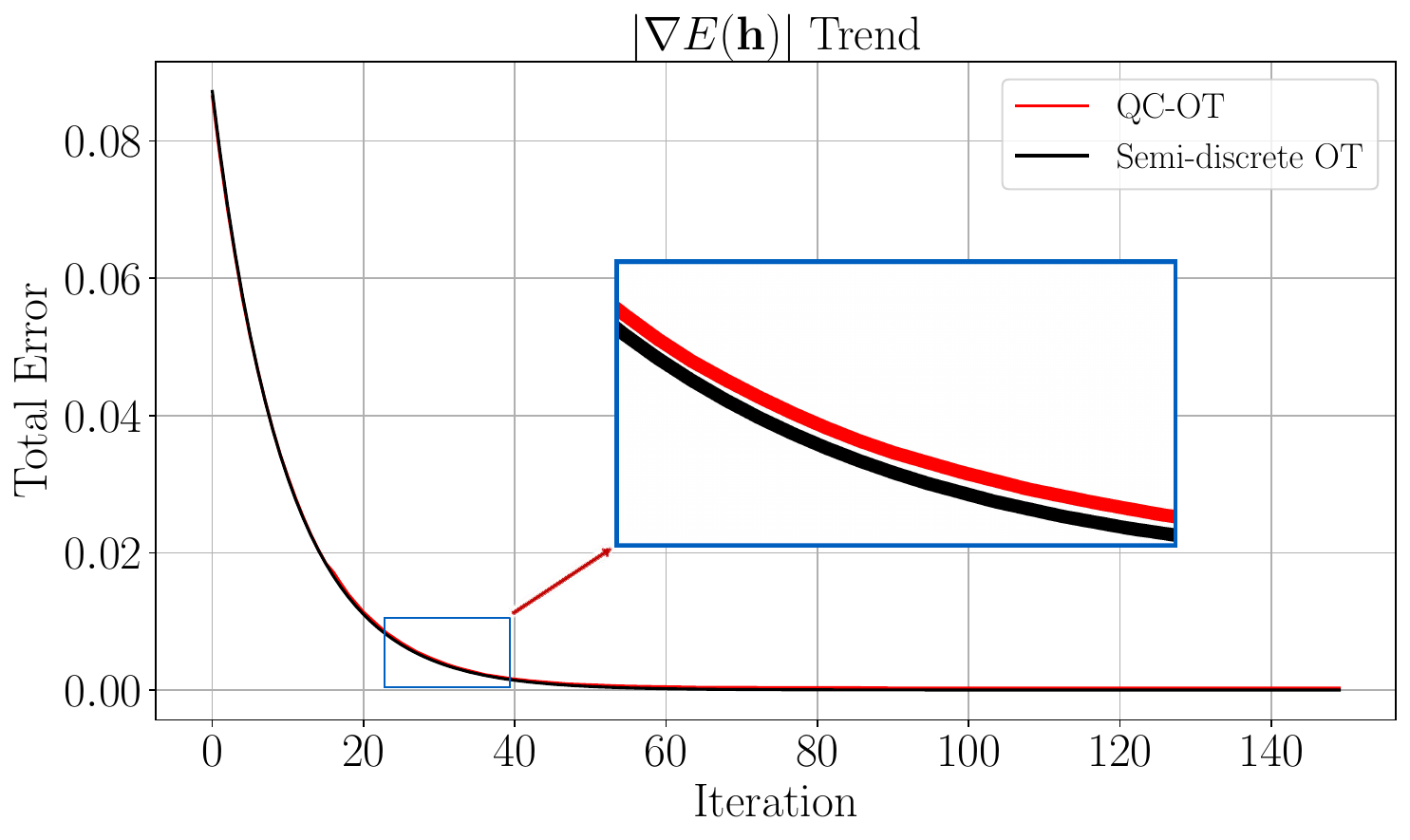}\\
(a)  $\lambda \leftarrow \frac{1}{2} \lambda$: empty power cell check &(b) $\lambda \leftarrow \frac{1}{2} \lambda$: $\widetilde{\Delta} \left(\left| \nabla E \right| \right) < 0$\\
\end{tabular}
\captionsetup{justification=justified, singlelinecheck=false}
\caption{Comparison between t-OT and sd-OT on the convergence trends of $|\nabla E(\textbf{h}) |$ of Fig. \ref{fig:t-OT_result_2Area}(e). Different update strategies of step length $\lambda$ are applied to the t-OT framework.} 
\label{EnergyTrendCicle}
\end{figure*}

\vspace{1mm}
\noindent\textbf{Convergence.}\label{CA}
We use the distribution of the value $\rho_{i} = |\omega_{i} - \nu_{i}|$  as the metric for evaluation of the convergence of the algorithms of t-OT and sd-OT. When $\rho_{i}=0$, the obtained measure reaches prescribed measure. 
Figure \ref{fig:t-OT_result_2Area} illustrates the distribution of $\omega$ and $\rho$ and it is clear that $\rho$ clusters around $0$ in t-OT and sd-OT, showing progressive convergence.
Figure \ref{EnergyTrendCicle} shows the convergence trends by $\left| \nabla E(\textbf{h}) \right|$ for both algorithms, where the convergence rate of t-OT is slightly lower than sd-OT.

In sd-OT, the strict convexity of the Brenier potential function is equivalent to the non-emptiness of Voronoi power cells during the computation process. This ensures the existence and uniqueness of the solution, as well as the favorable convergence property \citep{lei2021optimal}. In t-OT, however, the relaxation of the admissible space may result in empty power cells, causing the Brenier potential function to lose its strict convexity. This can potentially lead to the issues such as slower convergence, numerical instability, or oscillatory behavior. 
We adjust the update strategy of the height parameter as a complementary means to enhance the convergence performance. Unlike the update strategy of step length $\lambda$ based on checking for empty power cells in sd-OT, we update $\lambda \leftarrow \frac{1}{2} \lambda$ iteratively until the gradient variation of $\left| \nabla E \right|$ satisfies $\widetilde{\Delta} \left(\left| \nabla E \right| \right) < 0$.
Then we compare the two update strategies in t-OT framework. As shown in Fig. \ref{EnergyTrendCicle}, it is obvious that the original strategy may result in oscillation while the new strategy generates a smooth convergence trend. 
Therefore, by applying this auxiliary approach, our proposed algorithm can  exhibit an improved convergence performance.

\vspace{1mm}
\noindent\textbf{Running time.} 
Table \ref{tab:convergence-compa-all} shows the comparison of the  computational efficiency among t-OT, tt-OT and sd-OT, where two tests were conducted under the same settings of the predefined density functions and the threshold parameter $\epsilon=1e-5$. 
The experiments were performed using generic C++ on a Windows 10 64-bit platform with Intel 3.7GHz CPU and 64 GB of RAM.

\noindent\textit{Test 1: Running with the same iterations.}  It is obvious that t-OT has shorter running time and therefore is more efficient than sd-OT in each iteration.

\noindent\textit{Test 2: Running to the end.} Although the methods converge at different iteration numbers, t-OT achieves a higher convergence speed and tt-OT is slightly slower than t-OT with the same iteration number. 

\noindent The testing results are consistent to the algorithm design. t-OT relaxes the DT and convexity check operations, therefore, requiring less time; tt-OT needs QC correction in each iteration, therefore, slower than t-OT. 

With the above analysis, we conclude that t-OT approximates sd-OT, but preserves mesh topology and has superior computational efficiency. 

\begin{table}[h]
\caption{Comparison of running time (s) among sd-OT, t-OT and tt-OT.}\label{tab:convergence-compa-all}
\begin{tabular*}{\textwidth}{@{\extracolsep\fill}lccccccccc}
\toprule%
\multirow{2}{*}{} & \multicolumn{4}{@{}c@{}}{Test 1: Running with same iterations} & \multicolumn{5}{@{}c@{}}{Test 2: Running to the end} \\
\cmidrule{2-5}\cmidrule{6-10}%
 Mesh& Iter & sd-OT &t-OT &tt-OT& Iter &sd-OT&Iter&t-OT&tt-OT\\
\midrule
Fig. \ref{t-OT_input_2Area} (a)&25 & 0.365 &\textbf{0.108}&0.126& 61 & 0.764 & 112& \textbf{0.499}& 0.584 \\ 
Fig. \ref{fig:faces-conformal} (c) & 50 & 0.909& \textbf{0.319} &0.324& 80 & 1.42&79& \textbf{0.43}& 0.507 \\
Fig. \ref{Medical} (a)& 100 & 31.366 & \textbf{17.739} &18.287& 71 & 22.901&74 &\textbf{13.13}& 13.69\\
\bottomrule
\end{tabular*}
\end{table}

\subsection{Discussion}
We further discuss some intuitive questions on possible alternatives of the method and the resulting properties. 

\begin{figure*}[t]
\centering
\footnotesize
 \setlength{\fboxsep}{0pt}
\begin{tabular}{@{\hspace{1em}}m{0.22\linewidth}<{\centering}@{\hspace{1em}}m{0.22\linewidth}<{\centering}@{\hspace{1em}}m{0.22\linewidth}<{\centering}@{}}
\fbox{\includegraphics[width=0.19\textwidth]{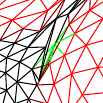}} &
\fbox{\includegraphics[width=0.19\textwidth]{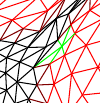}}&
\fbox{\includegraphics[width=0.19\textwidth]{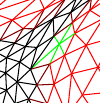}}\\
(a) Original Distortion&(b) Harmonic Filling&(c) QC Filling \\
	\end{tabular}
	\caption{Comparison of filling alternatives. 
 }
\label{fig:LocalTwistFilling}
\end{figure*}

\vspace{1mm}
\noindent\textbf{Harmonic map filling.} 
Note that the Beltrami coefficients capture the deformations from the original mesh up to the current position in the overall relaxed sd-OT process. Thus, the Beltrami coefficients for the unknown area, obtained by diffusion from boundary Beltrami coefficients, ensure that the deformations remain smooth on the original mesh. The QC mapping, implemented as a harmonic map with the auxiliary metric derived from the diffused Beltrami coefficients, preserves the deformations. If we directly diffuse the boundary vertex positions to the unknown area using the original mesh metric, we would lose the OT deformations for that area. 
Figure \ref{fig:LocalTwistFilling} illustrates an example on mesh structures of the twisted region under the methods of harmonic filling and QC filling with distortion control threshold $\varepsilon = 0.7$. 
While both methods eliminate local distortion, QC filling produces higher-quality triangular faces. We further test the two methods from the perspective of probability measures, using the metric of the measure ratio of the highlighted area over the total. As shown in Table \ref{tab:mu_harmonic}, QC filling preserves the OT measure better than harmonic filling in mesh correction. Therefore, QC filling is selected for our algorithm to 
enhance the stability and validity of the computation.

\begin{table}[t]
\captionsetup{width=1\linewidth}
\captionsetup{justification=justified, singlelinecheck=false}
\caption{Comparison of mesh correction methods by evaluating measure ratio.
}
\label{tab:mu_harmonic}

\begin{tabular}{@{}lccc@{}}
\toprule
Methods &Harmonic filling&QC filling&(prescribed measure)\\
\midrule
Measure ratio & 0.00384& \textbf{0.00527}& 0.00556 \\
\bottomrule
\end{tabular}
\end{table}

\vspace{1mm}
\noindent\textbf{Topology transfer to sd-OT.}
In sd-OT, based on the definition of the lower convex hull \(\pi_{i}^{*} := (p_{i},-h_{i})\), 
the potential adjacency relationships of both lower convex hull and power diagram remain \textit{largely} unchanged. As the power cell is not allowed to be empty, the adjacent power diagram undergoes a similar degree of change. With the influence of DT, the power cell has a relatively regular shape and the updated mass center of it won't go out of range. 
Transferring the original topology to the sd-OT result leads to a topology-preserving mapping $\check T:\check M \to  M_0$, but might generate significant angle distortions in $\check M$. 
This can lead to skinny and nearly degenerative triangles, as illustrated in Fig. \ref{fig:t-OT_result_2Area} (b). 
Such a mesh structure will cause instability in future geometric computation tasks, e.g., when taking the resulted mapping as a mesh parameterization. 
Furthermore, due to the correlation between the computation of discrete probability measure and topological structure, directly changing the topological structure to the original connectivity leads to inaccurate and unpredictable changes on probability measures of vertices. Therefore, directly copying the original topological structure to the final sd-OT map is not desired.

\vspace{1mm}
\noindent\textbf{Intrinsic geometry.} 
The proposed t-OT framework is based on the geometric variational sd-OT baseline and works for simply-connected domains. The resulted mesh of t-OT has relatively regular shape. When handling 3D surfaces, we first compute the conformal parameter domain and then use this domain for computing relaxed sd-OT with a prescribed measure, i.e., adapting the intrinsic geometric structure under the optimal transport conditions. 
The QC mapping is then incorporated to correct distortions that stray from the intrinsic geometry. Thus, the proposed strategy respects the intrinsic geometric structure as much as possible, while preserving the total probability measures and minimizing transportation cost.

\vspace{1mm}
\noindent\textbf{Controllable deformation.}
Beltrami coefficients represent local angle distortions (conformality) of the mapping, 
but can not directly control global scaling deformations, such as isotropic expansion or shrinkage where the discrete BC is close to zero. 
t-OT can intuitively customize mapping deformations by setting density functions (representing area distortions), such as applying higher density to regions that need to be magnified, and specifying the distortion threshold by the magnitude $|\mu|$ (representing angle distortions) to have control on the shape quality of triangular faces. Therefore, t-OT provides more flexibility on designing mesh parameterization and deformation editing. 

\section{Application}
\label{sec:experiment}
We apply the proposed t-OT and tt-OT framework to surface mesh parameterization, image mesh editing, and physical diffusion simulation. There are various strategies to set probability measures,  
$g(p_{i}), p_i\in V$, as follows:

\begin{itemize}
\item An area function of 3D triangular mesh to generate an area-preserving parameterization, which can be regarded as assigning a probably different scalar to every vertex such that each local neighborhood around vertex is deformed isotropically (see Fig. \ref{fig:area-preserve-BW}). Here, $g(p_{i})= A_s(p_i), p_{i} \in V$, where $A_s$ denotes the normalized  area of vertex power cell on surface. 

\item A uniform value to all vertices to equalize mesh density and generate a mesh-equalizing parameterization, which globally depends on topological connectivity but locally still respects the original geometry  (see Fig. \ref{fig:AreaEqualFaces_BW}). Here, $g(p_{i})= 1/N,  p_{i} \in V$, where $N=|V|$.  

\item A scalar function on regions of interests (ROIs), to generate corresponding local conformal (isotropic) expansion or shrinkage with global deformations in the parameterization 
(see Fig. \ref{CustomizedFaces}). Here,  
$g(p_{i})=  k_i * a_{i}, p_{i} \in  \hat{O}$, where $k_i$ is a constant scalar or a scalar function varying on vertices, $a_i$ denotes the original power cell area of vertex $p_i$, and  $\hat{O}$ denotes the vertex set in ROIs.

\item A density function defined by the given grayscale image. Here, $g(p_{i}) = k*a_i*(rgb + \delta), p_{i} \in V$, where $k$ is a constant scalar, $a_i$ denotes the original power cell area of vertex $p_i$, $rgb$ indicates the grayscale intensity, and $\delta$ is the perturbing term.

\end{itemize}

In all experiments below, we set the energy tolerance to $\epsilon = 1e-5$ and the distortion tolerance to $\varepsilon = 0.7$. 

\subsection{Surface Mesh Parameterization}\label{sec:5.1}
For surface cases, we display the distribution of the resulting measure $\omega$ in the t-OT results and use its difference with the specified metric $|\omega-\nu|$ to evaluate the convergence of the algorithm. All experimental cases verify that t-OT can achieve the prescribed measure $\nu$ with high accuracy. At the same time, we use check-board texture mapping from the t-OT parameter domain back to surface to visualize the deformations generated by t-OT. In the following, we list three types of t-OT mesh parameterization corresponding to the first three probability measure setting strategies.

\begin{figure*}[h]
\centering
\footnotesize
\setlength{\tabcolsep}{-1pt} 
\begin{tabular}{@{\hspace{-0.5em}}m{0.01\linewidth}<{\centering}@{\hspace{-1em}}m{0.23\linewidth}<{\centering}@{\hspace{-0.9em}}m{0.23\linewidth}
<{\centering}@{\hspace{0.6em}}m{0.23\linewidth}
<{\centering}@{\hspace{0.7em}}m{0.3\linewidth}<{\centering}@{}}
\rotatebox{90}{Dense Mesh}&
\includegraphics[height=0.23\textwidth]{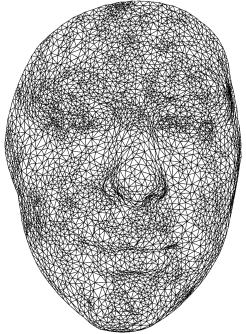}&
\includegraphics[height=0.23\textwidth]{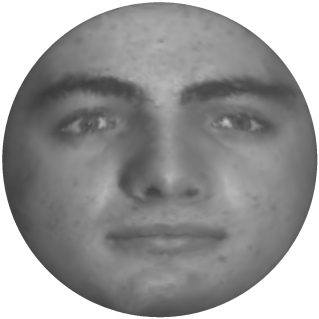}&
\includegraphics[height=0.23\textwidth]{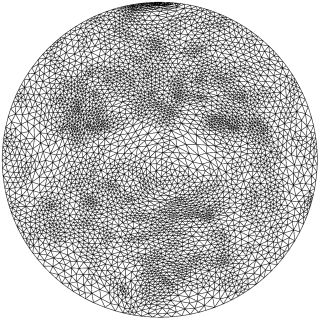}&
\includegraphics[height=0.23\textwidth]{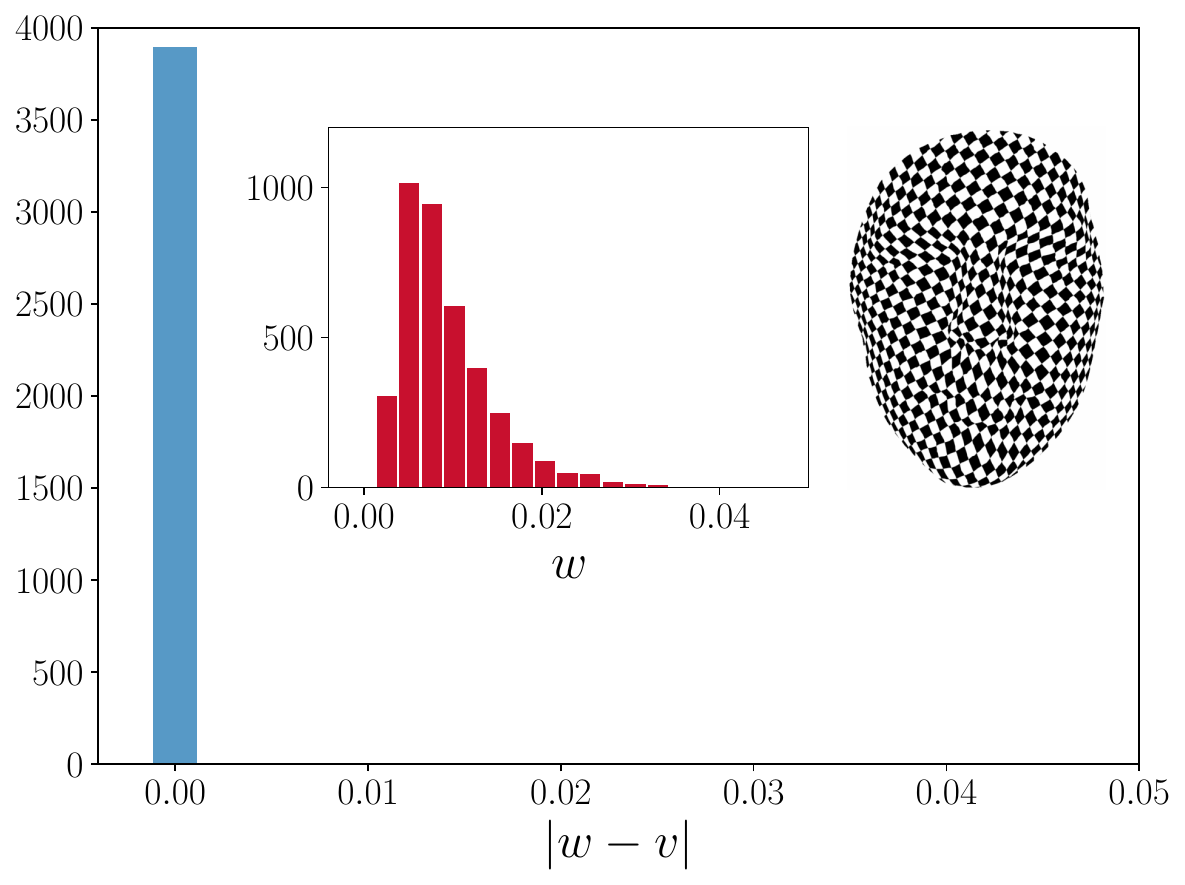}\\   
\rotatebox{90}{Coarse Mesh}&
\includegraphics[height=0.23\textwidth]{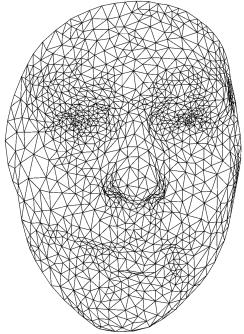}&
\includegraphics[height=0.23\textwidth]{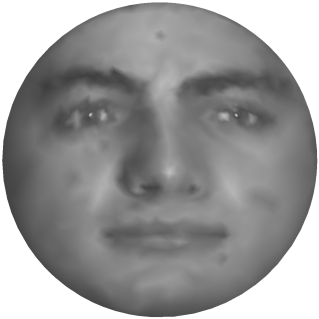}&
\includegraphics[height=0.23\textwidth]{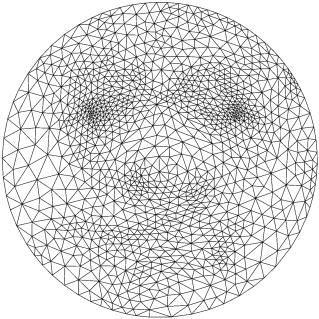}&
\includegraphics[height=0.23\textwidth]{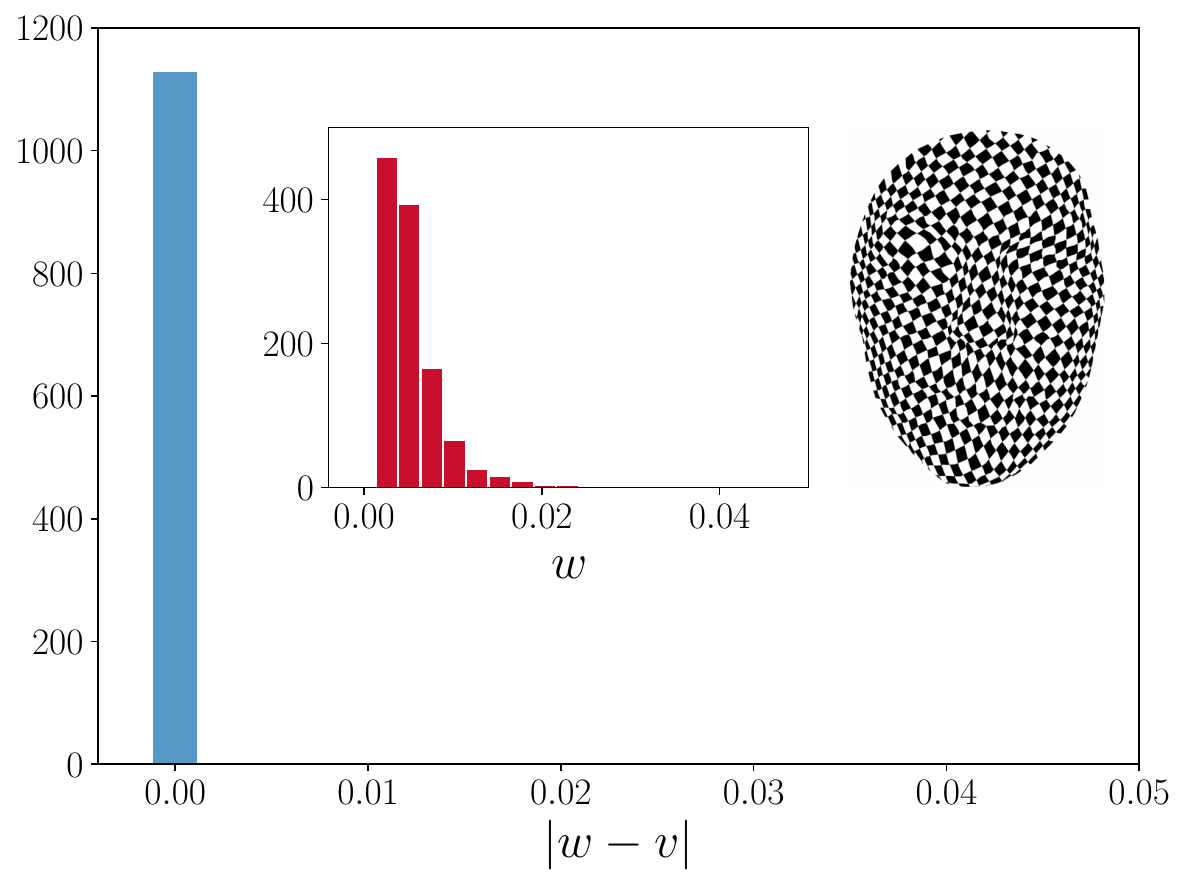}\\ 
& (a) 3D Surface&(b) Parameter Domain&  (c) Parameter Mesh&(d) Measure Distribution\\
\end{tabular}
\caption{Area-preserving mesh parameterization for  coarse and dense meshes.
} 
\label{fig:area-preserve-BW}
\end{figure*}

\vspace{1mm}
\noindent\textbf{Area-Preserving Parameterization.} 
We consider the local area of 3D surface to drive the mapping and design the probability measure in t-OT. 
Figure \ref{fig:area-preserve-BW} shows the area-preserving parameterizations by t-OT for the human face surface C in Fig. \ref{fig:faces-conformal} which exhibits the original 3D surface, 2D parameter domain, 2D parameter mesh, and checker-board texture mapping results of face C. For the same human face shown in (a), we use two mesh resolutions. The numbers of vertices and triangular faces $\left[\#v,\#f\right]$ of the coarse and dense meshes are $\left[ 1129, 2182 \right] $ and $\left[ 3896, 7628 \right] $, respectively. We first generate their intrinsic geometric representations by mapping them conformally to the 2D disk domain. The area-preserving t-OT, driven by 3D area probability measure, then produces very similar results which can be clearly visualized by the parameter domains with face texture in (b). That is due to the fact that they both respect the original area of the 3D surface, which is consistent in terms of both coarse and dense structure.

\begin{figure*}[h]
\centering
\footnotesize
\begin{tabular}{@{\hspace{-0.5em}}m{0.01\linewidth}<{\centering}@{\hspace{-1em}}m{0.23\linewidth}<{\centering}@{\hspace{-1em}}m{0.23\linewidth}
<{\centering}@{\hspace{0.7em}}m{0.23\linewidth}
<{\centering}@{\hspace{0.7em}}m{0.3\linewidth}<{\centering}@{}}
&
\includegraphics[height=0.23\textwidth]{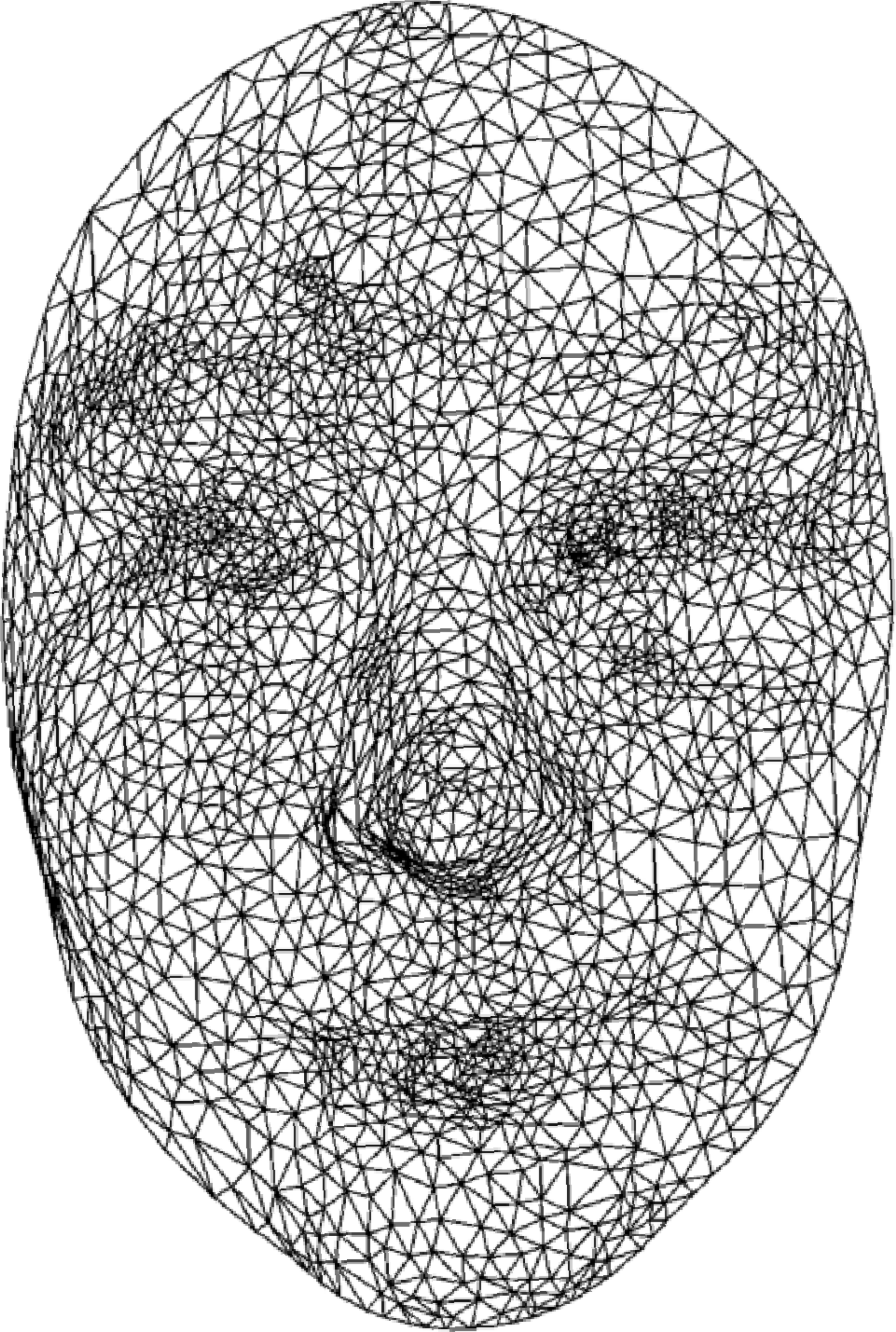}&
\includegraphics[height=0.23\textwidth]{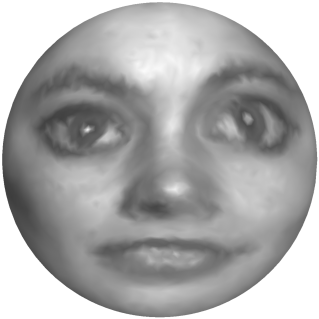} &\includegraphics[height=0.23\textwidth]{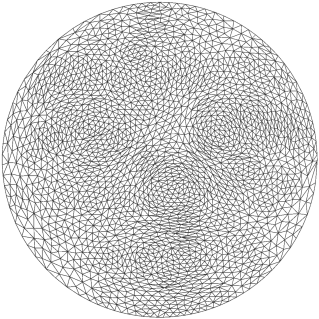}
&\includegraphics[height=0.23\textwidth]{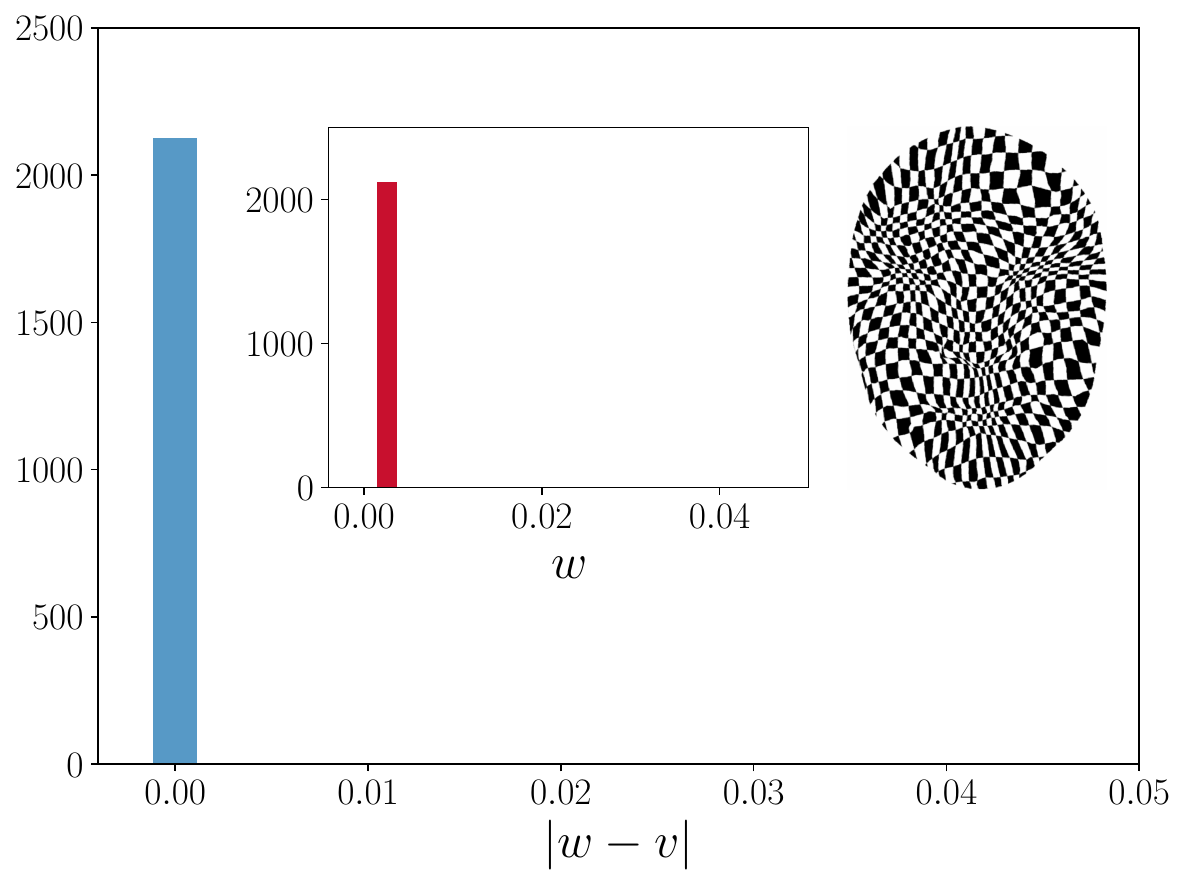}\\

&
\includegraphics[height=0.23\textwidth]{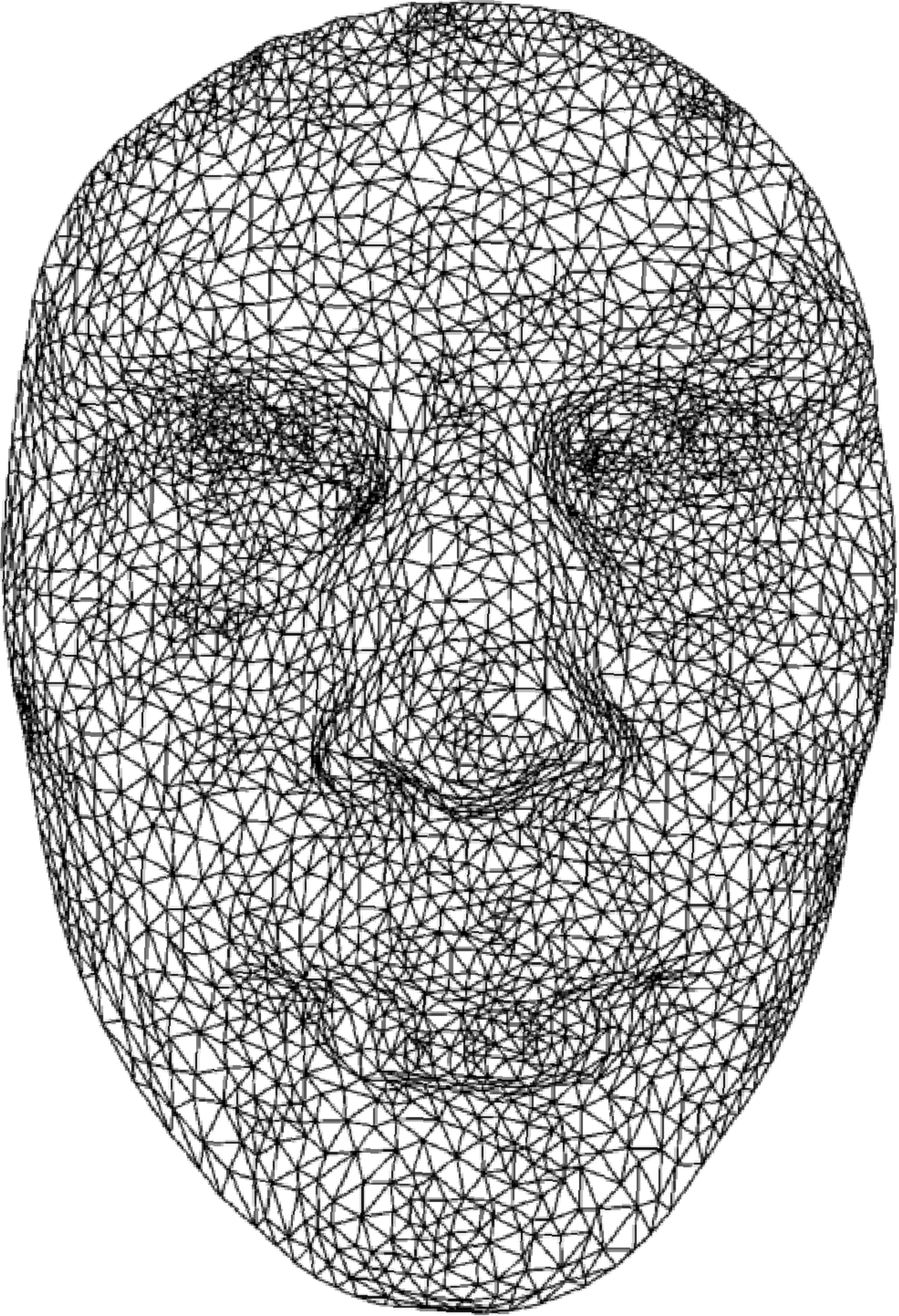}&	
\includegraphics[height=0.23\textwidth]{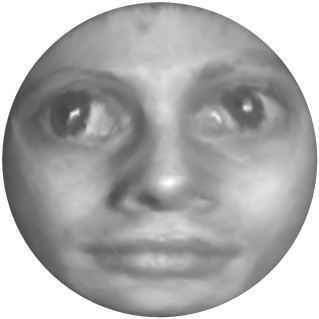}&
 \includegraphics[height=0.23\textwidth]{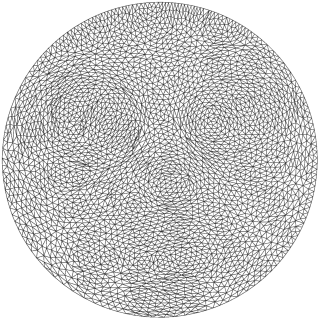}&
 \includegraphics[height=0.23\textwidth]{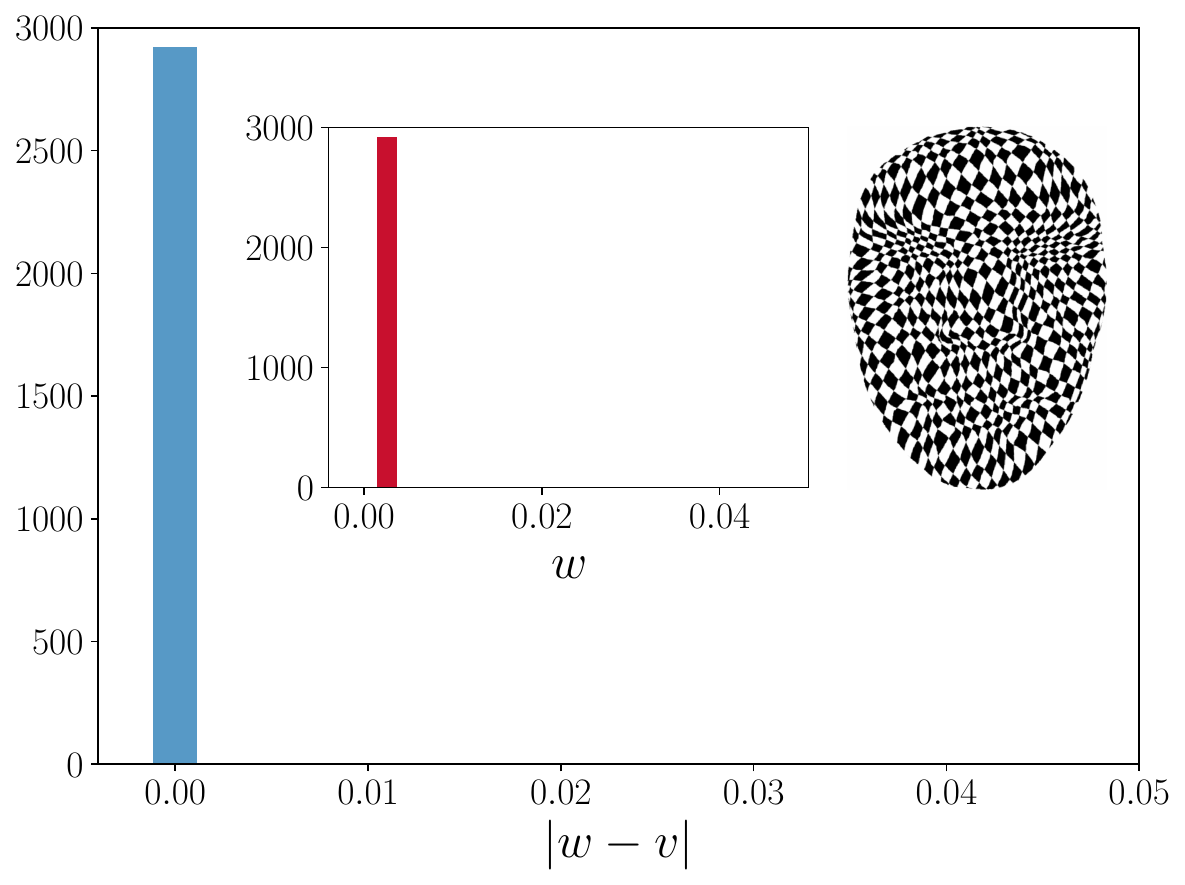}\\
& (a) 3D Surface&(b) Equalizing Domain&  (c) Parameter Mesh&(d) Measure Distribution\\
\end{tabular}
\caption{Mesh-equalizing parameterization for 3D human faces.
 } 
\label{fig:AreaEqualFaces_BW}
\end{figure*}

\begin{figure*}[h]
\centering
\footnotesize
\begin{tabular}{@{\hspace{0.5em}}m{0.01\linewidth}<{\centering}@{\hspace{0.7em}}m{0.23\linewidth}<{\centering}@{\hspace{1em}}m{0.23\linewidth}
<{\centering}@{\hspace{1em}}m{0.3\linewidth}<{\centering}@{}}
\rotatebox{90}{Dense Mesh}

&
\includegraphics[height=0.23\textwidth]{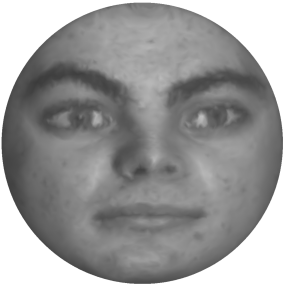}&
\includegraphics[height=0.23\textwidth]{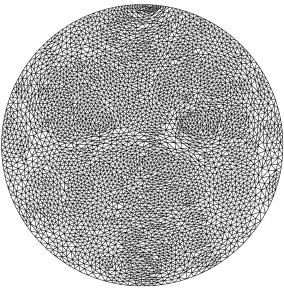}&
 \includegraphics[height=0.23\textwidth]{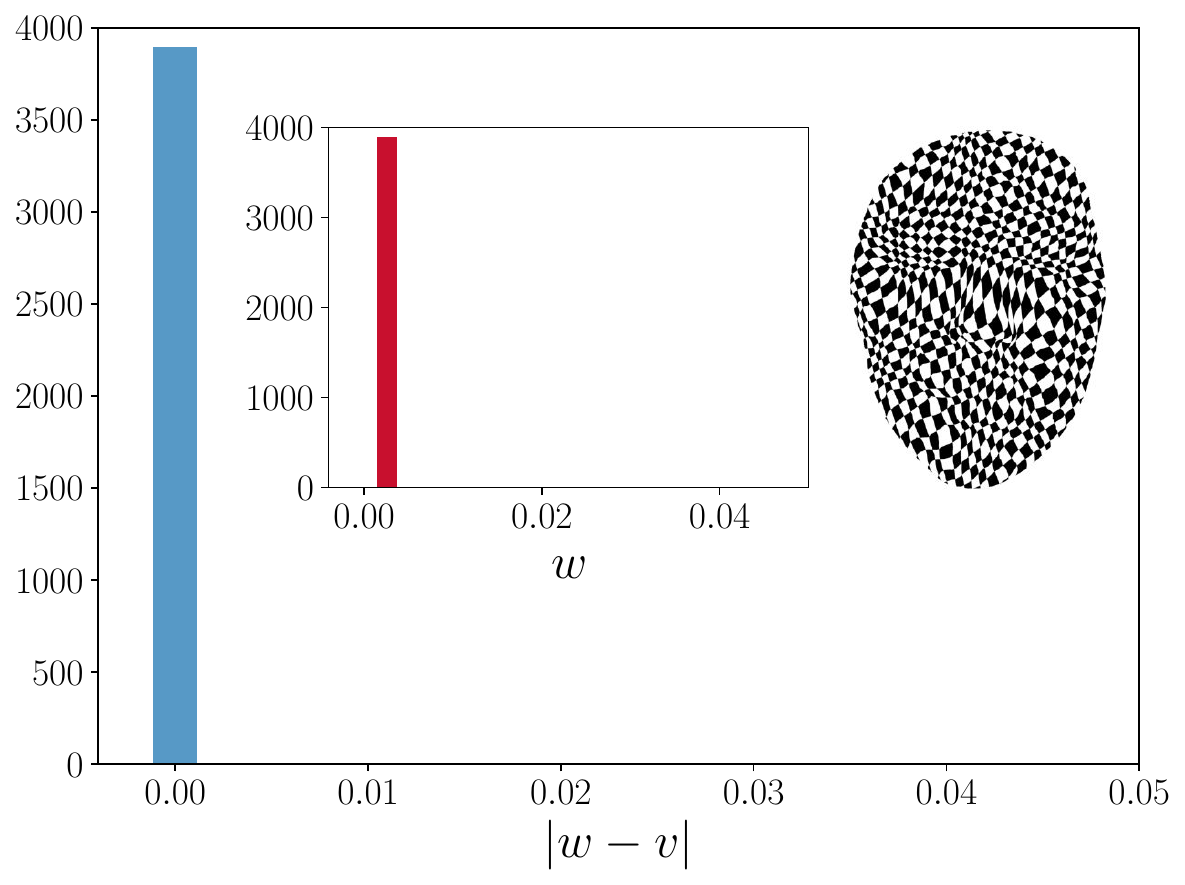}\\
\\

\rotatebox{90}{Coarse Mesh}
&
\includegraphics[height=0.23\textwidth]{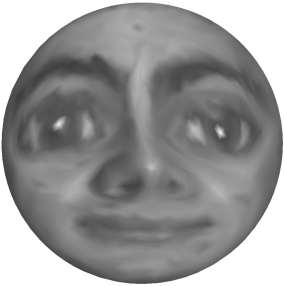}&
\includegraphics[height=0.23\textwidth]{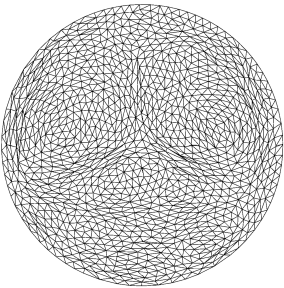}&
 \includegraphics[height=0.23\textwidth]{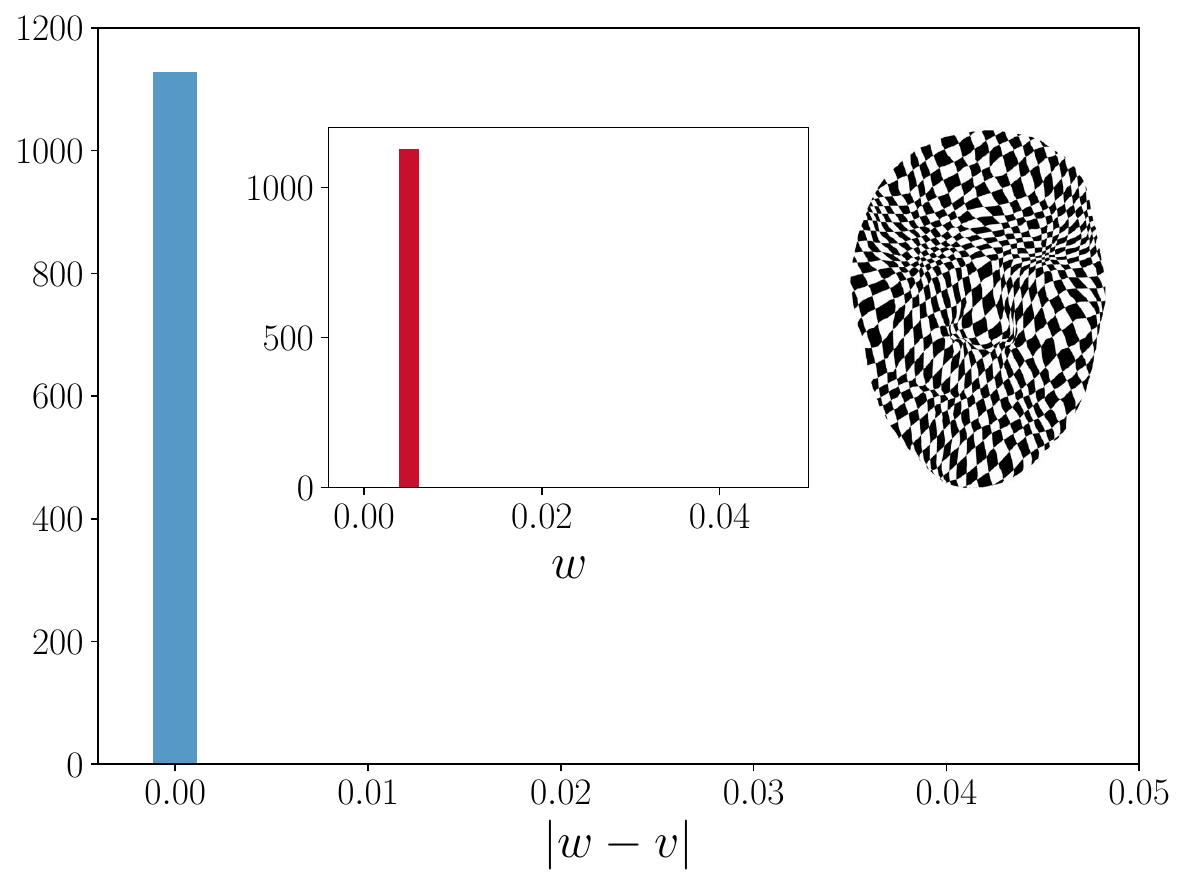}\\
 &(a) Equalizing Domain&  (b) Parameter Mesh&(c) Measure Distribution\\
\end{tabular}
\caption{Mesh-equalizing parameterization for 3D human faces.
 } 
\label{fig:AreaEqualFaces_DavidOnly}
\end{figure*}

\vspace{1mm}
\noindent\textbf{Mesh-Equalizing Parameterization.}
We apply the same probability measure to all vertices in t-OT to achieve mesh parameterization with an equalized mesh structure.
In this experiment, our mesh sizes are face A $\left[ 2127, 4146 \right]$ and face B $\left[ 2925, 5687 \right]$. 
As the density function in $\Omega$ is uniform, therefore, the measure of each vertex $\tau_{i}$ means the area of its power cell. 
Guided by the prescribed measure $\nu$, each mesh component is uniformly adjusted. Figure \ref{fig:AreaEqualFaces_BW} displays the mesh-equalizing parameterization results, demonstrating an improved uniformity compared to the original conformal parameterization.
The equalizing parameterization results for face C, as shown in Fig. \ref{fig:AreaEqualFaces_DavidOnly}, vary significantly across coarse and dense mesh structures due to the different distribution of vertices. 

\begin{figure*}[h]
\centering
\footnotesize
\begin{tabular}{@{\hspace{0em}}m{0.01\linewidth}<{\centering}@{\hspace{0.7em}}m{0.23\linewidth}<{\centering}@{\hspace{1em}}m{0.23\linewidth}
<{\centering}@{\hspace{1em}}m{0.3\linewidth}<{\centering}@{}}
&
\includegraphics[height=0.23\textwidth]{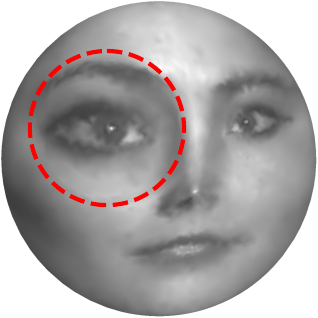} &
\includegraphics[height=0.23\textwidth]{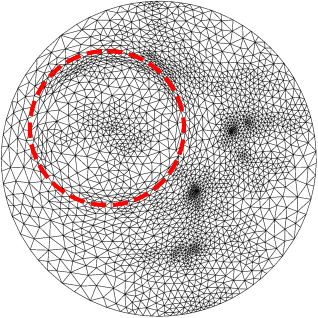} &
\includegraphics[height=0.23\textwidth]{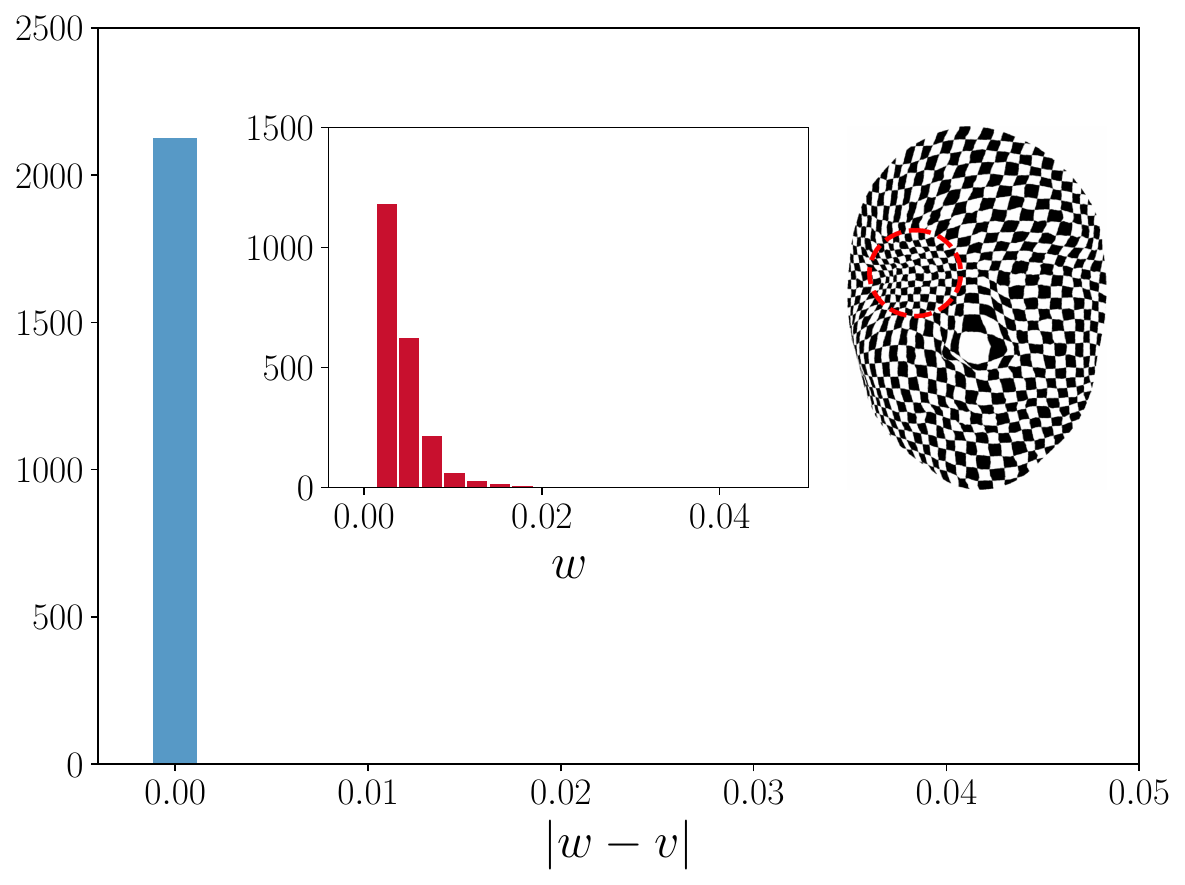} \\

&
\includegraphics[height=0.23\textwidth]{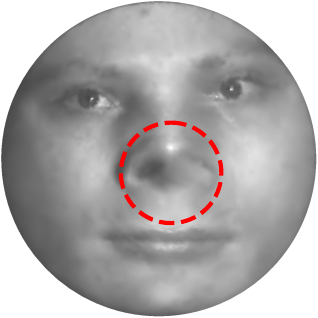}&
\includegraphics[height=0.23\textwidth]{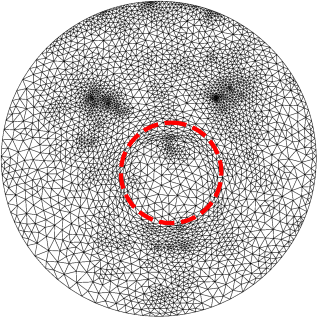}&
\includegraphics[height=0.23\textwidth]{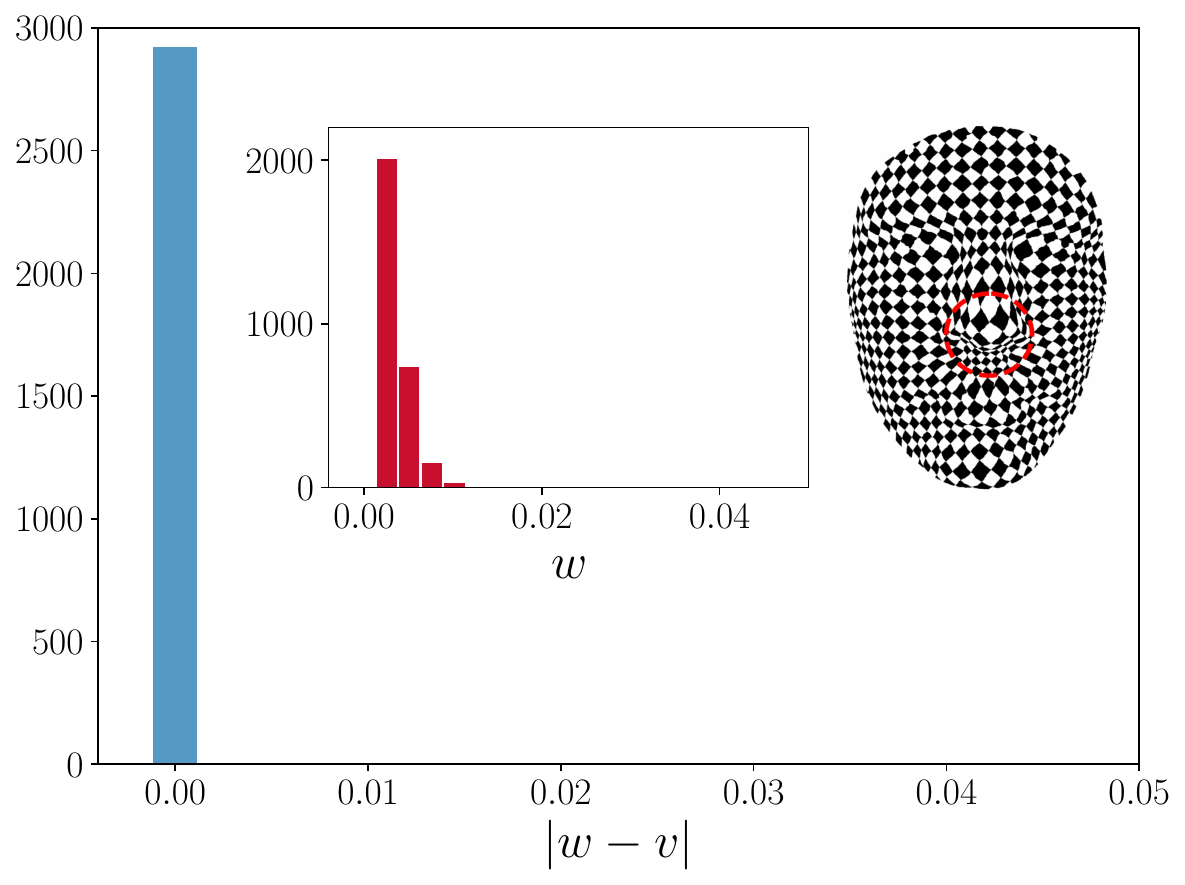}\\ 
&(a) Parameter Domain&(b) Parameter Mesh&(c) Measure Distribution\\
\end{tabular}
\caption{Customized mesh parameterization for 3D meshes in Fig. \ref{fig:faces-conformal} by t-OT algorithm.
   } 
\label{CustomizedFaces}
 \end{figure*}

\vspace{1mm}
\noindent\textbf{Customized Mesh Parameterization.}
We use the density function as probability measure in the t-OT method to customize deformation, which is an intuitive way to specify expansion or shrinkage, and the resulted mapping explicitly reflects our intuition. Figure \ref{CustomizedFaces} displays the customized mesh parameterization for human facial surfaces, with circles highlighting the ROIs. 
We set the scalar of the density function to $k=3$. 
The t-OT results demonstrate the expected deformations. 
 The results clearly show that local regions within them have expanded according to the prescribed scalar, meanwhile, the vertices surrounding the circles have shrunk. However, there is no much change in other regions where the scalar is regarded as 1.0.

\subsection{Image Mesh Editing}
For image cases, we just convert the image to a planar triangular mesh with square or rectangle boundary for further t-OT and tt-OT processing. 
If the gray value is considered in the design of probability measure, then the resulted t-OT generates image content-aware deformations. There are the following applications for image mesh editing. 

\begin{figure*}
 \centering
\footnotesize
 	\begin{tabular}{@{\hspace{-0.5em}}m{0.01\linewidth}<{\centering}@{\hspace{0.6em}}m{0.23\linewidth}<{\centering}@{\hspace{0.5em}}m{0.23\linewidth}<{\centering}@{\hspace{0.5em}}m{0.23\linewidth}<{\centering}@{\hspace{0.5em}}m{0.23\linewidth}<{\centering}@{}}
 		\rotatebox{90}{Heart}&
 		\includegraphics[width=0.23\textwidth]{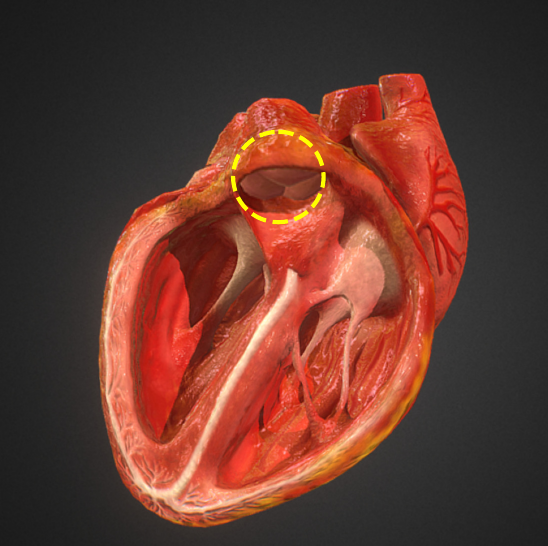} & 
 		\includegraphics[width=0.23\textwidth]{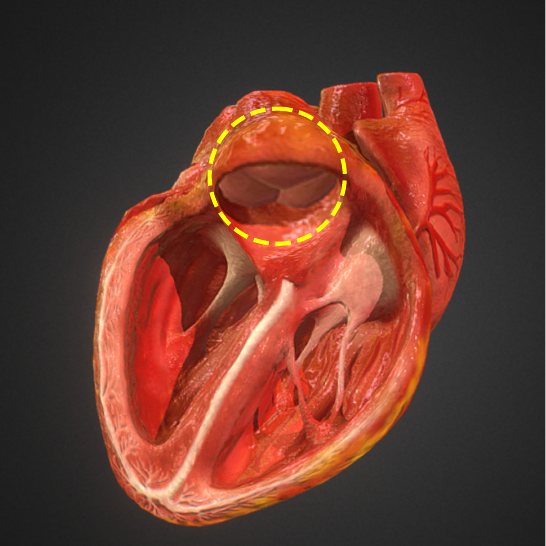}&
 		\includegraphics[width=0.23\textwidth]{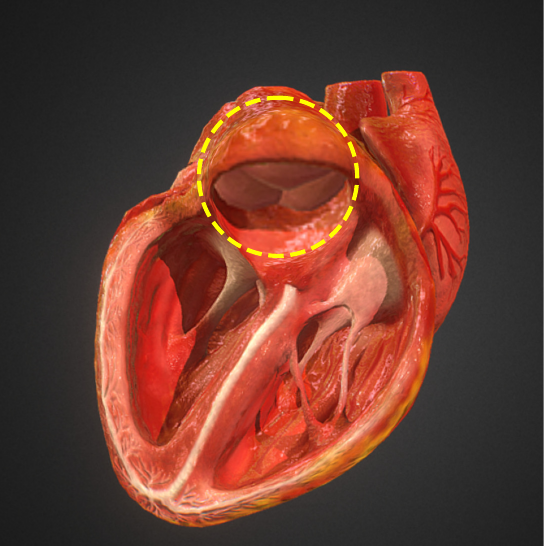}&
 		\rotatebox{90}{\includegraphics[width=0.23\textwidth]{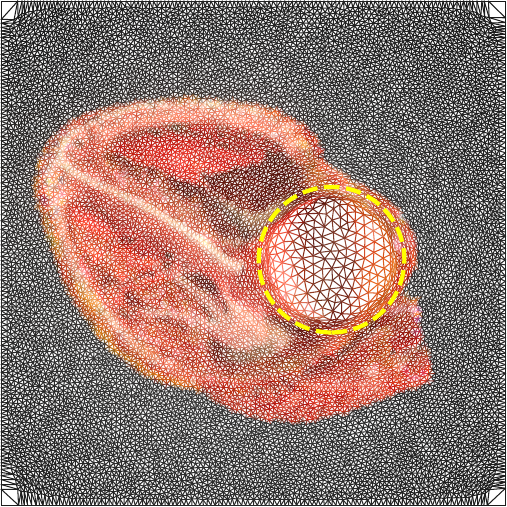}}\\
 		\rotatebox{90}{Cancer}&
 		\includegraphics[width=0.23\textwidth]{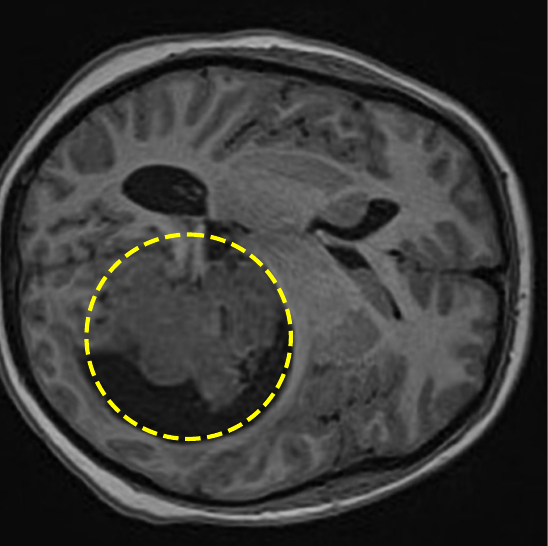} &
 		\includegraphics[width=0.23\textwidth]{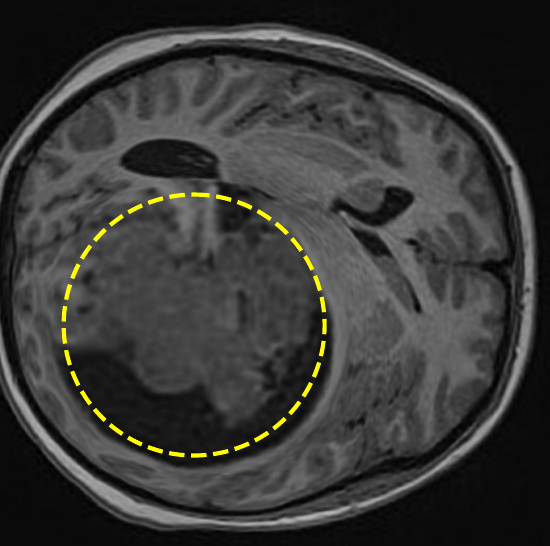}&	
 		\includegraphics[width=0.23\textwidth]{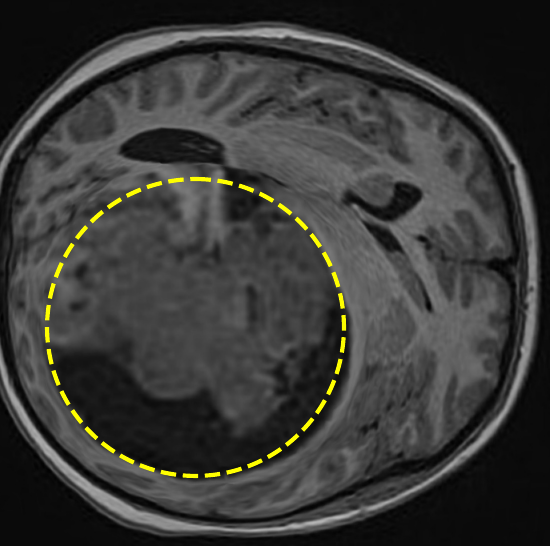} &
 		\includegraphics[width=0.23\textwidth]{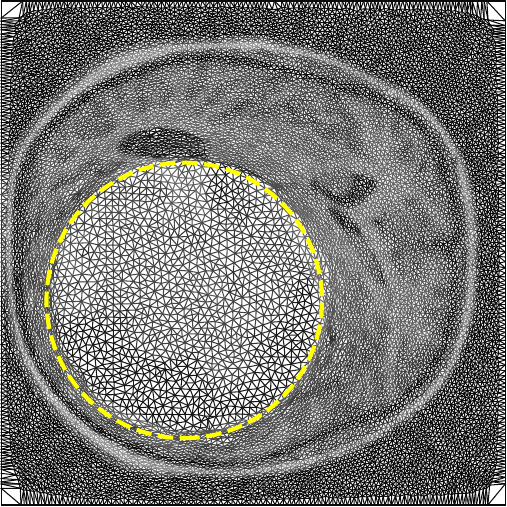} \\
 		\rotatebox{90}{Brain}&
 		\includegraphics[width=0.23\textwidth]{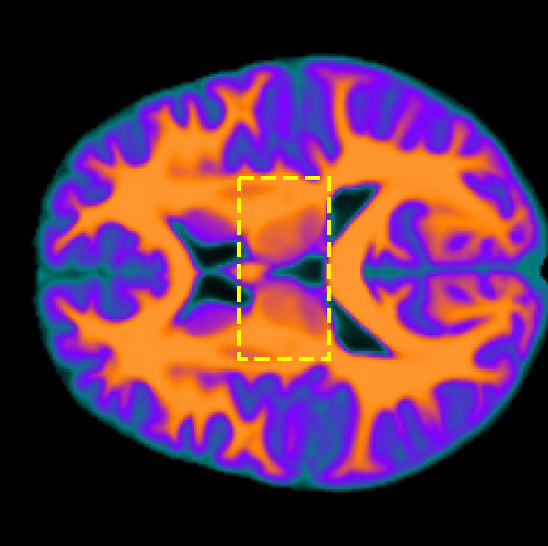} &
 		\includegraphics[width=0.23\textwidth]{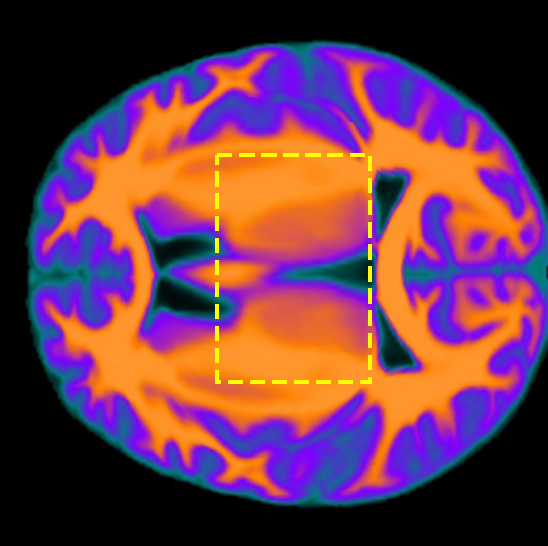}&	
 		\includegraphics[width=0.23\textwidth]{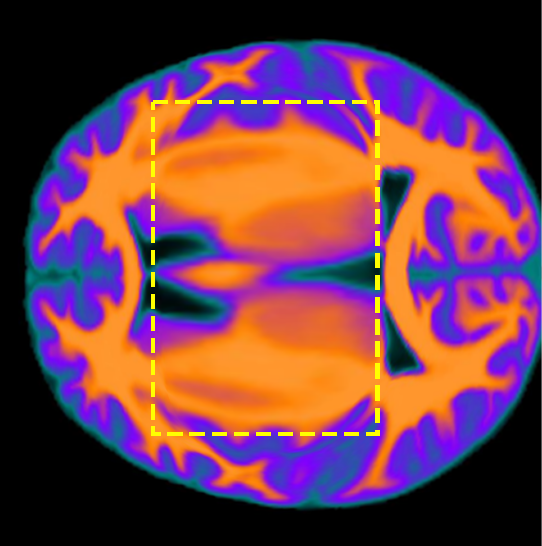}&
 		\reflectbox{\includegraphics[width=0.23\textwidth]{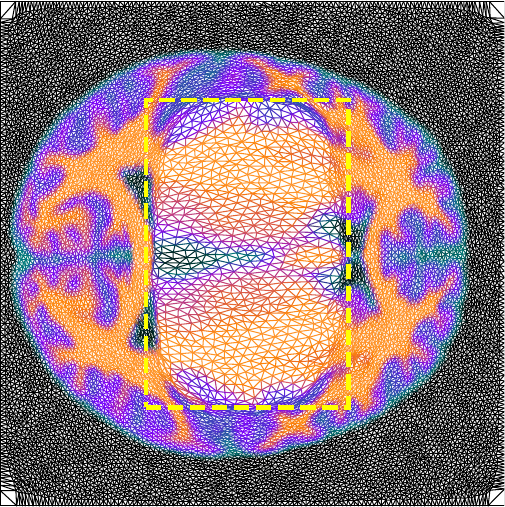}}  \\
 		&(a) Original& (b) $k=1.5$& (c) $k=2$ & (d) Mesh \\
 	\end{tabular}
 	\caption{Medical image magnifier with magnification factor $k = 1.5$ and $k=2$.
  } 
 	\label{Medical}
 \end{figure*}

\vspace{1mm}
\noindent\textbf{Medical Image Magnifier.}
We apply the scalar to the ROIs. Figure \ref{Medical} shows the original medical images in (a) and the ROIs labeled in a yellow circle (or rectangle). The density function are given as $g(p_{i})= k* a_{i}, p_{i} \in  \hat{O}$ with a scalar $k$. In this experiment, the mesh sizes for the three images are the same, i.e., $\left[ 13145, 25784 \right]$. Two constant scalars $k=1.5$ and $k=2$ are applied. As shown in Fig. \ref{Medical} (b-d), the t-OT results display expected deformations, with the ROIs magnified to the corresponding extents. As the magnification parameter \( k \) increases, the ROIs expand more, making biological features more visible to improve diagnosis. 
Therefore, the t-OT method can be effectively used as a medical image magnifier.

\begin{figure*}[t]
\centering
\footnotesize
\begin{tabular}{@{\hspace{-0.5em}}m{0.01\linewidth}<{\centering}@{\hspace{0.6em}}m{0.23\linewidth}<{\centering}@{\hspace{0.5em}}m{0.23\linewidth}<{\centering}@{\hspace{0.5em}}m{0.23\linewidth}<{\centering}@{\hspace{0.5em}}m{0.23\linewidth}<{\centering}@{}}
\rotatebox{90}{Cloud}&
\includegraphics[width=0.23\textwidth]{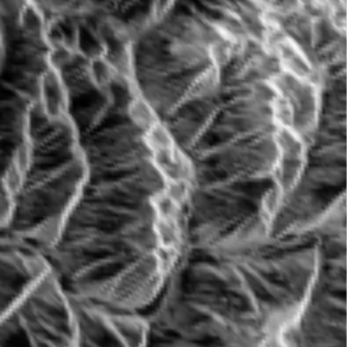}&
\includegraphics[width=0.23\textwidth]{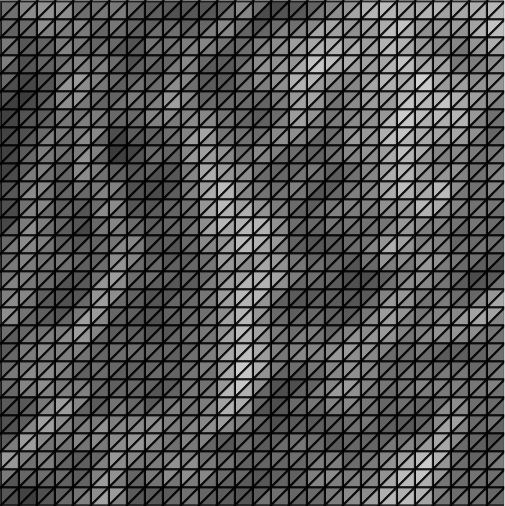}&
\includegraphics[width=0.23\textwidth]{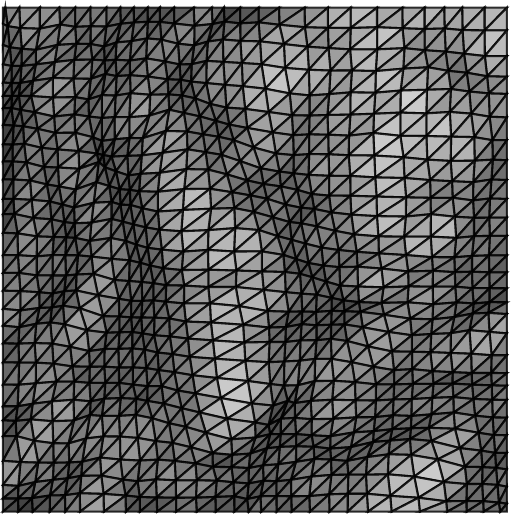}&
\includegraphics[width=0.23\textwidth]{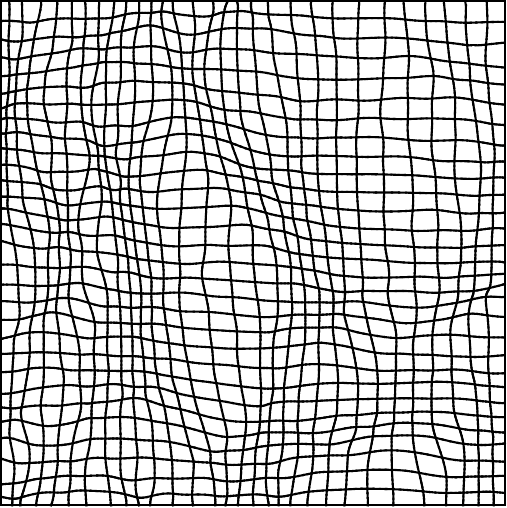}\\

\rotatebox{90}{Ring}&
\includegraphics[width=0.23\textwidth]{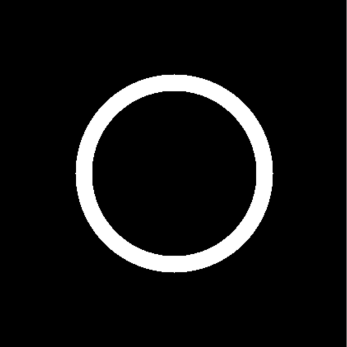}&
\includegraphics[width=0.23\textwidth]{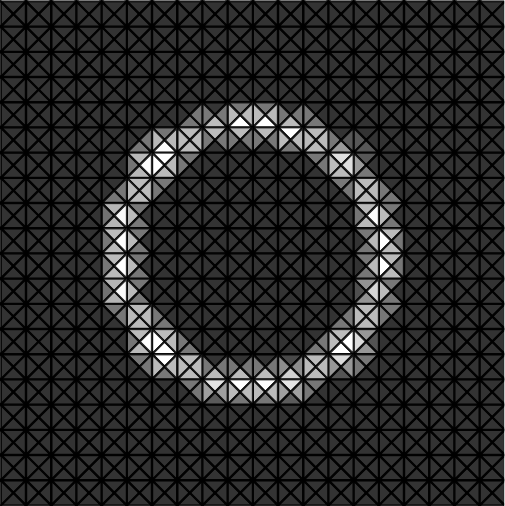}& 
\includegraphics[width=0.23\textwidth]{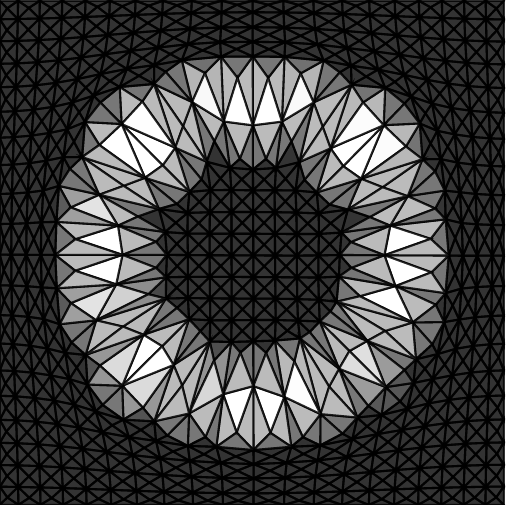}&
\includegraphics[width=0.23\textwidth]{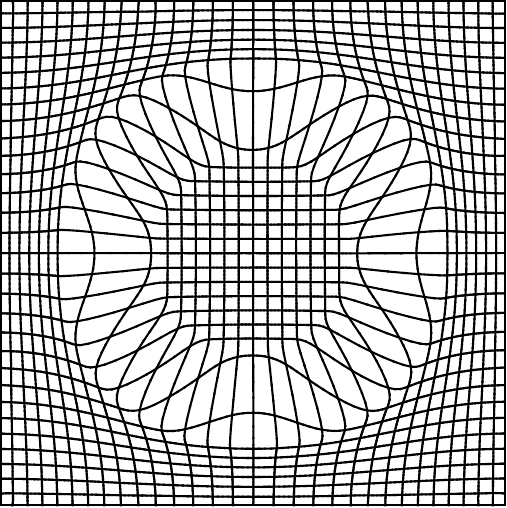}\\
\rotatebox{90}{Gaussian}&
\includegraphics[width=0.23\textwidth]{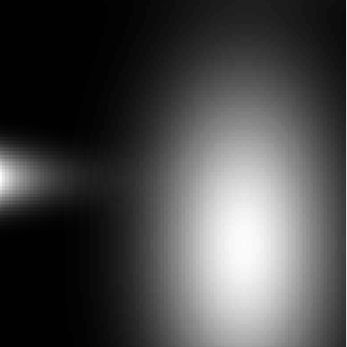}&
\includegraphics[width=0.23\textwidth]{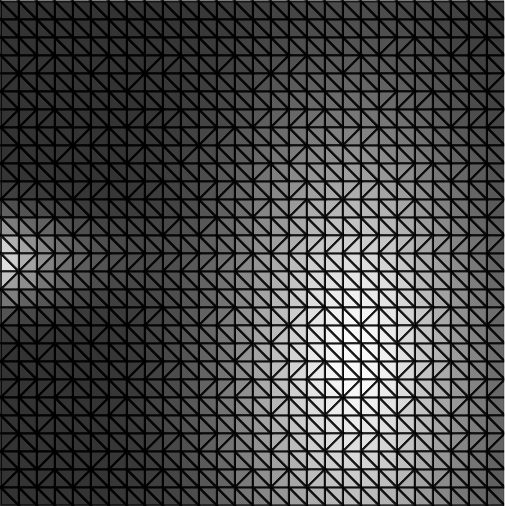}&
\includegraphics[width=0.23\textwidth]{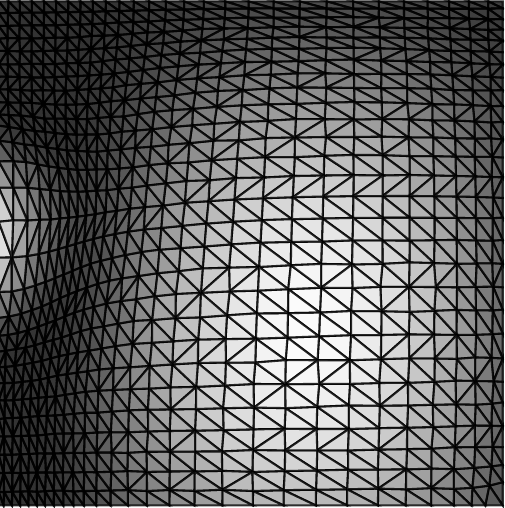}&
\includegraphics[width=0.23\textwidth]{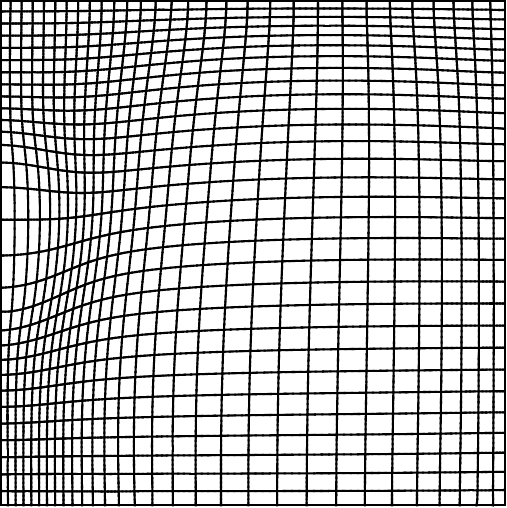}\\

&(a) Measure Image & (b) Initial Mesh & (c) t-OT Map& (d) Quad Mesh \\
\end{tabular}
\caption{Image-driven mesh generation with different $\delta=0.05,0.2,0.25$.} 
\label{fig:ImageDriven_OriginalData and Result}
\end{figure*}

\vspace{1mm}
\noindent\textbf{Image-Driven Mesh Generation.}
We use the intensity of grayscale image as the density function in t-OT for generating meshes, as shown in the last setting strategy of probability measure. Figure \ref{fig:ImageDriven_OriginalData and Result} shows the original grayscale images $I_i (i=1,2,3)$ and the corresponding triangular meshes. 
We set the scalar $k=1$ and the values of $\delta$ for $I_i(i=1,2,3)$ to $0.05$, $0.2$, and $0.25$, respectively. In this experiment, we have tried different discretizations of an image to a triangular mesh. 
Figure \ref{fig:ImageDriven_OriginalData and Result} gives the t-OT mapping results with the specified parameters as well as the quad mesh extracted based on the triangular results. The algorithm runs stable on different discretizations, giving similar results, as shown in Fig. \ref{fig:ImageDriven_OriginalData and Result}. The setting of probability measure determines that the bright areas expand and the dark areas shrink in the t-OT result. And, the brighter intensity generates more expansion. When the perturbing term $\delta$ increases, the darker area has relatively higher measure and thus shrinks less, as shown in (c-d).

\vspace{1mm}
\noindent\textbf{Physical Diffusion Simulation.}
In this experiment, we apply tt-OT into the original triangular mesh in Fig. \ref{fig:TemporalQCOT_TemporalResults} (a) by using the last setting strategy of probability measure with $k=1$ and $\delta = 0.1$.
Figure \ref{fig:TemporalQCOT_TemporalResults} illustrates the diffusion process. The first row shows the measure changes during the process, where the color-encoded height represents the density of probability measure $\psi _i = \frac{\nu _i}{\omega _i}$ and the volume represents the prescribed measure $\nu _i$ for each vertex. It is evident that as the time parameter $t$ increases, the height progressively decreases until all vertices have the same heights, i.e., visually the mesh has the same color everywhere and numerically the mesh has uniform density everywhere. This indicates that the total measure is preserved during the whole process and the volume at each vertex at $t=0$ eventually translates into the area of the power cell of each vertex on the result at $t=1$. 
We can specify the parameters in the image-driven density function to simulate the diffusion process and generate a sequence of meshes and images with progressive deformations. Then user could select the desired ones for further processing. 

In summary, t-OT and tt-OT provide comprehensive tools to control and edit image deformations and mesh structures. 
\begin{figure*}[h]
\centering
\footnotesize
 \setlength{\fboxsep}{0pt}
 \setlength{\tabcolsep}{2pt}
\begin{tabular}{ccccc}
\rotatebox{90}{\parbox{0.1602\textwidth}{\centering Measure}}&
\includegraphics[height=0.209\textwidth,width=0.23\textwidth]{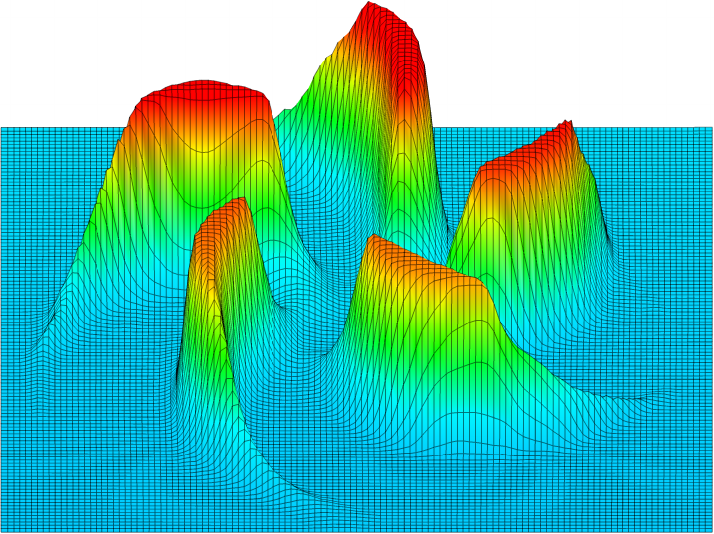}&
\includegraphics[height=0.1765\textwidth,width=0.23\textwidth]{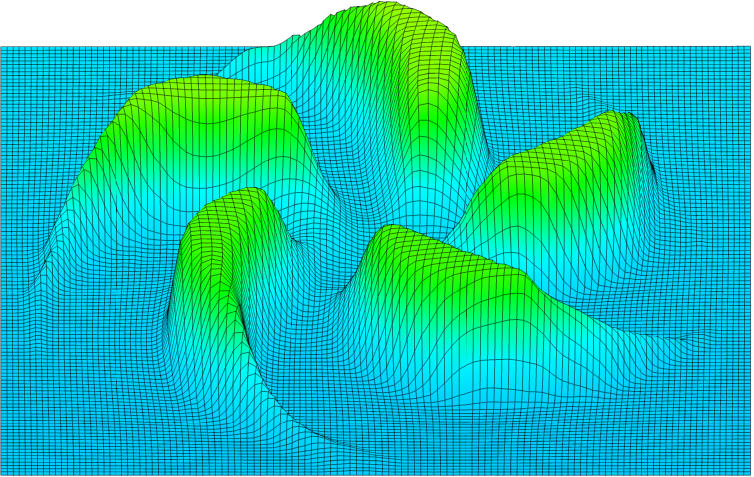}&
\includegraphics[height=0.1602\textwidth,width=0.23\textwidth]{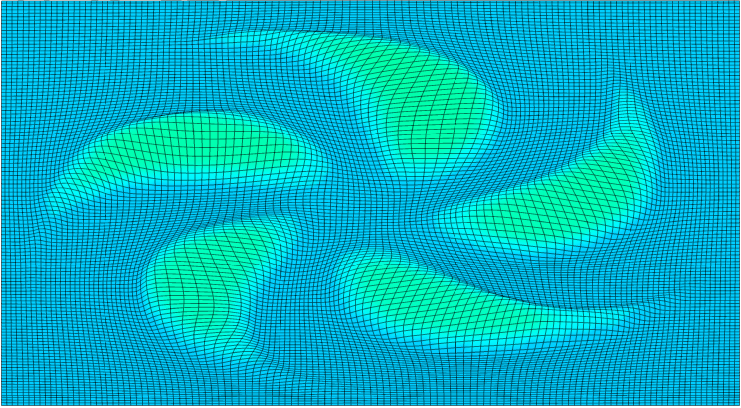}&
\includegraphics[height=0.1602\textwidth,width=0.23\textwidth]{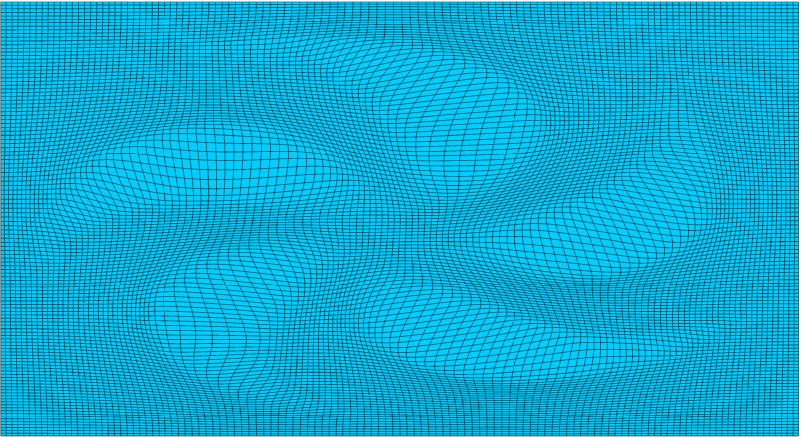}\\

\rotatebox{90}{\parbox{0.23\textwidth}{\centering Temporal Image}}&
\includegraphics[width=0.23\textwidth]{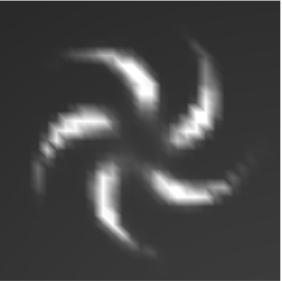}&
\includegraphics[width=0.23\textwidth]{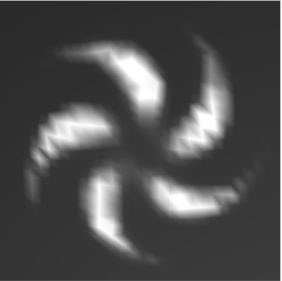}&
\includegraphics[width=0.23\textwidth]{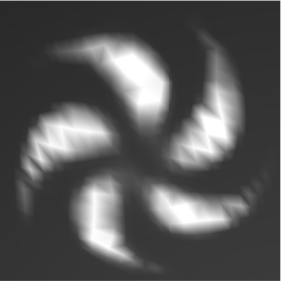}&
\includegraphics[width=0.23\textwidth]{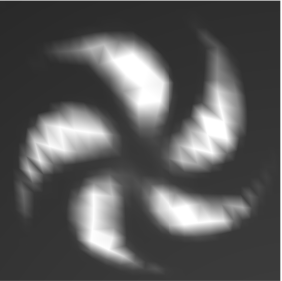}\\

\rotatebox{90}{\parbox{0.23\textwidth}{\centering Temporal Mesh}}&
\includegraphics[width=0.23\textwidth]{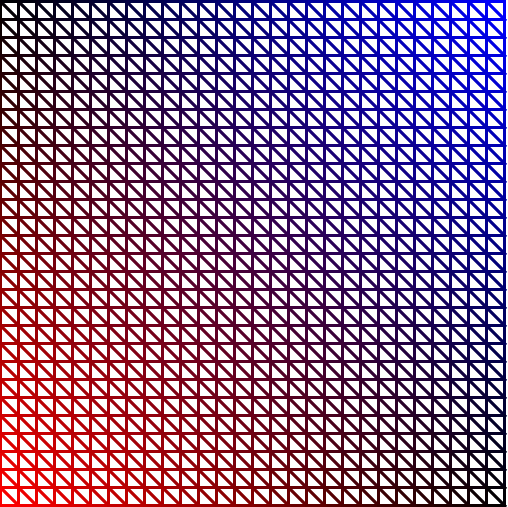}&
\includegraphics[width=0.23\textwidth]{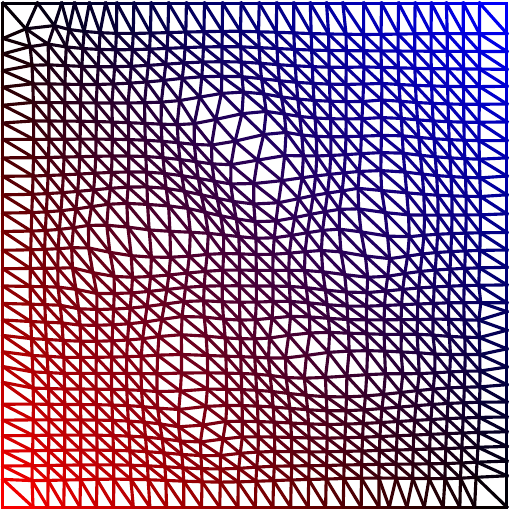}&
\includegraphics[width=0.23\textwidth]{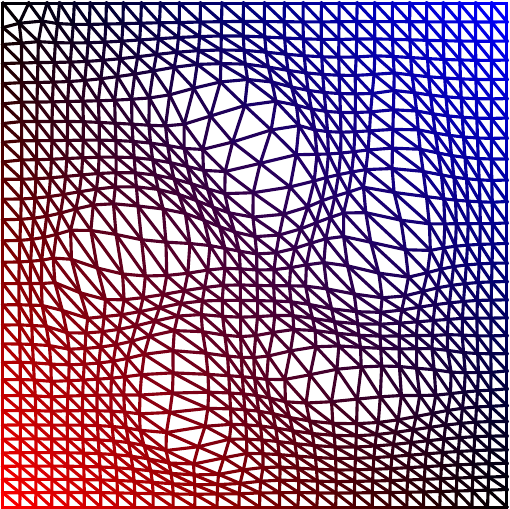}&
\includegraphics[width=0.23\textwidth]{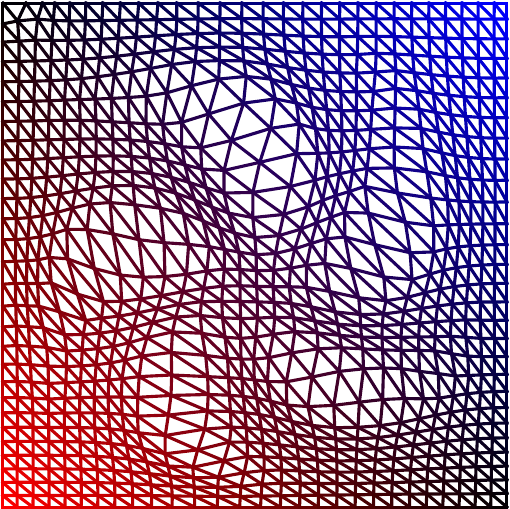}\\
\rotatebox{90}{\parbox{0.23\textwidth}{\centering Temporal Mesh}}&
\includegraphics[width=0.23\textwidth]{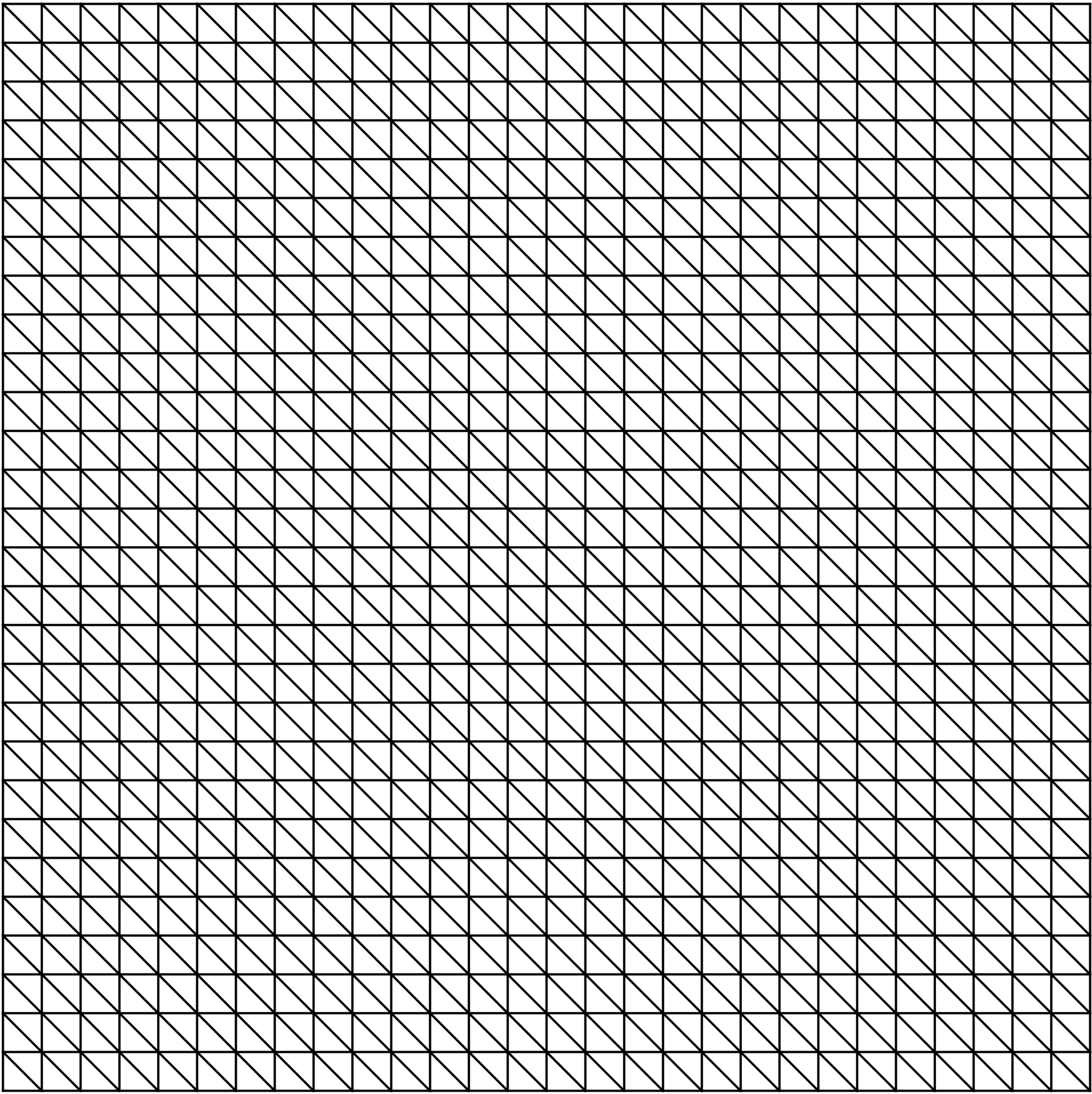}&
\includegraphics[width=0.23\textwidth]{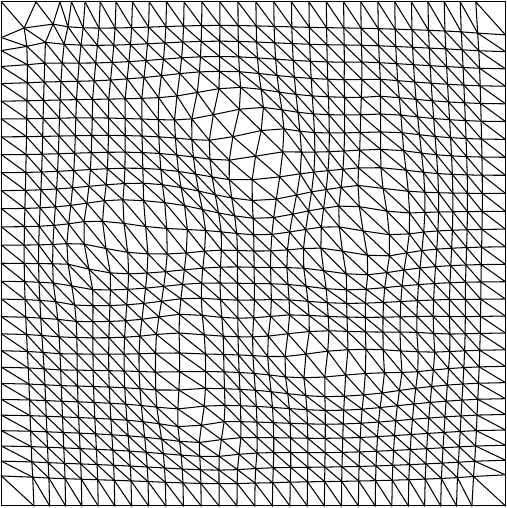}&
\includegraphics[width=0.23\textwidth]{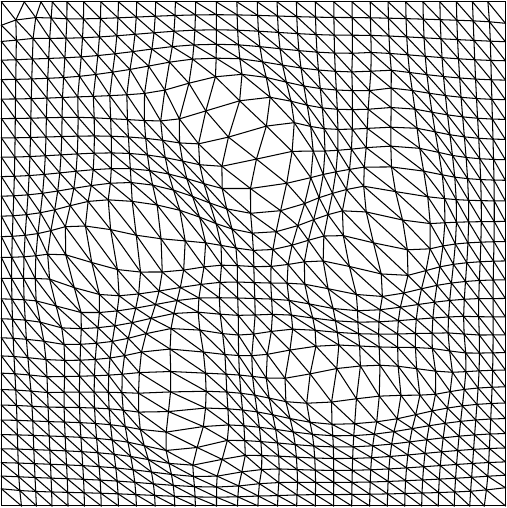}&
\includegraphics[width=0.23\textwidth]{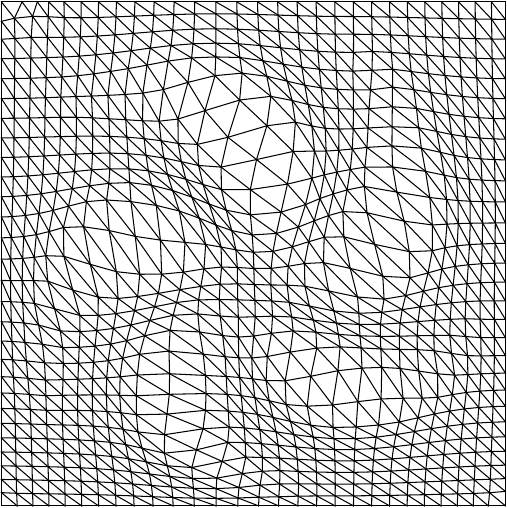}\\
&(a) $t = 0$&(b) $t = 0.1$&(c) $t = 0.5$&(d) $t = 1$\\

\end{tabular}
\caption{
The intermediate results generated by tt-OT with time parameter $t$. } 
\label{fig:TemporalQCOT_TemporalResults}
\end{figure*} 
\begin{figure*}[h]
\centering
\footnotesize
 \setlength{\fboxsep}{0pt}
 \setlength{\tabcolsep}{2pt}
\begin{tabular}{ccccc}
\rotatebox{90}{\parbox{0.1602\textwidth}{\centering Measure}}&
\includegraphics[height=0.209\textwidth,width=0.23\textwidth]{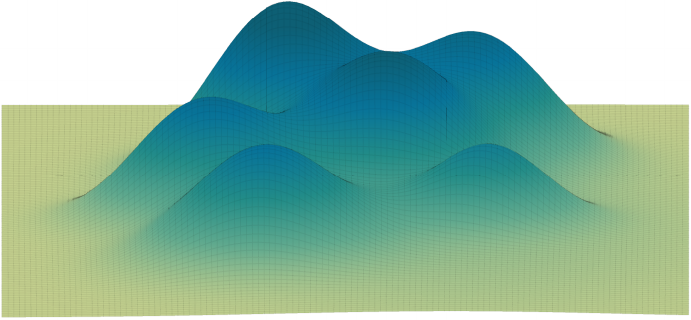}&
\includegraphics[height=0.156\textwidth,width=0.23\textwidth]{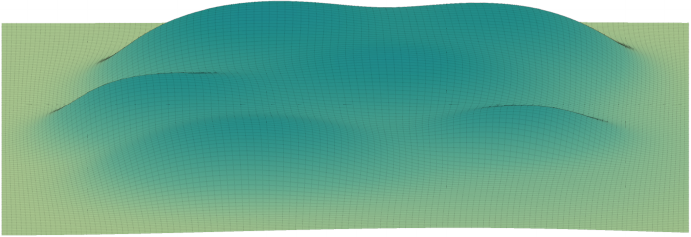}&
\includegraphics[height=0.150\textwidth,width=0.23\textwidth]{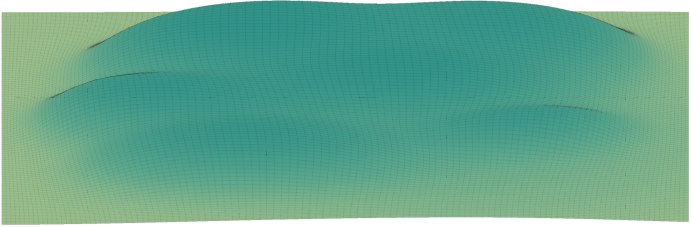}&
\includegraphics[height=0.144\textwidth,width=0.23\textwidth]{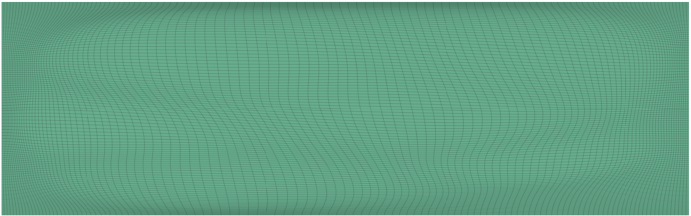}\\

\rotatebox{90}{\parbox{0.1602\textwidth}{\centering Measure}}&
\includegraphics[height=0.209\textwidth,width=0.23\textwidth]{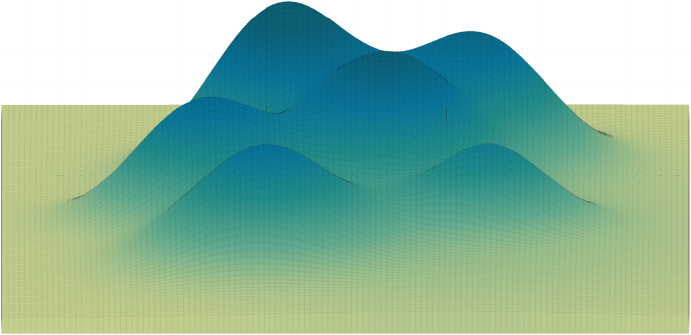}&
\includegraphics[height=0.158\textwidth,width=0.23\textwidth]{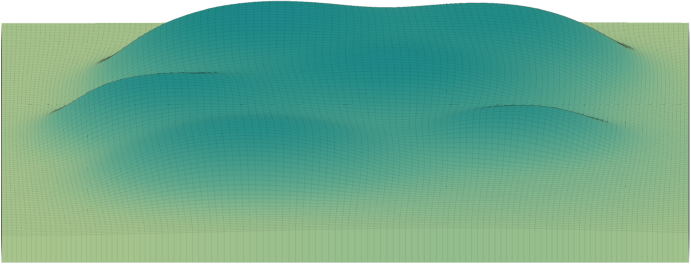}&
\includegraphics[height=0.152\textwidth,width=0.23\textwidth]{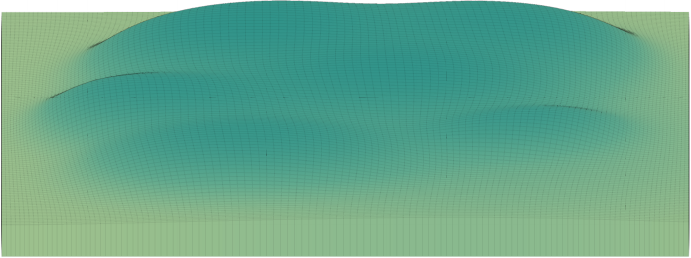}&
\includegraphics[height=0.146\textwidth,width=0.23\textwidth]{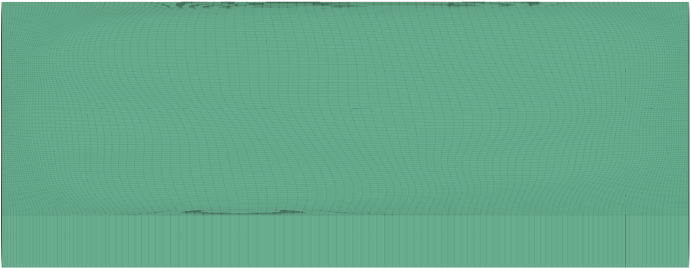}\\

\rotatebox{90}{\parbox{0.23\textwidth}{\centering Temporal Mesh}}&
\includegraphics[width=0.23\textwidth]{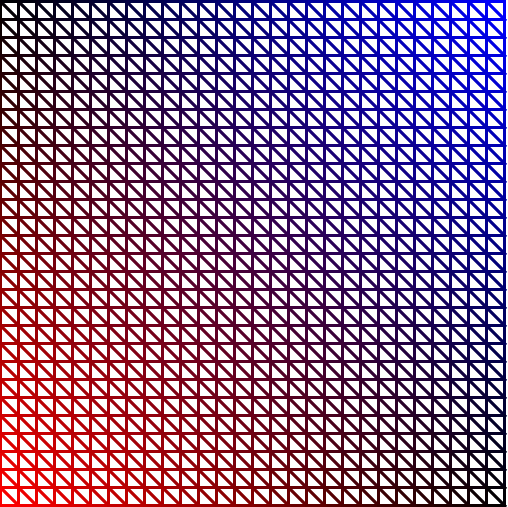}&
\includegraphics[width=0.23\textwidth]{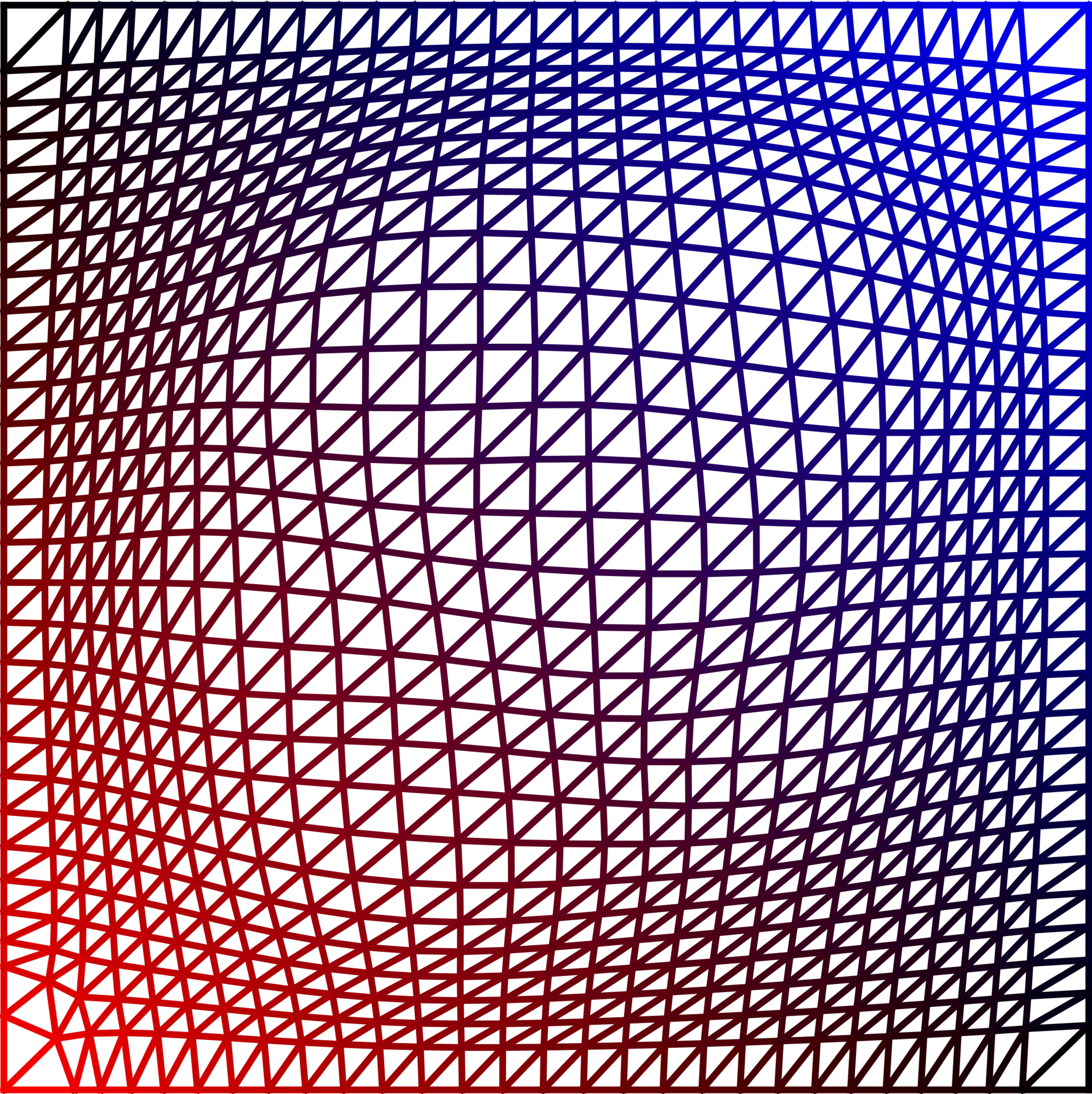}&
\includegraphics[width=0.23\textwidth]{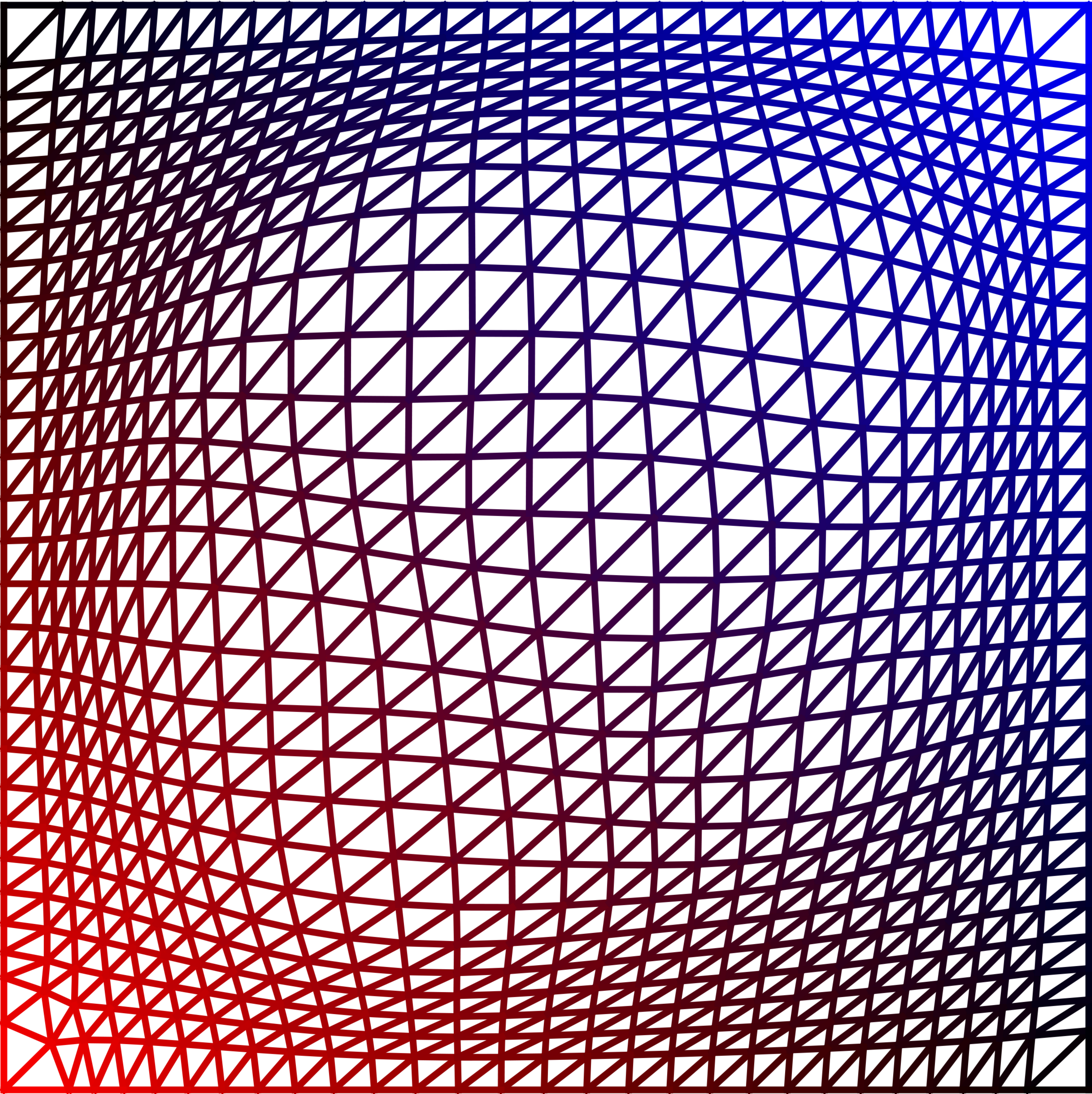}&
\includegraphics[width=0.23\textwidth]{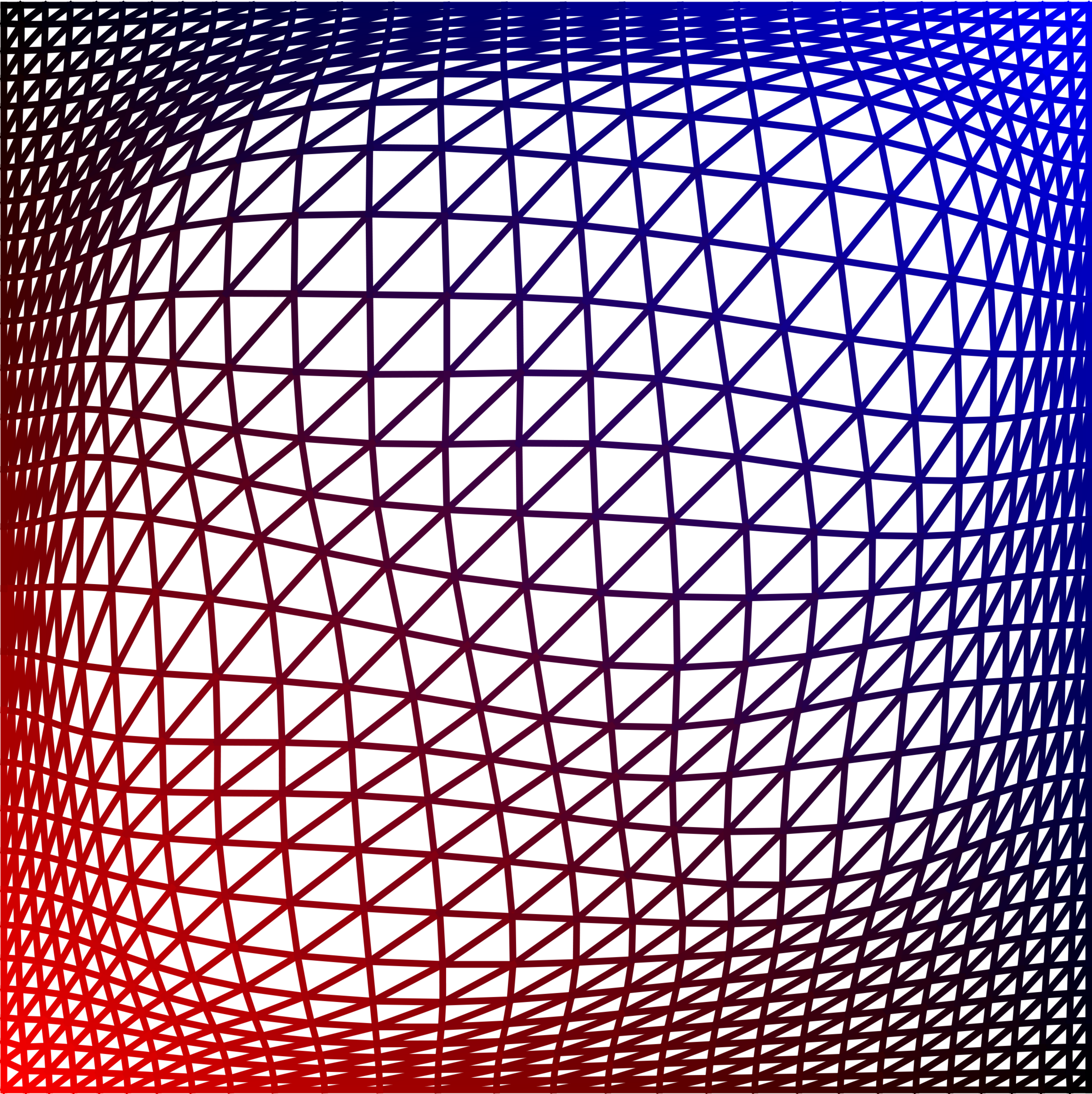}\\
\rotatebox{90}{\parbox{0.23\textwidth}{\centering Temporal Mesh}}&
\includegraphics[width=0.23\textwidth]{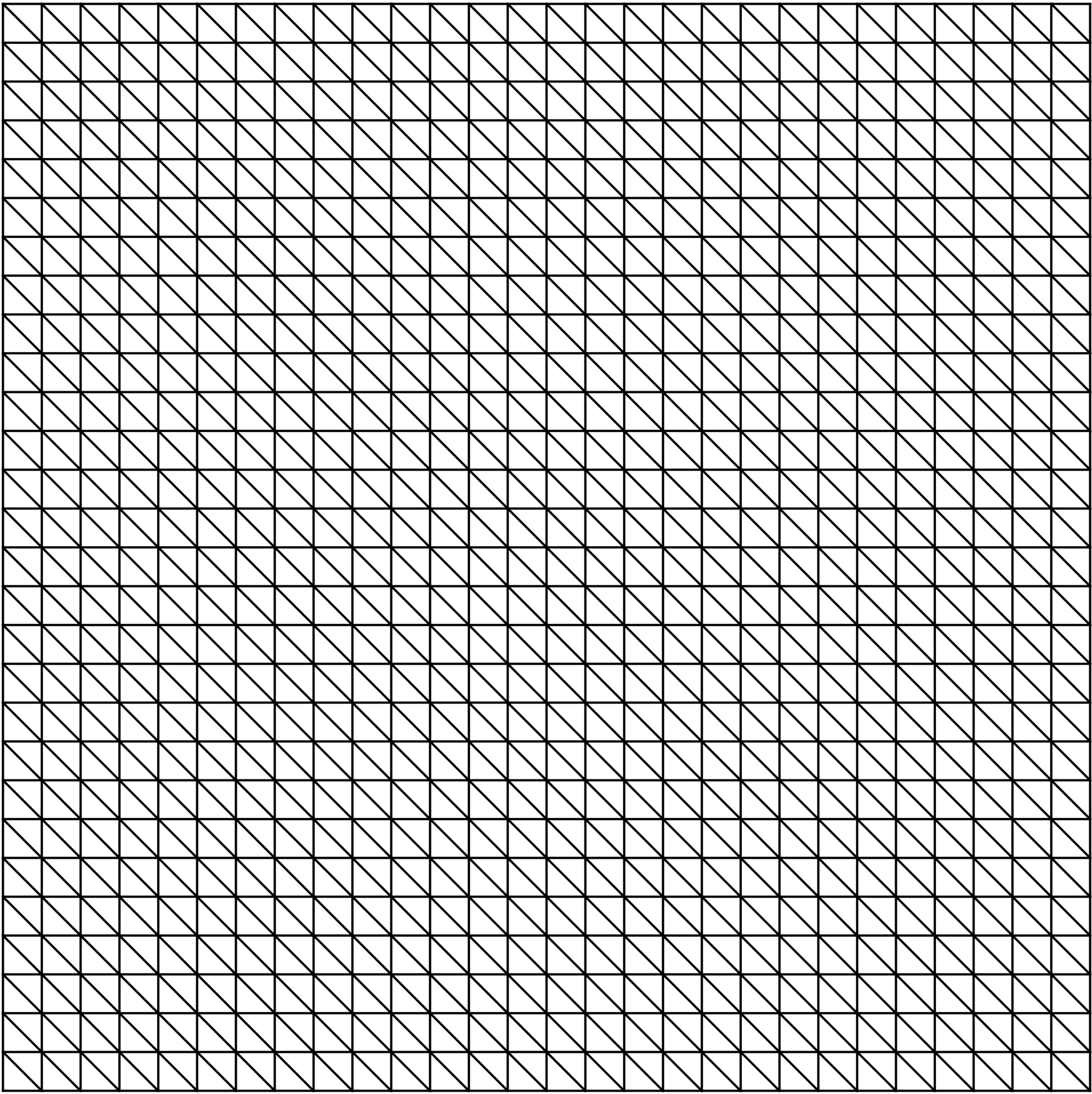}&
\includegraphics[width=0.23\textwidth]{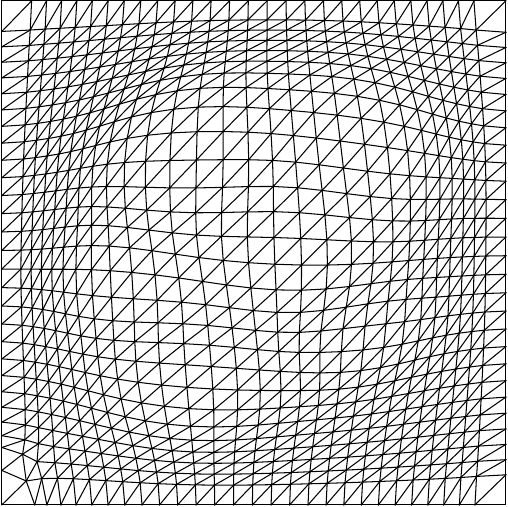}&
\includegraphics[width=0.23\textwidth]{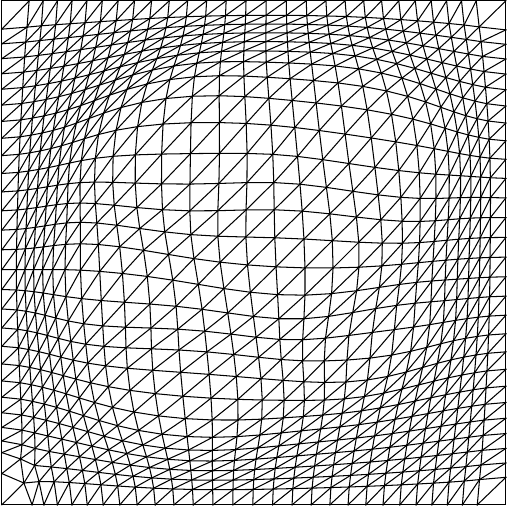}&
\includegraphics[width=0.23\textwidth]{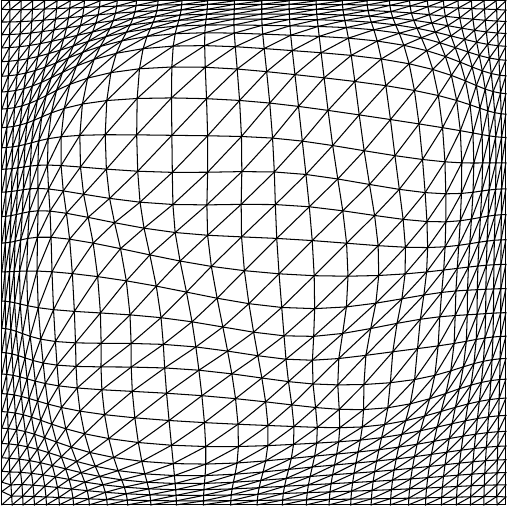}\\
&(a) $t = 0$&(b) $t = 0.1$&(c) $t = 0.5$&(d) $t = 1$\\

\end{tabular}
\caption{
The intermediate results generated by tt-OT with time parameter $t$. }
\label{fig:TemporalQCOT_TemporalResults}
\end{figure*} 

\section{Conclusion}
\label{sec:conclusion}
In this work, the concept of topology-preserving optimal transport (t-OT) is proposed.  We present the relaxed semi-discrete optimal transport by avoiding Delaunay triangulation operation and convexity check in the baseline semi-discrete optimal transport framework. Relaxed OT may cause unexpected distortions. Then we use quasiconformal geometric mapping to fill the ``bad'' regions without topology change. The proposed t-OT approximates the baseline sd-OT, but preserves topology, results higher quality OT result with the controllable distortions and runs faster. We further presented the temporal t-OT to study the dynamics of the transport. Thorough algorithm analysis and discussion are performed to enhance  understanding of the performance of the proposed methods. Experiments were conducted for surface mesh parameterization and image mesh editing, demonstrating high potential for general geometric processing tasks. In future, we will further explore the t-OT (tt-OT) framework for solving the registration problem between two domains where both original and prescribed measures are discrete, and create useful tools based on t-OT for shape analysis tasks in both engineering and medical imaging fields.

\textbf{Acknowledgments}

This work was supported in part by the National Key Research and Development Program of China (2021YFA1003002) and the National Natural Science Foundation of China (12090021 and 12090020).



\clearpage 





\printcredits



\begin{thebibliography}{49}
\expandafter\ifx\csname natexlab\endcsname\relax\def\natexlab#1{#1}\fi
\providecommand{\url}[1]{\texttt{#1}}
\providecommand{\href}[2]{#2}
\providecommand{\path}[1]{#1}
\providecommand{\DOIprefix}{doi:}
\providecommand{\ArXivprefix}{arXiv:}
\providecommand{\URLprefix}{URL: }
\providecommand{\Pubmedprefix}{pmid:}
\providecommand{\doi}[1]{\href{http://dx.doi.org/#1}{\path{#1}}}
\providecommand{\Pubmed}[1]{\href{pmid:#1}{\path{#1}}}
\providecommand{\bibinfo}[2]{#2}
\ifx\xfnm\relax \def\xfnm[#1]{\unskip,\space#1}\fi
\bibitem[{Ahlfors(2006)}]{ahlfors2006lectures}
\bibinfo{author}{Ahlfors, L.V.}, \bibinfo{year}{2006}.
\newblock \bibinfo{title}{Lectures on quasiconformal mappings}. volume~\bibinfo{volume}{38}.
\newblock \bibinfo{publisher}{American Mathematical Soc.}, \bibinfo{address}{Rhode Island}.
\bibitem[{Arjovsky et~al.(2017)Arjovsky, Chintala and Bottou}]{arjovsky2017wasserstein}
\bibinfo{author}{Arjovsky, M.}, \bibinfo{author}{Chintala, S.}, \bibinfo{author}{Bottou, L.}, \bibinfo{year}{2017}.
\newblock \bibinfo{title}{Wasserstein generative adversarial networks}, in: \bibinfo{booktitle}{International conference on machine learning}, \bibinfo{organization}{PMLR}. pp. \bibinfo{pages}{214--223}.
\bibitem[{Arroyo~Ohori et~al.(2015)Arroyo~Ohori, Ledoux and Stoter}]{arroyo2015evaluation}
\bibinfo{author}{Arroyo~Ohori, K.}, \bibinfo{author}{Ledoux, H.}, \bibinfo{author}{Stoter, J.}, \bibinfo{year}{2015}.
\newblock \bibinfo{title}{An evaluation and classification of n d topological data structures for the representation of objects in a higher-dimensional gis}.
\newblock \bibinfo{journal}{International Journal of Geographical Information Science} \bibinfo{volume}{29}, \bibinfo{pages}{825--849}.
\bibitem[{Benamou and Brenier(2000)}]{benamou2000computational}
\bibinfo{author}{Benamou, J.D.}, \bibinfo{author}{Brenier, Y.}, \bibinfo{year}{2000}.
\newblock \bibinfo{title}{A computational fluid mechanics solution to the monge-kantorovich mass transfer problem}.
\newblock \bibinfo{journal}{Numerische Mathematik} \bibinfo{volume}{84}, \bibinfo{pages}{375--393}.
\bibitem[{Bonneel et~al.(2011)Bonneel, Van De~Panne, Paris and Heidrich}]{bonneel2011displacement}
\bibinfo{author}{Bonneel, N.}, \bibinfo{author}{Van De~Panne, M.}, \bibinfo{author}{Paris, S.}, \bibinfo{author}{Heidrich, W.}, \bibinfo{year}{2011}.
\newblock \bibinfo{title}{Displacement interpolation using lagrangian mass transport}, in: \bibinfo{booktitle}{Proceedings of the 2011 SIGGRAPH Asia conference}, pp. \bibinfo{pages}{1--12}.
\bibitem[{Brenier(1991)}]{brenier1991polar}
\bibinfo{author}{Brenier, Y.}, \bibinfo{year}{1991}.
\newblock \bibinfo{title}{Polar factorization and monotone rearrangement of vector-valued functions}.
\newblock \bibinfo{journal}{Communications on pure and applied mathematics} \bibinfo{volume}{44}, \bibinfo{pages}{375--417}.
\bibitem[{Chen(2021)}]{chen2021spatiotemporal}
\bibinfo{author}{Chen, C.}, \bibinfo{year}{2021}.
\newblock \bibinfo{title}{Spatiotemporal imaging with diffeomorphic optimal transportation}.
\newblock \bibinfo{journal}{Inverse Problems} \bibinfo{volume}{37}, \bibinfo{pages}{115004}.
\bibitem[{Choi et~al.(2022)Choi, Giri and Kumar}]{choi2022adaptive}
\bibinfo{author}{Choi, G.P.}, \bibinfo{author}{Giri, A.}, \bibinfo{author}{Kumar, L.}, \bibinfo{year}{2022}.
\newblock \bibinfo{title}{Adaptive area-preserving parameterization of open and closed anatomical surfaces}.
\newblock \bibinfo{journal}{Computers in Biology and Medicine} \bibinfo{volume}{148}, \bibinfo{pages}{105715}.
\bibitem[{Clough et~al.(2020)Clough, Byrne, Oksuz, Zimmer, Schnabel and King}]{clough2020topological}
\bibinfo{author}{Clough, J.R.}, \bibinfo{author}{Byrne, N.}, \bibinfo{author}{Oksuz, I.}, \bibinfo{author}{Zimmer, V.A.}, \bibinfo{author}{Schnabel, J.A.}, \bibinfo{author}{King, A.P.}, \bibinfo{year}{2020}.
\newblock \bibinfo{title}{A topological loss function for deep-learning based image segmentation using persistent homology}.
\newblock \bibinfo{journal}{IEEE transactions on pattern analysis and machine intelligence} \bibinfo{volume}{44}, \bibinfo{pages}{8766--8778}.
\bibitem[{Estellers et~al.(2012)Estellers, Zosso, Lai, Osher, Thiran and Bresson}]{estellers2012efficient}
\bibinfo{author}{Estellers, V.}, \bibinfo{author}{Zosso, D.}, \bibinfo{author}{Lai, R.}, \bibinfo{author}{Osher, S.}, \bibinfo{author}{Thiran, J.P.}, \bibinfo{author}{Bresson, X.}, \bibinfo{year}{2012}.
\newblock \bibinfo{title}{Efficient algorithm for level set method preserving distance function}.
\newblock \bibinfo{journal}{IEEE Transactions on Image Processing} \bibinfo{volume}{21}, \bibinfo{pages}{4722--4734}.
\bibitem[{Flamary et~al.(2016)Flamary, Courty, Tuia and Rakotomamonjy}]{flamary2016optimal}
\bibinfo{author}{Flamary, R.}, \bibinfo{author}{Courty, N.}, \bibinfo{author}{Tuia, D.}, \bibinfo{author}{Rakotomamonjy, A.}, \bibinfo{year}{2016}.
\newblock \bibinfo{title}{Optimal transport for domain adaptation}.
\newblock \bibinfo{journal}{IEEE Trans. Pattern Anal. Mach. Intell} \bibinfo{volume}{1}, \bibinfo{pages}{2}.
\bibitem[{Floater and Hormann(2005)}]{floater2005surface}
\bibinfo{author}{Floater, M.S.}, \bibinfo{author}{Hormann, K.}, \bibinfo{year}{2005}.
\newblock \bibinfo{title}{Surface parameterization: a tutorial and survey}.
\newblock \bibinfo{journal}{Advances in multiresolution for geometric modelling} , \bibinfo{pages}{157--186}.
\bibitem[{Gu et~al.(2013)Gu, Luo, Sun and Yau}]{gu2013variational}
\bibinfo{author}{Gu, X.}, \bibinfo{author}{Luo, F.}, \bibinfo{author}{Sun, J.}, \bibinfo{author}{Yau, S.T.}, \bibinfo{year}{2013}.
\newblock \bibinfo{title}{Variational principles for minkowski type problems, discrete optimal transport, and discrete monge-ampere equations}.
\newblock \bibinfo{journal}{arXiv preprint arXiv:1302.5472} .
\bibitem[{Gu and Yau(2008)}]{gu2008computational}
\bibinfo{author}{Gu, X.D.}, \bibinfo{author}{Yau, S.T.}, \bibinfo{year}{2008}.
\newblock \bibinfo{title}{Computational conformal geometry}.
\newblock \bibinfo{journal}{(No Title)} .
\bibitem[{Haker et~al.(2004)Haker, Zhu, Tannenbaum and Angenent}]{haker2004optimal}
\bibinfo{author}{Haker, S.}, \bibinfo{author}{Zhu, L.}, \bibinfo{author}{Tannenbaum, A.}, \bibinfo{author}{Angenent, S.}, \bibinfo{year}{2004}.
\newblock \bibinfo{title}{Optimal mass transport for registration and warping}.
\newblock \bibinfo{journal}{International Journal of computer vision} \bibinfo{volume}{60}, \bibinfo{pages}{225--240}.
\bibitem[{Helfer et~al.(2013)Helfer, Springborn and Suris}]{helfer2013secondary}
\bibinfo{author}{Helfer, J.}, \bibinfo{author}{Springborn, B.}, \bibinfo{author}{Suris, Y.}, \bibinfo{year}{2013}.
\newblock \bibinfo{title}{The secondary fan of a punctured riemann surface} .
\bibitem[{Ho et~al.(2017)Ho, Nguyen, Yurochkin, Bui, Huynh and Phung}]{ho2017multilevel}
\bibinfo{author}{Ho, N.}, \bibinfo{author}{Nguyen, X.}, \bibinfo{author}{Yurochkin, M.}, \bibinfo{author}{Bui, H.H.}, \bibinfo{author}{Huynh, V.}, \bibinfo{author}{Phung, D.}, \bibinfo{year}{2017}.
\newblock \bibinfo{title}{Multilevel clustering via wasserstein means}, in: \bibinfo{booktitle}{International conference on machine learning}, \bibinfo{organization}{PMLR}. pp. \bibinfo{pages}{1501--1509}.
\bibitem[{Huang et~al.(2016)Huang, Guo, Kusner, Sun, Sha and Weinberger}]{huang2016supervised}
\bibinfo{author}{Huang, G.}, \bibinfo{author}{Guo, C.}, \bibinfo{author}{Kusner, M.J.}, \bibinfo{author}{Sun, Y.}, \bibinfo{author}{Sha, F.}, \bibinfo{author}{Weinberger, K.Q.}, \bibinfo{year}{2016}.
\newblock \bibinfo{title}{Supervised word mover's distance}.
\newblock \bibinfo{journal}{Advances in neural information processing systems} \bibinfo{volume}{29}.
\bibitem[{Izquierdo and Civera(2024)}]{izquierdo2024optimal}
\bibinfo{author}{Izquierdo, S.}, \bibinfo{author}{Civera, J.}, \bibinfo{year}{2024}.
\newblock \bibinfo{title}{Optimal transport aggregation for visual place recognition}, in: \bibinfo{booktitle}{Proceedings of the IEEE/CVF Conference on Computer Vision and Pattern Recognition}, pp. \bibinfo{pages}{17658--17668}.
\bibitem[{Kantorovich(2006)}]{kantorovich2006problem}
\bibinfo{author}{Kantorovich, L.V.}, \bibinfo{year}{2006}.
\newblock \bibinfo{title}{On a problem of monge}.
\newblock \bibinfo{journal}{Journal of Mathematical Sciences} \bibinfo{volume}{133}, \bibinfo{pages}{1383--1383}.
\bibitem[{Kusner et~al.(2015)Kusner, Sun, Kolkin and Weinberger}]{kusner2015word}
\bibinfo{author}{Kusner, M.}, \bibinfo{author}{Sun, Y.}, \bibinfo{author}{Kolkin, N.}, \bibinfo{author}{Weinberger, K.}, \bibinfo{year}{2015}.
\newblock \bibinfo{title}{From word embeddings to document distances}, in: \bibinfo{booktitle}{International conference on machine learning}, \bibinfo{organization}{PMLR}. pp. \bibinfo{pages}{957--966}.
\bibitem[{Lam et~al.(2015)Lam, Gu and Lui}]{lam2015landmark}
\bibinfo{author}{Lam, K.C.}, \bibinfo{author}{Gu, X.}, \bibinfo{author}{Lui, L.M.}, \bibinfo{year}{2015}.
\newblock \bibinfo{title}{Landmark constrained genus-one surface teichm{\"u}ller map applied to surface registration in medical imaging}.
\newblock \bibinfo{journal}{Medical image analysis} \bibinfo{volume}{25}, \bibinfo{pages}{45--55}.
\bibitem[{Lei et~al.(2019)Lei, Chen, Luo, Si and Gu}]{lei2019secondary}
\bibinfo{author}{Lei, N.}, \bibinfo{author}{Chen, W.}, \bibinfo{author}{Luo, Z.}, \bibinfo{author}{Si, H.}, \bibinfo{author}{Gu, X.}, \bibinfo{year}{2019}.
\newblock \bibinfo{title}{Secondary polytope and secondary power diagram}.
\newblock \bibinfo{journal}{Computational Mathematics and Mathematical Physics} \bibinfo{volume}{59}, \bibinfo{pages}{1965--1981}.
\bibitem[{Lei and Gu(2021)}]{lei2021optimal}
\bibinfo{author}{Lei, N.}, \bibinfo{author}{Gu, X.}, \bibinfo{year}{2021}.
\newblock \bibinfo{title}{Optimal transportation theory and computation}.
\bibitem[{L{\'e}vy(2006)}]{levy2006parameterization}
\bibinfo{author}{L{\'e}vy, B.}, \bibinfo{year}{2006}.
\newblock \bibinfo{title}{Parameterization of mesh-models: theory, implementation and applications}, in: \bibinfo{booktitle}{Proceedings of the 2006 ACM symposium on Solid and physical modeling}, pp. \bibinfo{pages}{171--171}.
\bibitem[{Liu and Ng(2024)}]{liu2024one}
\bibinfo{author}{Liu, C.}, \bibinfo{author}{Ng, M.K.}, \bibinfo{year}{2024}.
\newblock \bibinfo{title}{A one-step image retargeing algorithm based on conformal energy}.
\newblock \bibinfo{journal}{arXiv e-prints} , \bibinfo{pages}{arXiv--2402}.
\bibitem[{Lui et~al.(2013)Lui, Lam, Wong and Gu}]{lui2013texture}
\bibinfo{author}{Lui, L.M.}, \bibinfo{author}{Lam, K.C.}, \bibinfo{author}{Wong, T.W.}, \bibinfo{author}{Gu, X.}, \bibinfo{year}{2013}.
\newblock \bibinfo{title}{Texture map and video compression using beltrami representation}.
\newblock \bibinfo{journal}{SIAM Journal on Imaging Sciences} \bibinfo{volume}{6}, \bibinfo{pages}{1880--1902}.
\bibitem[{Lui et~al.(2012)Lui, Wong, Zeng, Gu, Thompson, Chan and Yau}]{lui2012optimization}
\bibinfo{author}{Lui, L.M.}, \bibinfo{author}{Wong, T.W.}, \bibinfo{author}{Zeng, W.}, \bibinfo{author}{Gu, X.}, \bibinfo{author}{Thompson, P.M.}, \bibinfo{author}{Chan, T.F.}, \bibinfo{author}{Yau, S.T.}, \bibinfo{year}{2012}.
\newblock \bibinfo{title}{Optimization of surface registrations using beltrami holomorphic flow}.
\newblock \bibinfo{journal}{Journal of scientific computing} \bibinfo{volume}{50}, \bibinfo{pages}{557--585}.
\bibitem[{Lyu et~al.(2023)Lyu, Choi and Lui}]{lyu2023bijective}
\bibinfo{author}{Lyu, Z.}, \bibinfo{author}{Choi, G.}, \bibinfo{author}{Lui, L.M.}, \bibinfo{year}{2023}.
\newblock \bibinfo{title}{Bijective density-equalizing quasiconformal map for multiply-connected open surfaces}.
\newblock \bibinfo{journal}{arXiv preprint arXiv:2308.05579} .
\bibitem[{Ma et~al.(2017)Ma, Lei, Chen, Su and Gu}]{ma2017robust}
\bibinfo{author}{Ma, M.}, \bibinfo{author}{Lei, N.}, \bibinfo{author}{Chen, W.}, \bibinfo{author}{Su, K.}, \bibinfo{author}{Gu, X.}, \bibinfo{year}{2017}.
\newblock \bibinfo{title}{Robust surface registration using optimal mass transport and teichm{\"u}ller mapping}.
\newblock \bibinfo{journal}{Graphical models} \bibinfo{volume}{90}, \bibinfo{pages}{13--23}.
\bibitem[{Monge(1781)}]{monge1781memoire}
\bibinfo{author}{Monge, G.}, \bibinfo{year}{1781}.
\newblock \bibinfo{title}{M{\'e}moire sur la th{\'e}orie des d{\'e}blais et des remblais}.
\newblock \bibinfo{journal}{Mem. Math. Phys. Acad. Royale Sci.} , \bibinfo{pages}{666--704}.
\bibitem[{Mullen et~al.(2008)Mullen, Tong, Alliez and Desbrun}]{mullen2008spectral}
\bibinfo{author}{Mullen, P.}, \bibinfo{author}{Tong, Y.}, \bibinfo{author}{Alliez, P.}, \bibinfo{author}{Desbrun, M.}, \bibinfo{year}{2008}.
\newblock \bibinfo{title}{Spectral conformal parameterization}, in: \bibinfo{booktitle}{Computer Graphics Forum}, \bibinfo{organization}{Wiley Online Library}. pp. \bibinfo{pages}{1487--1494}.
\bibitem[{Pan et~al.(2019)Pan, Han, Chen, Tang and Jia}]{pan2019deep}
\bibinfo{author}{Pan, J.}, \bibinfo{author}{Han, X.}, \bibinfo{author}{Chen, W.}, \bibinfo{author}{Tang, J.}, \bibinfo{author}{Jia, K.}, \bibinfo{year}{2019}.
\newblock \bibinfo{title}{Deep mesh reconstruction from single rgb images via topology modification networks}, in: \bibinfo{booktitle}{Proceedings of the IEEE/CVF International Conference on Computer Vision}, pp. \bibinfo{pages}{9964--9973}.
\bibitem[{Peyr{\'e} et~al.(2019)Peyr{\'e}, Cuturi et~al.}]{peyre2019computational}
\bibinfo{author}{Peyr{\'e}, G.}, \bibinfo{author}{Cuturi, M.}, et~al., \bibinfo{year}{2019}.
\newblock \bibinfo{title}{Computational optimal transport: With applications to data science}.
\newblock \bibinfo{journal}{Foundations and Trends{\textregistered} in Machine Learning} \bibinfo{volume}{11}, \bibinfo{pages}{355--607}.
\bibitem[{Solomon et~al.(2015)Solomon, De~Goes, Peyr{\'e}, Cuturi, Butscher, Nguyen, Du and Guibas}]{solomon2015convolutional}
\bibinfo{author}{Solomon, J.}, \bibinfo{author}{De~Goes, F.}, \bibinfo{author}{Peyr{\'e}, G.}, \bibinfo{author}{Cuturi, M.}, \bibinfo{author}{Butscher, A.}, \bibinfo{author}{Nguyen, A.}, \bibinfo{author}{Du, T.}, \bibinfo{author}{Guibas, L.}, \bibinfo{year}{2015}.
\newblock \bibinfo{title}{Convolutional wasserstein distances: Efficient optimal transportation on geometric domains}.
\newblock \bibinfo{journal}{ACM Transactions on Graphics (ToG)} \bibinfo{volume}{34}, \bibinfo{pages}{1--11}.
\bibitem[{Stankovi{\'c} et~al.(2020)Stankovi{\'c}, Mandic, Dakovi{\'c}, Brajovi{\'c}, Scalzo, Li, Constantinides et~al.}]{stankovic2020data}
\bibinfo{author}{Stankovi{\'c}, L.}, \bibinfo{author}{Mandic, D.}, \bibinfo{author}{Dakovi{\'c}, M.}, \bibinfo{author}{Brajovi{\'c}, M.}, \bibinfo{author}{Scalzo, B.}, \bibinfo{author}{Li, S.}, \bibinfo{author}{Constantinides, A.G.}, et~al., \bibinfo{year}{2020}.
\newblock \bibinfo{title}{Data analytics on graphs part iii: Machine learning on graphs, from graph topology to applications}.
\newblock \bibinfo{journal}{Foundations and Trends{\textregistered} in Machine Learning} \bibinfo{volume}{13}, \bibinfo{pages}{332--530}.
\bibitem[{Su et~al.(2015)Su, Wang, Shi, Zeng, Sun, Luo and Gu}]{su2015optimal}
\bibinfo{author}{Su, Z.}, \bibinfo{author}{Wang, Y.}, \bibinfo{author}{Shi, R.}, \bibinfo{author}{Zeng, W.}, \bibinfo{author}{Sun, J.}, \bibinfo{author}{Luo, F.}, \bibinfo{author}{Gu, X.}, \bibinfo{year}{2015}.
\newblock \bibinfo{title}{Optimal mass transport for shape matching and comparison}.
\newblock \bibinfo{journal}{IEEE transactions on pattern analysis and machine intelligence} \bibinfo{volume}{37}, \bibinfo{pages}{2246--2259}.
\bibitem[{Sun et~al.(2024)Sun, Suresh, Ro, Beirami, Jain and Yu}]{sun2024spectr}
\bibinfo{author}{Sun, Z.}, \bibinfo{author}{Suresh, A.T.}, \bibinfo{author}{Ro, J.H.}, \bibinfo{author}{Beirami, A.}, \bibinfo{author}{Jain, H.}, \bibinfo{author}{Yu, F.}, \bibinfo{year}{2024}.
\newblock \bibinfo{title}{Spectr: Fast speculative decoding via optimal transport}.
\newblock \bibinfo{journal}{Advances in Neural Information Processing Systems} \bibinfo{volume}{36}.
\bibitem[{Ta et~al.(2022)Ta, Tu, Lu and Wang}]{ta2022quantitative}
\bibinfo{author}{Ta, D.}, \bibinfo{author}{Tu, Y.}, \bibinfo{author}{Lu, Z.L.}, \bibinfo{author}{Wang, Y.}, \bibinfo{year}{2022}.
\newblock \bibinfo{title}{Quantitative characterization of the human retinotopic map based on quasiconformal mapping}.
\newblock \bibinfo{journal}{Medical image analysis} \bibinfo{volume}{75}, \bibinfo{pages}{102230}.
\bibitem[{Tong et~al.(2022)Tong, Wolf and Krishnaswamy}]{tong2022fixing}
\bibinfo{author}{Tong, A.}, \bibinfo{author}{Wolf, G.}, \bibinfo{author}{Krishnaswamy, S.}, \bibinfo{year}{2022}.
\newblock \bibinfo{title}{Fixing bias in reconstruction-based anomaly detection with lipschitz discriminators}.
\newblock \bibinfo{journal}{Journal of Signal Processing Systems} \bibinfo{volume}{94}, \bibinfo{pages}{229--243}.
\bibitem[{Tu et~al.(2020)Tu, Ta, Gu, Lu and Wang}]{tu2020diffeomorphic}
\bibinfo{author}{Tu, Y.}, \bibinfo{author}{Ta, D.}, \bibinfo{author}{Gu, X.D.}, \bibinfo{author}{Lu, Z.L.}, \bibinfo{author}{Wang, Y.}, \bibinfo{year}{2020}.
\newblock \bibinfo{title}{Diffeomorphic registration for retinotopic mapping via quasiconformal mapping}, in: \bibinfo{booktitle}{2020 IEEE 17th International Symposium on Biomedical Imaging (ISBI)}, \bibinfo{organization}{IEEE}. pp. \bibinfo{pages}{687--691}.
\bibitem[{Wong et~al.(2023)Wong, Li, Tam, Yuen, Au, Chan, Li and Lui}]{wong2023quasiconformal}
\bibinfo{author}{Wong, M.H.}, \bibinfo{author}{Li, M.}, \bibinfo{author}{Tam, K.M.}, \bibinfo{author}{Yuen, H.M.}, \bibinfo{author}{Au, C.T.}, \bibinfo{author}{Chan, K.C.C.}, \bibinfo{author}{Li, A.M.}, \bibinfo{author}{Lui, L.M.}, \bibinfo{year}{2023}.
\newblock \bibinfo{title}{A quasiconformal-based geometric model for craniofacial analysis and its application}.
\newblock \bibinfo{journal}{Axioms} \bibinfo{volume}{12}, \bibinfo{pages}{393}.
\bibitem[{Yau and Gu(2016)}]{Yau2016ComputationalCG}
\bibinfo{author}{Yau, S.T.}, \bibinfo{author}{Gu, X.}, \bibinfo{year}{2016}.
\newblock \bibinfo{title}{Computational conformal geometry}.
\newblock \URLprefix \url{https://api.semanticscholar.org/CorpusID:117563466}.
\bibitem[{Yu et~al.(2020)Yu, Yang, Roth, Bai, Zhang, Yuille and Xu}]{yu2020c2fnas}
\bibinfo{author}{Yu, Q.}, \bibinfo{author}{Yang, D.}, \bibinfo{author}{Roth, H.}, \bibinfo{author}{Bai, Y.}, \bibinfo{author}{Zhang, Y.}, \bibinfo{author}{Yuille, A.L.}, \bibinfo{author}{Xu, D.}, \bibinfo{year}{2020}.
\newblock \bibinfo{title}{C2fnas: Coarse-to-fine neural architecture search for 3d medical image segmentation}, in: \bibinfo{booktitle}{Proceedings of the IEEE/CVF Conference on Computer Vision and Pattern Recognition}, pp. \bibinfo{pages}{4126--4135}.
\bibitem[{Zeng and Gu(2011)}]{zeng2011registration}
\bibinfo{author}{Zeng, W.}, \bibinfo{author}{Gu, X.D.}, \bibinfo{year}{2011}.
\newblock \bibinfo{title}{Registration for 3d surfaces with large deformations using quasi-conformal curvature flow}, in: \bibinfo{booktitle}{CVPR 2011}, \bibinfo{organization}{IEEE}. pp. \bibinfo{pages}{2457--2464}.
\bibitem[{Zeng et~al.(2012)Zeng, Lui, Luo, Chan, Yau and Gu}]{zeng2012computing}
\bibinfo{author}{Zeng, W.}, \bibinfo{author}{Lui, L.M.}, \bibinfo{author}{Luo, F.}, \bibinfo{author}{Chan, T.F.C.}, \bibinfo{author}{Yau, S.T.}, \bibinfo{author}{Gu, D.X.}, \bibinfo{year}{2012}.
\newblock \bibinfo{title}{Computing quasiconformal maps using an auxiliary metric and discrete curvature flow}.
\newblock \bibinfo{journal}{Numerische Mathematik} \bibinfo{volume}{121}, \bibinfo{pages}{671--703}.
\bibitem[{Zeng et~al.(2009)Zeng, Luo, Yau and Gu}]{zeng2009surface}
\bibinfo{author}{Zeng, W.}, \bibinfo{author}{Luo, F.}, \bibinfo{author}{Yau, S.T.}, \bibinfo{author}{Gu, X.}, \bibinfo{year}{2009}.
\newblock \bibinfo{title}{Surface quasi-conformal mapping by solving beltrami equations}, in: \bibinfo{booktitle}{Mathematics of Surfaces XIII: 13th IMA International Conference York, UK, September 7-9, 2009 Proceedings}, \bibinfo{organization}{Springer}. p. \bibinfo{pages}{391}.
\bibitem[{Zeng et~al.(2014)Zeng, Ming~Lui and Gu}]{zeng2014surface}
\bibinfo{author}{Zeng, W.}, \bibinfo{author}{Ming~Lui, L.}, \bibinfo{author}{Gu, X.}, \bibinfo{year}{2014}.
\newblock \bibinfo{title}{Surface registration by optimization in constrained diffeomorphism space}, in: \bibinfo{booktitle}{Proceedings of the IEEE Conference on Computer Vision and Pattern Recognition}, pp. \bibinfo{pages}{4169--4176}.
\bibitem[{Zeng and Yang(2014)}]{zeng2014surface2}
\bibinfo{author}{Zeng, W.}, \bibinfo{author}{Yang, Y.J.}, \bibinfo{year}{2014}.
\newblock \bibinfo{title}{Surface matching and registration by landmark curve-driven canonical quasiconformal mapping}, in: \bibinfo{booktitle}{Computer Vision--ECCV 2014: 13th European Conference, Zurich, Switzerland, September 6-12, 2014, Proceedings, Part I 13}, \bibinfo{organization}{Springer}. pp. \bibinfo{pages}{710--724}.

\end{thebibliography}



\end{document}